\documentclass{amsart}

\usepackage{graphicx}
\usepackage{color}
\usepackage{framed}
\usepackage{amsthm}
\usepackage{amscd}
\usepackage{amsmath}
\usepackage{amsfonts}
\usepackage{amssymb}
\usepackage{psfrag}
\usepackage{amstext}
\usepackage{hyperref}

\hypersetup{ colorlinks=true, pdftitle={Metric stability for random walks}, pdfsubject={Stability of metric properties of random walks by asymptotically small perturbations and its applications to one-dimensional renormalization theory}, pdfauthor={Carlos Gustavo  Moreira and Daniel Smania},pdfkeywords={random walk, renormalization, stability, rigidity, conjugacy, absolutely continuous, universality, Feigenbaum, period-doubling, Fibonacci, transience, recurrence, strongly quasisymmetric, multifractal, Hausdorff, dimension, wild attractor, group extension, skew-product}}

\begin{document}

\author{Carlos Gustavo  Moreira and Daniel Smania }

\address{Instituto de Matem\'atica Pura e Aplicada (IMPA) \\
   Estrada Dona Castorina, 110 \\
               Jardim Bot\^anico \\ Rio de Janeiro-RJ \\  CEP 22460-320
 \\ Brazil.}
\email{gugu@impa.br}
\urladdr{\href{http://www.impa.br/~gugu}{www.impa.br/$\sim$gugu/}}

\address{Departamento de Matem\'atica \\
   ICMC/USP - S\~ao Carlos \\
               Caixa Postal 668 \\ S\~ao Carlos-SP \\ CEP 13560-970 \\ Brazil.}
\email{smania@icmc.usp.br} \urladdr{\href{http://www.icmc.usp.br/~smania}{www.icmc.usp.br/$\sim$smania/}}

\date{\today}
\title[Metric stability  for random walks]{Metric stability for random walks \\
\tiny{with applications in renormalization theory}}

\begin{abstract} Consider deterministic random walks $F\colon
 I \times \mathbb{Z}\rightarrow I \times \mathbb{Z}$, defined by $F(x,n)=(f(x),\psi(x)+n)$,  where $f$
 is an expanding Markov map on the interval $I$ and
 $\psi\colon I \rightarrow \mathbb{Z}$. We study the
 universality (stability)  of ergodic (for instance, recurrence and transience),
  geometric and multifractal properties in the class of
   perturbations of the type
 $\tilde{F}(x,n)=(f_n(x),\tilde{\psi}(x,n)+n)$ which are topologically conjugate with $F$ and $f_n$ are expanding  Markov  maps exponentially close to $f$ when
 $|n|\rightarrow \infty$. We give applications of these results in the study
 of the regularity of conjugacies between
  (generalized) infinitely renormalizable maps of the interval and the existence of wild attractors for one-dimensional maps.  \end{abstract}


\thanks{We thank the referre for the very useful, detailed and precise comments and suggestions. D.S. was partially supported by CNPq 470957/2006-9, 310964/2006-7,  472316/03-6,  303669/2009-8 and FAPESP 03/03107-9,  2008/02841-4, 2010/08654-1, and  C. G. M. was   partially supported by CNPq. The authors would like to thank the hospitality of ICMC-USP and IMPA}

\maketitle
\newcommand{\co}{\mathbb{C}}
\newcommand{\incl}[1]{i_{U_{#1}-Q_{#1},V_{#1}-P_{#1}}}
\newcommand{\inclu}[1]{i_{V_{#1}-P_{#1},\co-P}}
\newcommand{\func}[3]{#1\colon #2 \rightarrow #3}
\newcommand{\norm}[1]{\left\lVert#1\right\rVert}
\newcommand{\norma}[2]{\left\lVert#1\right\rVert_{#2}}
\newcommand{\hiper}[3]{\left \lVert#1\right\rVert_{#2,#3}}
\newcommand{\hip}[2]{\left \lVert#1\right\rVert_{U_{#2} - Q_{#2},V_{#2} -
P_{#2}}}
\newtheorem{prop}{Proposition}[section]
\newtheorem{lem}[prop]{Lemma}

\newtheorem{rem}[prop]{Remark}

\newtheorem*{mth}{Main Theorem}
\newtheorem{thm}{Theorem}

\newtheorem{cor}[prop]{Corollary}

\setcounter{tocdepth}{1}
\tableofcontents

\section{Introduction}

\subsection{Metric stability for random walks} In the study of a dynamical system, some of the most
important questions concerns  the stability of their dynamical
properties under (most of the) perturbations: how robust are
they?

Here we are mainly interested in the stability of metric
(measure-theoretical) properties of dynamical systems. A
well-known example is given by ($C^2$)  expanding maps on
the circle: this is a class stable under perturbations and all of
them have an absolutely continuous and ergodic invariant probability
satisfying certain decay of correlations estimatives. In
particular, in the measure theoretical sense, most of the orbits
are dense in  the phase space.

Now let us  study a slightly more complicated situation: consider a
$C^2$ Markov almost onto expanding map of the interval $f\colon I \rightarrow
I$ with bounded distortion  and large images (see Section 2 for details) and let $\psi \colon I \rightarrow \mathbb{Z}$ be a function
which is constant in each interval of the Markov partition of $f$.
We can define $F\colon I \times \mathbb{Z} \rightarrow I \times
\mathbb{Z}$ as

$$F(x,n):= (f(x), \psi(x) + n).$$

The second entry of $(x,n)$ will be called its {\bf state}. We
also assume that

\begin{equation}\label{cond1}
\inf \psi > -\infty
\end{equation}
and that $F$ is topologically mixing.

\begin{figure}
\centering \psfrag{f}{$f$}

\psfrag{Fx}[][][0.8]{$F$}
 \psfrag{Fy}[][][0.8]{$F$}
 \psfrag{x}[][][0.8]{$x$}
 \psfrag{y}[][][0.8]{$y$}
\psfrag{fx}[][][0.8]{$f(x)$}
 \psfrag{fy}[][][0.8]{$f(y)$}
 \psfrag{a}[][][0.8]{$i$}
 \psfrag{b}[][][0.8]{$i+1$}
 \psfrag{c}[][][0.8]{$i+2$}
 \psfrag{d}[][][0.8]{$i-1$}
 \psfrag{e}[][][0.8]{$i-2$}
 \psfrag{psix}[][][0.8]{$\psi(x)=-1$}
 \psfrag{psiy}[][][0.8]{$\psi(y)=1$}

\includegraphics[width=1.0\textwidth]{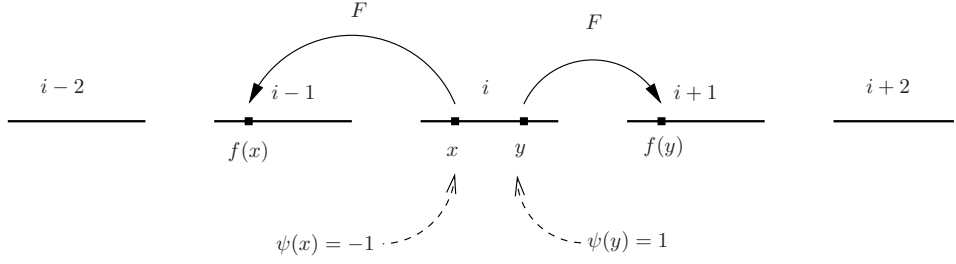}
\caption{A deterministic random walk}
\end{figure}

The map $F$ is  refereed to  in literature in many ways: as a
"skew-product between $f$ and the translation on the group
$\mathbb{Z}$", a "group extension of $f$", or even a
"deterministic random walk generated by $f$", and its metric
behavior is very well studied: for instance, are most the orbits
recurrent? Everything depends on the {\bf mean drift }
$$M = \int \psi d\mu,$$
where $\mu$ is the absolutely continuous invariant probability of
$f$ (the function $\psi$ will be called  {\bf drift function}).
Indeed, note that
$$F^n(x,i)= (\ f^n(x)\ ,\ i + \sum_{k=0}^{n-1}\psi(f^k(x))\ ).$$
By the Birkhoff Ergodic Theorem
$$
\lim_{n \rightarrow \infty}\ \frac{\pi_2(F^n(x,i))- \pi_2(x,i)}{n}
= \lim_{n \rightarrow \infty}\ \frac{1}{n}
\sum_{k=0}^{n-1}\psi(f^k(x))= M.$$ for almost every $x \in I$
(here $\pi_2(x,n):=n$). In particular if $M \neq 0$ then almost every point $(x,i)
\in I \times \mathbb{Z}$ is {\bf transient}: in other words  we
have
$$\lim_{n\rightarrow \infty} |\pi_2(F^n(x,i))|= \infty.$$
 So most of the points are not recurrent.

On the other hand, if $M=0$, most of points are
recurrent (see Guivarc'h \cite{guivarch}): by the Central Limit Theorem for expanding maps (here we need to assume that $\psi$ is not constant and $f \in On$: see Section 2) of the
interval
$$sup_{\epsilon \in \mathbb{R}} \ | \mu( x \in I\colon \frac{\sum_{k=0}^{n-1}\psi(f^k(x))}{\sigma \sqrt{n}} \leq \
\epsilon) - \frac{1}{\sqrt{2\pi}}\int_{-\infty}^{\epsilon}
e^{-\frac{u^2}{2}} \ du | \leq \frac{C}{\sqrt{n}},$$
Given $\delta >0$  we can easily obtain, taking
$\epsilon = n^{-1/4}$ and applying  Borel-Cantelli Lemma, that

\begin{equation}\label{recu1}  \mu(A_{+}):=\mu( x \in I \colon \ \limsup_{n\rightarrow \infty}
 \frac{\sum_{k=0}^{n-1}\psi(f^k(x))}{ \sqrt[2 + \delta]{n}}=\infty)\geq \frac{1}{2},\end{equation}

\begin{equation} \mu(A_{-}):=\mu( x \in I \colon \ \liminf_{n\rightarrow \infty}
 \frac{\sum_{k=0}^{n-1}\psi(f^k(x))}{ \sqrt[2 + \delta]{n}}=-\infty)\geq \frac{1}{2}.\end{equation}

Clearly $A_{+}$ and $A_{-}$ are invariant sets: the ergodicity of
$f$ implies that \begin{equation}\label{recu2} \mu(A_{+}~\cap~A_{-})~=~1.\end{equation} 

By the conditions on $\psi$ in Eq. (\ref{cond1}),  that $f$ is expanding with  distortion control
 and  that $F$ is transitive, we can easily conclude that almost every point in
$I\times \mathbb{Z}$ is a  $F$-recurrent point.

 Note that the  random walk $F$ is a  dynamical system quite similar to
expanding circle maps: $F$ is an expanding map, with good
bounded distortion properties; but the lack of compactness of the
phase space allows the non-existence of an invariant
probability absolutely continuous with respect to the Lebesgue measure on $I\times \mathbb{Z}$. Moreover, in general  the random walk  is not even
recurrent and  the recurrence property lost its stability: given a
recurrent random walk $(f,\psi)$, it is possible to obtain a
transient random walk by  changing a little bit $f$ and $\psi$.

Since the non compactness of the phase  space seems to be the
origin of the lack of stability of  recurrence and transience
properties, a natural question is to ask if such properties are
stable by compact perturbations. The answer is yes. Indeed, as we are going to see
in Theorems \ref{sttr}-\ref{strec}, the transience and recurrence
are preserved even by non-compact perturbations which decreases
fast away from state $0$. For instance,we can choose  perturbations like

$$\tilde{F}(x,n)=(f_n(x),\psi(x)+n),$$
where, for some $\lambda \in [0, 1)$,
\begin{equation}\label{decay}|f_n-f|_{C^3} \leq
\lambda^{|n|}.\end{equation} The notations and conventions are more or less
obvious: we postponed the rigorous definitions to the next
section.

With respect to the stability of transience and recurrence, there
is a previous quite elegant result by R. L. Tweedie \cite{tw}: if
$p_{ij}$ are the transition probabilities of a Markov chain on
$\mathbb{Z}$, then any perturbation $\tilde{p}_{ij}$ so that
$$   (1+\epsilon_i)^{-1}p_{ij} \leq \tilde{p}_{ij} \leq p_{ij}
(1+\epsilon_i), \ j \neq i,$$ and
$$\prod_{i=0}^{\infty} (1+\epsilon_i) < \infty$$
preserves the recurrence or transience of the original Markov
chain. But Tweedie argument does not seem to work in our setting.
Our result coincides with Tweedie result in the very special case
where $f$ and $f_n$ are linear Markov maps and $\epsilon_i \sim
C\lambda^{|i|}$.

In the transient case we can tell a little more: there will be a
conjugacy between the original random walk $f$ and its
perturbation which is a martingale strongly quasisymmetric
 map (for short, mSQS-map) with respect to certain dynamically defined set of partitions.
Unlike the usual class of one-dimensional quasisymmetric
functions, which does not share many of most interesting
properties of  higher dimensional quasisymmetric maps,  the
one-dimensional mSQS-maps are
 much closer to their high-dimensional cousins, as quasiconformal maps in dimension $2$. For  instance,  they are absolutely
continuous.

We also study the behavior of the Hausdorff dimension of dynamically defined sets: Denote by $\Omega_+(F)$ the set of points which have non-negative states along the positive orbit by $F$. We prove that $\Omega_+(F)$ has Hausdorff dimension strictly smaller than one if and only if $\Omega_+(\tilde{F})$ has dimension less than one for all perturbation satisfying Eq. (\ref{decay}). Furthermore we give a variational characterization for the Hausdorff dimension $HD(\Omega_+(F))$ as the minimum of $HD(\Omega_+(\tilde{F}))$, where $\tilde{F}$ runs on the set of such perturbations. For these results we study of the stability of the multifractal spectrum of the random walk $F$ under those perturbations.

\subsection{Applications to (generalized) renormalization theory}An unimodal map is a map with an unique  critical point. Under reasonable conditions (real-analytic maps with negative Schwarzian derivative and  non-flat critical point) two non renormalizable unimodal maps with the same topological entropy are indeed topologically conjugated. A key question in one-dimensional dynamics is about the regularity of the conjugacy: is it H\"older? Is it absolutely continuous? Since Dennis Sullivan work in the 80's the quasisymmetry of the conjugacy became a very useful tool to obtain deep results in one-dimensional dynamics.  Lyubich proved that under the reasonable condition above the conjugacy between two non renormalizable unimodal maps is quasisymmetric. Later on, the density of the hyperbolic maps in the real quadratic family was proved verifying  the quasisymmetry of the conjugacies for all combinatorics, including infinitely renormalizable ones.

Note that quasisymmetric maps are not, in general, absolutely continuous: they do not even preserve (in general) sets of Hausdorff dimension one. Are the conjugacy between unimodal maps absolutely continuous? The answer is no: M.~Martens and W.~de~Melo \cite{mm} proved that under the reasonable conditions above an absolutely continuous conjugacy is actually $C^\infty$, provided the unimodal maps
\begin{itemize}
\item[ ] \
\item[{\it i.}] {\it do not have  a periodic attractor,}\\
\item[{\it ii.}]{\it  are not infinitely renormalizable, }\\
\item[{\it iii.}] {\it do not have a wild attractor (the topological and measure-theoretical attractor must coincide).} \\
\end{itemize}

Since we can change the eigenvalues of the periodic points of maps preserving its topological class, and the eigenvalues are preserved by $C^1$ conjugacies, we conclude that in general a conjugacy between unimodal maps is not absolutely continuous.

Condition i. is clearly necessary. This work (Theorem \ref{apl1}) shows that the Condition ii.  is necessary proving that the conjugacy between two arbitrary Feigenbaum unimodal maps with same critical order is {\em always }  absolutely continuous . Actually the conjugacy is martingale strongly quasisymmetric with respect to a set of dynamically defined partitions.

Condition iii. is never violated when the critical point is quadratic. But for certain topological classes of unimodal maps wild attractor appears when the order of the critical point increases: Fibonacci maps are the simplest kind of such maps \cite{bkns}\cite{bruin}. We are going to prove (Theorem \ref{apl2}) that a Fibonacci map with even order has a wild attractor if and only if all  Fibonacci maps with the same even order are conjugated to each other by an absolutely continuous mapping (in particular all these Fibonacci maps have a wild attractor). So Condition iii. is necessary.

To show that conditions i. and ii. are necessary, the (generalized) renormalization theory for unimodal maps and the study of perturbations of transient and recurrent random walks are going to be crucial. Feigenbaum and Fibonacci unimodal maps  admit  induced maps which are  essentially  perturbations of  deterministic random walks (Section \ref{applications}). In the Fibonacci case the transience of this random walk is equivalent to the existence of a wild attractor. Random walks associated to a Feigenbaum map will always be  transient.

For both Feigenbaum and Fibonacci maps there are infinitely many periodic points (indeed in the Fibonacci case the periodic points are also dense in the maximal invariant set). It is well known that the conjugacy between critical circle maps with same irrational rotation number and satisfying certain Diophantine condition is absolutely continuous, but we think that these are the first interesting examples of a similar phenomena for maps with many periodic points.

\section{Expanding Markov maps, random walks and its perturbations}

In this article we will deal with maps
$$F\colon I \times \mathbb{Z} \rightarrow I \times \mathbb{Z}$$
which are piecewise $C^2$ diffeomorphisms, which means that there is a partition $\mathcal{P}^0$  of $I \times \mathbb{Z}$ so that each element $J \in \mathcal{P}^0$ is an open interval where $F|_{\overline{J}}$ is a $C^2$ diffeomorphism. Denote $I_n=I\times \{ n \}$. Denote by $m$ the Lebesgue measure in the in $I\times\mathbb{Z}$, that is, if $A\subset I\times\mathbb{Z}$ is a Borelian set then
$$m(A)=\sum_n m_I(\pi(A\cap I_n)),$$
where $m_I$ is the Lebesgue measure in the interval $I$ and $\pi(n,x)=x$.

If $A_J$ denotes the unique affine transformation which maps the interval $J$ to $[0,1]$ and preserves orientation,  then  define, for each  $J \in \mathcal{P}^0$,  $$\tau_J^F:= A_{J}\circ F^{-1} \circ A_{F(J)}^{-1}.$$

 Throughout this article we will assume that $F$  satisfies some of the following properties:

\begin{itemize}
\item[\ ] \ \\

\item {\bf Markovian (Mk)}:  For each $J \in \mathcal{P}^0$,
$F(J)$ is a connected union of elements in $\mathcal{P}^0$. In particular we can write $F(x,n)=(f_n(x), n + \psi(x,n))$, where $f_n\colon I \rightarrow I$ is a piecewise $C^2$ diffeomorphism relative to the partition $\mathcal{P}^0_n:=\{J \in \mathcal{P}^0\colon \ J \subset I_n   \}$ and  $\psi\colon I \times \mathbb{Z} \rightarrow \mathbb{Z}$, called the {\bf drift function}, is constant on each element of $\mathcal{P}^0$.\\

\item {\bf Lower Bounded Drift (LBD)} $F$ is Markovian and $\min \psi > -\infty$.\\

\item {\bf Large Image (LI)}: $F$ is Markovian and there exists $\delta > 0$ so that for each $J \in \mathcal{P}^0$ we have $|F(J)|\geq \delta$.\\

\item {\bf Onto (On)}: $F$ is Markovian and for each $J \in \mathcal{P}^0$ we have $F(J)=I^n$, for some $n \in \mathbb{Z}$.\\

\item {\bf Bounded Distortion (BD)}: There exists $C > 0$ so that every $J \in \mathcal{P}^0_n$ and map $\tau_J$ is a $C^2$ function satisfying
$$\sup_{J} \Big| \frac{D^2\tau_J}{(D\tau_J)^2} \Big| \leq C.$$ \\

\item {\bf Strong Bounded Distortion (sBD)}: There exists $C > 0$ so that every $J \in \mathcal{P}^0_n$ and map $\tau_J$ is a $C^2$ function satisfying
$$\sup_{J} \Big| \frac{D^2\tau_J}{(D\tau_J)^2} \Big| \leq C|J|.$$  \\

\item {\bf Expansivity (Ex)}: If $J \in \mathcal{P}^0_n:=\{ J \in \mathcal{P}^0\colon \ J \subset I_n\}$, denote $\phi_J:=f_n^{-1}|_{f_n(J)}$. Then  either $\phi_J$ can be extended to a function in a $\delta$-neighborhood of $J$ so that $$ S\phi_J > 0,$$ where $S\phi_J$ denotes the Schwarzian derivative of $\phi_J$,  or there exists $\theta \in (0,1)$ so that for all $n$ and 
$J \in \mathcal{P}^0_n$ we have $$|\phi'_J| < \theta$$
on $f_n(J)$. \\

\item {\bf Regularity a (Ra)}: There exists $N \in \mathbb{N}$, $\delta > 0$ and $C > 0$ with the following properties:  the intervals in $\mathcal{P}^0_n$ are positioned in $I_n=[a,b]$ in such way that the complement of  $$\bigcup_{J \in \mathcal{P}^0_n}\  int \ J$$ contains at most $N$  accumulation points $$c_1^n < c_2^n < \dots < c^{n}_{i_n},$$ with $i_n\leq N$, which are in the interior of $I_n$.  Furthermore $|c^n_{i+1}-c^n_i|\geq \delta$ and $|a-c_1^n|$, $|b-c_{i_n}^n|\geq \delta$.    Moreover, given $P$ and $Q \in \mathcal{P}^0_n$ so that $\overline{P} \cap \overline{Q} \neq \phi$ then
\begin{equation} \label{prop_rb}\frac{1}{C} \leq \frac{|P|}{|Q|} \leq C.\end{equation}
\item {\bf Regularity b (Rb)}:  Assume $Ra$. There exists $C > 0$, $\lambda \in (0,1)$, $\delta > 0$ so that for each $1< i<i_n$ we can find a point $$d_i^n \in (c^n_i,c^n_{i+1}),$$
which does not belong to any $P \in \mathcal{P}^0_n$, and $$\min \{|c^n_{i+1}-d_i^n|, |d_i^n- c^n_{i}| \} \geq \delta$$ with the following property: If $J$ is a connected component of  $$I_n \setminus \{d_i^n, c^n_j\}_{i,j}$$
then  we can enumerate the set $$\{P\}_{P \in \mathcal{P}^0_n,\  P \subset J}=\{J_i\}_{i \in  \mathbb{N}}$$ in such way that  $\partial J_{i} \cap \partial J_{i+1} \neq \phi$ for each $i$  and   $$\frac{|J_{i+j}|}{|J_{i}|}\leq  C\lambda^j$$ for $i \geq  0$, $j > 0$.\\

\item {\bf Good Drift (GD)}: , if $\psi$ is the  drift function of the random walk  then for each $n\in \mathbb{Z}$ there exists $x$ such that $\psi(x,n) > 0$. Moreover there exists $\gamma \in (0,1)$ and $C > 0$ so that for every $k\geq 0$
$$m(\{(x,n) \ s.t. \ \psi(x,n)\geq k   \}) \leq   C\gamma^k.$$ \\

\item {\bf Transitive (T)}: $F$ has a dense orbit.
\end{itemize}

For convenience of the notation if for instance  $F$ is Markovian and it has Bounded Distortion, we will write $F \in Mk+BD$.

A {\bf deterministic random walk} (or simply random walk) is a map $$F \in Mk+LBD+LI+Ex+BD+GD.$$ It is generated by the pair $(\{ f_n\},\psi)$ if
$$F(x,n):= (f_n(x),\psi(x,n)+n).$$

When $f_n=f \in Mk$ and $\psi(x,n)=\psi(x)$, we say that $F$ is the
{\bf  homogeneous deterministic random walk} generated by the pair $(f,\psi)$. There is a large
literature about such random walks. We will sometimes assume  the following property:

\begin{itemize}
\item[\ ] \ \\
\item {\bf Almost Onto (aO)}: For every $i, j \in \Lambda$ there exists a finite sequence $i=i_0,i_1,i_2,\dots,i_{n-1},i_n=j \in \Lambda$ so that  $$f(I_{i_k})\cap f(I_{i_{k+1}})\neq \emptyset$$ for each $k < j$. \\
\end{itemize}

Denote $\pi(x,n):= \pi_2(x,n):=n$. A  random walk is called {\bf transient} if for almost every $(x,n) \in I \times \mathbb{Z}$
$$\lim_{k\rightarrow \infty} |\pi_2(F^k(x,n))| = \infty,$$
and it is {\bf recurrent} if for almost every $(x,n) \in I \times \mathbb{Z}$
$$\# \{ k \colon \pi_2(F^k(x,n))=n \} = \infty.$$
Making use of usual bounded distortion tricks it is easy to show that every $F \in Mk+LI+Ex+BD+T$ is either recurrent or transient.

A (topological) {\bf perturbation} of a random walk is a random walk $\tilde{F}$, generated by a pair $(\{\tilde{f}_i\},\tilde{\psi})$, so that $H\circ F=  \tilde{F}\circ H$ for some homeomorphism $$H\colon I \times \mathbb{Z} \rightarrow I \times \mathbb{Z}$$
which preserves states: $\pi_2(H(x,i))=i$.

Define $\mathcal{P}^n(F) := \vee_{i=0}^{n-1} F^{-i}\mathcal{P}^0(F)$. If $F$ and $\tilde{F}$ are random walks and $H$ is a topological conjugacy that preserves states between $F$ and $\tilde{F}$, then for each interval $L$ such that $L\subset J \in \mathcal{P}^{n-1}(F)$, define

$$dist_n(L):= \sup_{y \in L} \Big|\ln \frac{D\tilde{F}^n(H(y))}{DF^n(y)} \Big|,$$
Similarly,  define
$$dist_n(x):= \Big|\ln \frac{D\tilde{F}^n(H(x))}{DF^n(x)} \Big|$$
and
$$dist_\infty(x):= \sup_n dist_n(x).$$

Another kind of random walk which will have a central role in our
results are those which are {\bf asymptotically small} perturbations: these are  perturbations $(\{\tilde{f}_i\},\tilde{\psi})$ of a deterministic  random walk $(\{f_i\},\psi)$ such that there exists $\lambda \in (0,1)$  and $C >0$ satisfying either
\begin{equation}\label{asymp} |\log \frac{D\tilde{F}(H(p))}{DF(p)}| \leq C\lambda^{|\pi_2(p)|},\end{equation}
if $\psi$ is bounded, or
\begin{equation}\label{asymp2} |\log \frac{D\tilde{F}(H(p))}{DF(p)}| \leq C\lambda^{\pi_2(p)},\end{equation}
for $\pi_2(p)\geq 0$ and $D\tilde{F}(H(p))=DF(p)$ otherwise, if $\psi$ has only a lower bound.

It is easy to see that properties $Ra$, $Rb$ and $GD$ are invariant by asymptotically small perturbations (if we allow to change the constants described in these properties).

Let $F=(\{f_i\},\psi)$ be a random walk, where $\psi$ is Lebesgue integrable on compact subsets of $I \times \mathbb{Z}$.  We say that $F$ is {\bf strongly transient} if $K > 0$ and
$$\mathbb{E}(\psi\circ F^n | \mathcal{P}^{n-1}(F)) > K$$
for every $n\geq 1$.  We will also say that  $F$ is $K$-strongly transient. Here we are considering   conditional expectations relative to the  Lebesgue  measure. As the notation suggest, every strongly transient random walk is transient. Moreover we have the following large deviations result:

\begin{prop}\label{largedeviationssthomo}\label{bru2} Let  $F=(f,\psi)\in On+sBD+Ra+Rb$ be a homogeneous random walk with positive mean drift. Let $K := \int \psi \ dm > 0$.  Then $F$ is transient and  for every small $\epsilon >  0$  there exist $\lambda \in [0,1)$  and $C > 0$  so that for each $P \in \mathcal{P}^0$ we have
$$m( p \in P \colon \ \pi_2(F^n(p))-\pi_2(p)  < (K-\epsilon) n )\leq C\lambda^n |P|.$$\end{prop}

\begin{prop}\label{largedeviationsst}\label{bru} Every $K$-strongly transient random walk $F\in Ra+Rb$ is transient. Furthermore for every small $\epsilon > 0$ there exist $\lambda \in [0,1)$  and $C > 0$ so that for each $P \in \mathcal{P}^0$ we have
$$m( p \in P \colon \ \pi_2(F^n(p))-\pi_2(p)  < (K -  \epsilon) n )\leq C\lambda^n |P|.$$
\end{prop}

We will postpone the proof of Propositions \ref{largedeviationsst} and \ref{largedeviationssthomo} to Section \ref{sectiontransience}.

\begin{rem}{\rm  By the Birkhoff Ergodic Theorem it is easy to see that  a sufficiently high iteration of a homogeneous random walk with positive mean drift is strongly transient (see the proof of Proposition \ref{homstr} for details). }
\end{rem}

\section{Statements of results}

\subsection{Stability of transience}

\begin{thm}[Stability of Transience I]\label{sttr} Assume that
the random walk $F$ defined by the pair $(\{f_i\},\psi)$ is
strongly transient. Then every asymptotically small perturbation $G$ of $F$ is also transient.
Indeed there is a topological conjugacy between $F$ and $G$  which
is an absolutely continuous map and  preserves the states.
\end{thm}

We have a similar theorem for all transient homogeneous random
walks:
\begin{thm}[Stability of Transience II]\label{sttrho}\label{abscont}  Suppose that the homogeneous random walk $F$ defined by
the pair $(f,\psi)$ has positive mean drift. Then every asymptotically small perturbation of $F$ is topologically conjugated to $F$ by
an absolutely continuous map which preserves the states.
\end{thm}

We can be more precise regarding the regularity of the conjugacy
if the drift is non-negative:

Let $\mathcal{A}_0,\  \mathcal{A}_1, \  \cdots
, \mathcal{A}_n, \  \mathcal{A}_{n+1}, \cdots $ be a succession of partitions by intervals of $I \times \mathbb{Z}$, such that $\mathcal{A}_{n+1}$ refines $\mathcal{A}_{n}$ and  whose union generates the Borelian
algebra of $\sqcup_n I_n$. We say that $h\colon \sqcup_n I_n \rightarrow \sqcup_n I_n$ is a {\bf martingale strongly quasisymmetric (mSQS)} map with respect to the {\bf stochastic basis } $\cup_n \mathcal{A}_n$  if there exist $C > 0$ and $\alpha \in (0,1]$ so that
$$   \frac{m(h(B))}{|h(J)|} \leq C \left( \frac{m(B)}{|J|}  \right )^\alpha$$
for all Borelian $B \subset J \in \cup_n \mathcal{A}_n$, and the same inequality holds replacing $h$ by $h^{-1}$ and $\cup_n \mathcal{A}_n$ by $\cup_n h(\mathcal{A}_n)$.

\begin{thm}[Strongly quasisymmetric rigidity]\label{sqr} Let $F$ be either a strongly transient random walk  or a
transient homogeneous random walk with positive mean drift. Moreover assume in both cases that $\psi\geq 0$.  Then every asymptotically small perturbation $G$ of
$F$ is topologically conjugated to $F$ by
an absolutely continuous map $h$ which preserves the states.
Furthermore $h$ on $\cup_{i\geq 0} I_i$ is a martingale strongly quasisymmetric mapping with respect
to the stochastic basis $\cup_i \mathcal{P}^i.$
\end{thm}

\subsection{Stability of recurrence}

In the recurrent case, we are going to restrict ourselves to the
stability of the metric properties of homogeneous random walks under asymptotically
small perturbations: it is easy to see that the recurrence
is not stable by perturbations which are not asymptotically small. Nevertheless

\begin{thm}[Stability of Recurrence]\label{strec} Suppose that  $F \in On+T$ is a
recurrent homogeneous random walk generated by the pair
$(f,\psi)$. Then every asymptotically small perturbation of $F$ is also recurrent.
\end{thm}

If $p$ is a periodic point with prime period $n$ then $DF^n(p)$ is called the spectrum of the periodic point $p$.
Note that we can not expect, as in the transient case,  an
absolutely continuous conjugacy which preserves states between $F$
and $G$,  once asymptotic small perturbations do not
preserve (in general) the spectrum of the periodic points and:

\begin{prop}[Rigidity]\label{rigidity} Suppose that the random walk $F \in On$ generated by
a pair $(\{f_i  \}_i,\psi)$ is recurrent.
If there is an absolutely continuous conjugacy which preserves
states $H$ between $F$ and a random walk $G$, then $H$ is $C^1$ in
each state. In particular the spectrum of the corresponding
periodic points of $F$ and $G$ are the same.
\end{prop}

The reader should compare this result with similar results by Shub and Sullivan \cite{ss} for expanding maps on the circle and de Melo and Martens \cite{mm} for unimodal maps.

\subsection{Stability of the multifractal spectrum}

Let $F$ be a random walk and denote

$$\Omega_+(F):=\{p \colon \ \pi_2(F^jp)\geq 0, \ for \ j\geq 0    \},$$

$$\Omega_+^k(F):= \{(x,k) \colon \ \pi_2(F^j(x,k))\geq 0, \ for \ j\geq 0    \}$$

and
$$\Omega_{+\beta}^k(F):= \{  (x,k)  \in  \Omega_+^k\ s.t\ \ \underline{\lim}_{\ n} \ \frac{\pi_2(F^n(x,k))}{n}\geq\beta\}$$

\begin{thm}\label{multi} Let $F\in Ra + Rb+ On$ be a random walk. Then,  for all $k \in \mathbb{Z}$ and $\beta >  0$ the Hausdorff dimension $HD(\Omega_{+\beta}^k)$ is
invariant
 by asymptotically small perturbations.
\end{thm}

We will need

\begin{prop}\label{ms}\label{inv} Let $F \in Ra+Rb+On$ be a homogeneous random walk.  Then
$$HD(\Omega_+^k(F))= \lim_{\beta\rightarrow 0^+} HD(\Omega_{+\beta}^k(F)).$$
\end{prop}
and as a consequence of Theorem \ref{multi} and Proposition \ref{ms}:
\begin{thm}\label{omega} Let $F\in Ra+Rb+On$ be a homogeneous random walk. If $G$ is an
asymptotically small perturbation of $F$ then
\begin{equation}\label{ineq} HD(\Omega_{+}^k(G)) \geq HD(\Omega_{+}^k(F)).\end{equation}
\end{thm}

We can not replace the inequality in Eq. (\ref{ineq}) by an
equality. Indeed, even if $HD(\Omega_{+}^k (F)) < 1$, we have that
$sup \ HD(\Omega_{+}^k(G)) =1$, where the supremum is taken on all
asymptotically small perturbations $G$ of $F$. Nevertheless:

\begin{thm}\label{menor} Let $F \in Ra+Rb+On+T$ be the homogeneous random walk generated by the  pair $(f,\psi)$.
Consider  $M=\int \psi d\mu$, where $\mu$ is the unique
absolutely continuous invariant measure of $f$.
\begin{itemize}
\item[-] If $M >0$ then  for all
asymptotically small perturbations $G$ of $F$ we have $m
(\Omega_{+}(G))> 0$.\\
\item[-] If $M =0$ then  for all
asymptotically small perturbations $G$ of $F$ we have $HD
(\Omega_{+}(G))=1$ but $m
(\Omega_{+}(G))= 0$.\\
\item[-] If $M < 0$ then  for all
asymptotically small perturbations $G$ of $F$ we have $HD
(\Omega_{+}(G))< 1$.
\end{itemize}
\end{thm}

\begin{rem}{\rm  Since the authors are more familiar with deterministic rather than stochastic terminology, we stated and proved  the results in this work  for determinist random walks. However  we believe that the above results could be easily translated to the theory  of  chains with  complete connections (g-measures, chains of infinite order) and one-sided shifts on an infinite alphabet. }
\end{rem}

\subsection{Applications to renormalization
theory of  one-dimensional maps}

\begin{thm}\label{apl1} Let $f$ and $g$ be unimodal maps which are infinitely
renormalizable with the same bounded combinatorial type and even
critical order. Then the continuous conjugacy $h$ between $f$ and
$g$ is a strongly quasisymmetric mapping with respect to a certain
stochastic basis of intervals  $\mathcal{P}$.\end{thm}

The set of intervals $\mathcal{P}$ is defined using a map induced
 by $f$. See the details in Section \ref{apl}.
 
 \begin{rem} {\rm D. Sullivan \cite{ds}\cite{ms}  show  that on the assumptions of Theorem \ref{apl1} the conjugacy $h$ is a quasisymmetric map. However it is known that  quasymmetric maps on the real line are not in general absolutely continuous maps. }\end{rem}

Let $\mathcal{F}_d$ be  the class of
 analytic maps with negative Schwarzian  derivative which  are infinitely renormalizable in the Fibonacci sense with even critical order $d$ (see
 Section \ref{aplf} for definitions).  If $f$ is a Fibonacci map, denote by $J_{\mathbb{R}}(f)$ the maximal invariant set of $f$.  Let $\mathcal{F}_d^{uni}$ be the class of Fibonacci {\it unimodal}  maps with negative Schwarzian derivative.

\begin{thm}[Metric Universality]\label{juliathm} For each even critical order $d$,  $d\geq 4$, one of the
following statements holds:
\begin{itemize}
\item $HD(J_{\mathbb{R}}(f)) < 1$, for all $f \in \mathcal{F}_d$.

 \item $HD(J_{\mathbb{R}}(f))= 1$ and $m(J_{\mathbb{R}})=0$ for all $f \in
\mathcal{F}_d$.

\item $HD(J_{\mathbb{R}}(f))= 1$ and $f$ has a wild attractor (in
particular, $m(J_{\mathbb{R}}(f))> 0$) for all $f \in \mathcal{F}_d$
\end{itemize}
\end{thm}

\begin{thm}[Measurable Deep Point]\label{deep} Let $f \in \mathcal{F}_d$, where $d\geq4$ is an even integer, and assume that $0$ is its critical point. If $J_{\mathbb{R}}(f)$ has positive Lebesgue measure then there exists $\alpha > 0$ and $C > 0$ so that
$$m(x \in (-\delta,\delta)\colon \ x \not\in J_{\mathbb{R}}(f))\leq C\delta^{1+\alpha}.$$

\end{thm}
\begin{rem} {\rm Indeed $\alpha$ can be taken depending only on $d$. }\end{rem}

\begin{thm}\label{apl2}For each even critical order $d$, $d\geq 4$,    the following statements are
equivalent:
\begin{enumerate}
\item There exists $f \in \mathcal{F}_d$ such that $m(J_{\mathbb{R}}(F)) > 0$.
\item There exists $f \in \mathcal{F}_d$ with a wild attractor.
\item There exist maps $f,g \in \mathcal{F}_d^{uni}$ which are
conjugated by a continuous absolutely continuous maps $h$, but $f$ has a periodic point $p$ whose eigenvalue is different from the eigenvalue of  the periodic point $h(p)$ of $g$. \item All maps in
$\mathcal{F}_d$ have wild attractors. \item All maps in
$\mathcal{F}_d^{uni}$ can be conjugated with each other by an
absolutely continuous conjugacy.
\end{enumerate}
\end{thm}

\section{Preliminaries}

\subsection{Probabilistic tools.}
We are going to collect here a handful of probabilistic tools which are going to be useful along the article. A good reference for these results is \cite{broise}.

Most of the probabilistic results in dynamical systems (large deviation, central limit theorem) assumes the  observables is quite regular: usual regularity assumptions are either Holder continuity or bounded variation. Fix $(f,\psi) \in On+ sBD+Ra+Rb+GD$. Then $f$ has a unique absolutely continuous invariant probability $\mu$. Moreover this invariant measure is ergodic (see \cite[page 29]{broise}). We are interested in $\mathcal{P}^0$-measurable observables with integer values which do not have such regularity. Fortunately   this is almost true: Denote by  $ \mathcal{O}(f)$ the class of $\mathcal{P}^0$-measurable functions $\phi\colon I \rightarrow \mathbb{Z}$ so that
\begin{itemize}
\item[ ] \
\item[-] $\phi \in L^2(\mu)$,\\

\item[-] If $P$ denotes the Perron-Frobenius-Ruelle operator of $f$, then $P\phi$ has bounded variation.\\

\end{itemize}
Then  $\psi \in \mathcal{O}(f)$. Up to simple modifications in the proofs in \cite{broise}, we have

\begin{prop}[Large Deviations Theorem \cite{broise}]\label{ldt} Suppose $(f,\psi) \in On+ sBD+Ra+Rb+GD$. For every $\psi \in \mathcal{O}(f)$ and $\epsilon > 0$ there exists $\gamma \in (0,1)$ and $C \geq 0$ so that
$$\mu(\{x \in I \colon |\frac{1}{n}  \sum_{i=0}^{n-1}\psi(f^i(x)) - \int \psi d\mu|\geq \epsilon   \}) \leq C \gamma^n$$
\end{prop}

\begin{prop}[Proposition 6.1 of \cite{broise}] Suppose $(f,\psi) \in On+ sBD+Ra+Rb+GD$. For every $\psi \in \mathcal{O}(f)$ the limit
$$\sigma^2 := \lim_{n\rightarrow \infty} \int \left(   \frac{1}{\sqrt{n}} \sum_{k=0}^{n-1}\psi(f^k(x)) \right) ^2  d\mu$$ exists. Furthermore $\sigma^2=0$ if and only if there exists a function $\alpha \in L^2(\mu)$ so that $$\psi = \alpha\circ f - \alpha. $$
\end{prop}
and
\begin{prop}[Central Limit Theorem: Theorem 8.1 in \cite{broise}]\label{clt} Suppose $(f,\psi) \in On+Mk+ sBD+Ra+Rb+GD$.  For every $\psi \in \mathcal{O}(f)$ so that $\sigma^2\neq 0$ we have  \begin{equation}\label{cltt} sup_{\epsilon \in \mathbb{R}} \ | \mu( x \in I\colon \frac{\sum_{k=0}^{n-1}\psi(f^k(x))}{\sigma \sqrt{n}} \leq \
\epsilon) - \frac{1}{\sqrt{2\pi}}\int_{-\infty}^{\epsilon}
e^{-\frac{u^2}{2}} \ du | \leq \frac{C}{\sqrt{n}},\end{equation}
\end{prop}

Indeed we are going to see that the assumption $\sigma^2 \neq 0$ is very weak: to this end we  need the following result:

\begin{prop}[Theorem 3.1 in \cite{ad}]\label{cocycle}  Let $f\colon \cup_i I_i \rightarrow I$ be a map in Mk + BD + Ex + Ra + Rb. Let $\psi \colon \cup_i I_i \rightarrow \mathbb{S}^1$ be a $\mathcal{P}_0$-measurable function. If $$\psi = \frac{\alpha\circ f}{\alpha},$$ where $\alpha$ is measurable, then $\alpha$ is $\mathcal{P}^\star$-measurable, where $\mathcal{P}^\star$ is the finest partition of $I$ so that  $f(I_i)$ is included in an atom of $\mathcal{P}^\star$ for each $i \in \Lambda$.
\end{prop}

\begin{prop}\label{cco} Let $\psi \colon \cup_i I_i \rightarrow \mathbb{Z}$ be a $\mathcal{P}^0$-measurable function. If $\psi = \alpha\circ f - \alpha$, where $\alpha$ is measurable, then $\alpha$ is constant  on $f(I_i)$, for each $i \in \Lambda$.
\end{prop}
\begin{proof} Note that we can assume that $\alpha(x) \in \mathbb{Z}$, for every $x$. Indeed, the relation $\psi = \alpha\circ f - \alpha$ implies that the function $\beta(x)=\alpha(x) \mod 1$ is $f$-invariant, so we can replace $\alpha$ by $\alpha -\beta$, if necessary.  Fix an irrational number $\gamma$. Then
$$e^{2\pi \gamma \psi(x)i} = \frac{e^{2\pi \gamma \alpha(f(x))i}}{e^{2\pi \gamma \alpha(x)i}},$$
so by Proposition \ref{cocycle} we have that $e^{2\pi \gamma \alpha(x)i}$ is a $\mathcal{P}^\star$-measurable function. Since $j \in \mathbb{Z} \rightarrow e^{2\pi \gamma j i} \in \mathbb{S}^1$ is one-to-one, we get that $\alpha$ is $\mathcal{P}^\star$-measurable.
\end{proof}

A Markov map $f$ is almost onto if and only if $\mathcal{P}_0^\star= \{ I\}$, so

\begin{cor}On the conditions of Proposition \ref{cco}, if $f$ is almost onto then $\alpha$ is constant.\end{cor}

\begin{cor} \label{sigma_nz}For every nonconstant $\psi \in \mathcal{O}(f)$ we have that $\sigma^2 \neq 0$. In particular the Central Limit Theorem as given in Eq. (\ref{cltt}) holds for every non-constant $\psi$.\end{cor}

Let  $\mathcal{A}_0 \subset \mathcal{A}_1 \subset \mathcal{A}_2 \subset \dots $ be an increasing sequence of $\sigma$-subalgebras of a probability  space $(\Omega,\mathcal{A},\mu)$. A {\bf martingale difference sequence } is a sequence of functions $\psi_n\colon \Omega \rightarrow \mathbb{R}$, where $\psi_n$ is $\mathcal{A}_n$-measurable for $n\geq 1$, so that
$$\mathbb{E}(\psi_n | \mathcal{A}_{n-1}) =0$$
for every $n$. Here $\mathbb{E}(\psi | \mathcal{B})$ denotes de conditional expectation of $\psi$ relative to the sub-algebra $\mathcal{B}$. When $\mathcal{B}$ is generated by atoms $\{J_i  \}_i$ then $\mathbb{E}(\psi | \mathcal{B})$ is the function defined as
$$\mathbb{E}(\psi | \mathcal{B})(x)= \frac{1}{\mu(J_i)} \int_{J_i} \psi \ d\mu$$
for every $x \in J_i$.

The following Proposition is the classic Azuma-Hoeffding inequality: see, for instance Exercise E14.2 in \cite{williams}:

\begin{prop}[Azuma-Hoeffding inequality]\label{azuma} Let $\psi_n$ be a martingale difference sequence  and furthermore assume that
$$|| \psi_i ||_{\infty} = c_i < \infty.$$
Define
$$\psi:= \sum_{i=1}^{n} \psi_i.$$
Then

$$\mu(x \in \Omega\colon \ |\psi - \mathbb{E}(\psi)| > t ) \leq 2\exp ({-\frac{t^2}{2\sum_{i=1}^{n} c_i^2}}).$$

\end{prop}

\subsection{How to construct asymptotically small perturbations.}

As we will see in the next Proposition, it is easy to construct asymptotically small  perturbations of a random walk:

\begin{prop}\label{identity} \label{how}Let $F$ and $G$ be random walks satisfying the properties $LI$, $Ex$, $sBD$, $Ra$ and $Rb$, where $G$ is a topological perturbation of $F$. Assume that there exist $C >0$ and $\lambda \in (0,1)$ with the following properties: if $I^n_j$ is as in properties $Ra$ and $Rb$, then
\begin{itemize}
\item[ ] \
\item[i.] For every $I^n_j \in \mathcal{P}^0_n$ we have
$$| log \frac{|I^n_{j+1}|}{|I^n_j|}\frac{|H(I^n_j)|}{|H(I^n_{j+1})|}|  \leq C\lambda^{|n|+|j|}.$$\\

\item[ii.] For every $J \in \mathcal{P}^0_n$ we have
$$ |\tau_{J}^F - \tau_{H(J)}^G|_{C^2} \leq C\lambda^{|n|}.$$\\

\item[iii.] If $I_i^n=[a_i^n,b_i^n]$ then
$$\max_i \max \{|a_i^n - H(a_i^n)|, |b_i^n - H(b_i^n)|  \} \leq C\lambda^{|n|}.$$\\

\item[iv.] Either $\psi$ is a bounded funtion or $\psi$ has a lower bound and $F = G$ on $\cup_{n<0}I_n$.

\end{itemize}
Then $G$ is an asymptotically small perturbation of $F$. Furthermore there exist $\beta \in [0,1)$ and $C > 0$ so that
$$|H(p)-p|\leq C\beta^{|\pi_2(p)|}.$$
\end{prop}
\begin{proof}
We will assume that $\psi$ is bounded: the other case is analogous.  Consider $(x,n) \in I\times \mathbb{Z}$ and $(y,n)=H(x,n)$.  Denote $(x_i,n_i):=F^i(x,n)$,  $(y_i,n_i):=G^i(y,n)$.

 Denote $\delta_i = |y_i-x_i|$ and $\tilde{\delta}_i=|A_{G(H(J_{i-1}))}(y_i)-A_{F(J_{i-1})}(x_i)|$. Here $(x_i,n_i)~\in~J_i \in~\mathcal{P}^0$. It is easy to conclude, using iii. and property  $LI$, that
 \begin{equation}\label{til} \tilde{\delta}_i \leq \frac{\delta_i}{|F(J_{i-1})|} + C\lambda^{|n_i|}\end{equation}
and making use of ii. to get
$$|\tau_{H(J_{i-1})}^G(A_{G(H(J_{i-1}))}(y_i))-\tau_{J_{i-1}}^F(A_{F(J_{i-1})}(x_i))| \leq D\tau^F_{J_{i-1}}(z_i) \frac{\delta_i}{|F(J_{i-1})|} +C\lambda^{|n_i|}. $$ Here $z_i \in [0,1]$. Since $D\tau^F_{J_{i-1}}(z_i)|F(J_{i-1})|/|J_{i-1}|\leq \lambda$ (property $Ex$), we get, using again $iii.$ \begin{equation} \label{rec} \delta_{i-1} \leq \lambda \delta_i + C\lambda^{|n_i|}.\end{equation}

Because $\psi$ is bounded, $|n_{i+1}-n_i|\leq B=\max |\psi|$. So if $i< n/2B$ then $|n_i| > |n_0|/2$. Since $\delta_{[\frac{n}{2B}]} \leq 1$, Eq. (\ref{rec}) implies
\begin{equation} \label{add2} |H(x,n)-(x,n)|=|y_0-x_0|\leq C\lambda^{\frac{|n|}{2}}.\end{equation}
In particular, by Eq. (\ref{til}) and property ii., we have
\begin{equation} \label{normalizado} |D\tau_{H(J_0)}^G(A_{G(H(J_0))}(y_1))-D\tau_J^F(A_{F(J_0)}(x_1))|\leq C\lambda^{\frac{|n|}{2}}.\end{equation}
By $Ra+Rb$ there exists $\theta \in (0,1)$ so that
\begin{equation}\label{lowerb} \theta^{|i|}\leq |I^n_i|.\end{equation}
Let $i$ be so that  $J=I^n_i$.

{\it Case A.} $|i|\geq |n/2|(\log \lambda/\log \theta)$: Due to i. and iii. and  property $Ra$, there exists $C > 0$ so that
$$ |\log \frac{|H(I^n_i)|}{|I^n_i|}| \leq C \lambda^n.$$ Together  with $sBD+LI$ and $iii.$, this implies that for every $p \in I^n_i$, with $|i|\geq |n/2|(\log \lambda/\log \theta)$,  we have
$$|\log \frac{DG(H(p))}{DF(p)}| \leq C\lambda^{\frac{|n|}{2}\frac{\log \lambda}{\log \theta}}.$$

{\it Case B.} $|i| < |n/2|(\log \lambda/\log \theta)$: In this case, by iii. and Eq. (\ref{lowerb}) we have
$$\log \frac{|H(I^n_i)|}{|I^n_i|} \leq  C\frac{|H(b^n_i)-b^n_i| +  |H(a^n_i)-a^n_i|}{|b^n_i - a^n_i|} \leq C\lambda^{\frac{|n|}{2}}. $$
Now using Eq. (\ref{add2}) and Eq. (\ref{normalizado}) we can easilly obtain
$$|\log \frac{DG(H(p))}{DF(p)}|\leq C\lambda^{\frac{|n|}{2}}.$$
\end{proof}

\section{Stability of transience}\label{sectiontransience}

We will begin this section with the large deviations result to transient homogeneous random walks and strongly transient random walks:

\begin{proof}[\bf Proof of Proposition \ref{largedeviationssthomo}] Let $P \in \mathcal{P}^0(F)$ be such that $F(P)=I_\ell$. By Proposition \ref{ldt} we have that for every $\epsilon > 0$ there exist $\gamma < 1$ such that 
$$m( p \in I_\ell \colon \  \pi_2(F^{n-1}p)-\pi_2(p)  <  (K-\frac{\epsilon}{2} ) \  (n-1) \leq C \gamma^{n-1},$$
for every $\ell$.  By the property BD we have
\begin{equation} \label{est_P} m( p \in P  \colon \  \pi_2(F^{n}p)-\pi_2(F(p))  <  (K-\frac{\epsilon}{2} )\  (n-1)) \leq C \gamma^{n-1} |P|.\end{equation}
Denote by $\Lambda_P^n$ the set in the l.h.s. of Eq. (\ref{est_P}). Let $n_0$ be such that $\min \psi > -\epsilon (n_0-1)/2 -\epsilon + K$. Then for $n\geq n_0$ we have 
\begin{equation} \label{est_P2} \tilde{\Lambda}_P^n:=\{ p \in P  \colon \  \pi_2(F^{n}p)-\pi_2(p)  <  (K -\epsilon ) \  n \}\subset \Lambda_P^n.\end{equation}
Indeed 
$$\pi_2(F^{n}p)-\pi_2(p) < (K-\epsilon) n$$ implies $$\pi_2(F^{n}p)-\pi_2(F(p))  <  (K -\epsilon/2  ) \  (n -1).$$ So
$$m( p \in P \colon \  \pi_2(F^{n}p)-\pi_2(p)  <  (K -\epsilon ) \  n) \leq C_2 \gamma^n |P|$$
for every $n$.
This completes the proof.
\end{proof}

\begin{proof}[\bf Proof of Proposition \ref{largedeviationsst}]  Fix $\epsilon > 0$ small. We intend to apply the Azuma-Hoeffding inequality, but since $\psi$ is not necessarily bounded, we need to make  some adjustments first: Fix $P \in \mathcal{P}^0(F)$ and define $\mathcal{F}_0:=\{P\}$ and $\mathcal{F}_n:=\{Q\}_{Q\subset P,\  Q \in \mathcal{P}^n(F)}$. Since $F\in GD$,  by the usual distortion control tricks for $F$, we can find $M > \min \psi$ such that $\alpha(x):= \min \{\psi(x),M  \}$ satisfies
\begin{equation}\label{perturbacao} \mathbb{E}(\alpha\circ F^n | \mathcal{F}_{n-1}) \geq K-\epsilon/4 \end{equation}
for every $n\geq 1$.  Here we are considering   conditional expectations relative to the probability
$$\mu_P(A):= \frac{m(A)}{|P|},$$
where $m$ is the Lebesgue  measure.

Define the martingale difference sequence

$$\Psi_n:= \alpha\circ F^n - \mathbb{E}(\alpha\circ F^n | \mathcal{F}_{n-1}).$$
Of course $||\Psi_n||_\infty \leq M$, if $M$ is large enough. By the Azuma-Hoeffding inequality we have

$$m( p \in P \colon \ |\sum_{i=1}^{n} \Psi_i(p) | > t ) \leq 2 \exp( -\frac{t^2}{2nM^2})|P|.$$
Taking $t=\epsilon n/4$ we obtain
\begin{equation}\label{dcorr}  m( p \in P \colon \ |\sum_{i=1}^{n} \Psi_i(p) | > \frac{\epsilon}{4} \ n ) \leq 2 \exp( -\frac{\epsilon^2 n}{32M^2})|P|.\end{equation}

Since

$$\pi_2(F^{n+1}p)-\pi_2(F(p))=  \sum_{i=1}^n \psi(F^i(p)) \geq \sum_{i=1}^{n} \alpha(F^i(p)) = \sum_{i=1}^n \Psi_i(p) + \sum_{i=1}^n \mathbb{E}(\alpha\circ F^i|\mathcal{F}_{i-1})(x)$$
$$\geq \sum_{i=1}^n \Psi_i(p)  + (K-\epsilon/4) n.$$
Due Eq. (\ref{dcorr}), this implies that
$$m( p \in P \colon \  \pi_2(F^{n}p)-\pi_2(F(p))=\sum_{i=1}^{n-1} \psi(F^i(p))  <  (K -\epsilon/2  ) \  (n -1)) \leq C_1 \exp( -\frac{\epsilon^2 n}{32M^2}) |P|.$$
Let $n_0$ be such that $\min \psi > -\epsilon (n_0-1)/2 -\epsilon + K$. Then for $n\geq n_0$ we have that  $$\pi_2(F^{n}p)-\pi_2(p) < (K-\epsilon) n$$ implies $$\pi_2(F^{n}p)-\pi_2(F(p))  <  (K -\epsilon/2  ) \  (n -1).$$ So
$$m( p \in P \colon \  \pi_2(F^{n}p)-\pi_2(p)  <  (K -\epsilon ) \  n) \leq C_2 \exp( -\frac{\epsilon^2 n}{32M^2}) |P|$$
for every $n$.
This completes the proof.
\end{proof}

\begin{prop}\label{prws}\label{homstr} Let $F$ be either a homogeneous random walk with positive mean drift or a strongly transient  random walk. Then any asymptotically small perturbation $G$ of $F$ has the following property: there exists $\lambda \in [0,1)$, $C > 0$  and $\tilde{K} > 0$  so that for every $P \in \mathcal{P}^0(G)$
$$m(p \in P \colon \ \sum_{i=0}^{n-1}\psi(G^i(p)) < \tilde{K} n  )\leq C\lambda^n |P|.$$
In particular $G$ is also transient.
\end{prop}
\begin{proof}
We will carry out the proof assuming the strong transience: the homogeneous case with positive mean drift  is analogous:
Fix $\epsilon > 0$. Let $\tilde{\delta}_1 > 0$ be small enough such that

$$(1-\tilde{\delta_1})(K-\epsilon) + \tilde{\delta_1} \min \psi> K-2\epsilon.$$

Due the bounded distortion of $G$, there exists $\delta_1 > 0$ such that  for every $n\geq 1$ and every $P \in \mathcal{P}^{n-1}(G)$,  interval $Q \subset G^n(P)$,  and set $A\subset Q$   satisfying
$$\frac{m(A)}{m(Q)} \geq 1 -\delta_1$$
we have
\begin{equation}\label{distg}\frac{m(P\cap G^{-n}A)}{m(P\cap G^{-n}Q)} \geq 1 -\tilde{\delta}_1.\end{equation}

 By Proposition \ref{largedeviationsst} we have 
\begin{equation}\label{unpum} m( p \in P \colon \  \sum_{i=0}^{n-1} \psi(F^i(p))  <  (K -\epsilon) n\ for \ some \ n\geq n_0) \leq C_1
\exp( -C_2n_0)|P|,\end{equation}
for every $P \in \mathcal{P}^0_j(F)$. Since $G$ is an asymptotically small perturbation, Eq. (\ref{asymp}) implies that
\begin{equation}\label{uv} m( p \in H(P) \colon \  \sum_{i=0}^{n-1} \psi(G^i(p))  <  (K -\epsilon) n\ for \ some \ n\geq n_0) \leq C_3
\exp( -C_4n_0) |H(P)|\end{equation}
provided that $P \in \mathcal{P}^0_j(F)$, $j \geq 2\ |\min \psi|\ n_0$. Indeed, the set in the l.h.s. of Eq. (\ref{unpum}) can the written as the pairwise disjoint union of the sets $\Delta_j$, $j\geq n_0$, where $\Delta_j$ is defined as
$$\{ p \in P \colon \  \sum_{i=0}^{n-1} \psi(F^i(p))  \geq   (K -\epsilon) n\ for \ every \ n_0\leq n< k   \ and \   \sum_{i=0}^{k-1} \psi(F^i(p)) <   (K -\epsilon)k \}$$
So by Eq. (\ref{asymp}) and Eq. (\ref{asymp2}) we have that 
$$dist_k(p)\leq C n_0 \lambda^{\ |\min \psi|\ n_0}+  \sum_{i=n_0+1}^{\infty} C \lambda^{(K -\epsilon) i} \leq  \tilde{C} < \infty$$
for every $ p\in \Delta_k$, $k\geq n_0$, and $j \geq 2\ |\min \psi|\ n_0$.  In particular 
\begin{equation}\label{cmeas} m(H(\Delta_k))\leq \tilde{C} m(\Delta_k).\end{equation}
Note that the set in the l.h.s. of Eq. (\ref{uv}) is the pairwise disjoint union of $H(D_j)$. Since $P \in \mathcal{P}^0_j$ we have $m(P)\leq C m(H(P))$, so from   Eq. (\ref{cmeas}) we obtain  Eq. (\ref{uv}).

In particular there exists $n_0=n_0(\delta_1)$ such that for every $P \in \mathcal{P}^0_j(G)$, $j \geq 2\ |\min \psi|\ n_0$,  we have
\begin{equation} \label{rato} m(\tilde{\Omega}_P)\geq  (1-\delta_1)|P|,\end{equation}
where
$\tilde{\Omega}_P$ is the set of points
 $p \in P$ such that  $\pi_2(G^n(p))\geq |\min \psi|\ n_0$ for all $n\geq 0$ and $\pi_2(G^n(p))-\pi_2(p)\geq (K -\epsilon)n$ for all  $n\geq n_0.$

By the $GD$ condition, there exists $n_1$ such that for $n \geq n_1$ we have

$$m(p \in P\colon \  there \ exists \ i\leq n \ s.t. \ \psi(F^i(p))\geq n)\leq \frac{\delta_1}{4}    $$

By Eq (\ref{unpum}) there exists $n_2>n_1$ such that

\begin{equation}\label{unpdois} m( p \in P \colon \  \sum_{i=0}^{n_2-1} \psi(F^i(p))  >  (K -\epsilon) n_2) \geq (1-\frac{\delta_1}{4})|P|.\end{equation}

So
\begin{equation}\label{unptres} m( p \in P \colon \  \sum_{i=0}^{n_2-1} \psi(F^i(p))  >  (K -\epsilon) n_2 \  and \ \psi(F^i(p))<  n_2 \ for \ every \ i\leq n_2) \end{equation}
$$\geq (1-\frac{\delta_1}{2})|P|.$$

Note that for  $p$ in the set in Eq (\ref{unptres})  we have $\pi_2(G^i(p))-\pi_2(p)\leq  (n_2)^2$ for every  $i\leq n_2$. Since $G$ is an asymptotically small perturbation of $F$, this observation and Eq. (\ref{unptres}) implies that there exists $n_3 >> (n_2)^2$ such that for $P \in \mathcal{P}^0_j(G)$, with $j\leq - n_3$, we have

\begin{equation} m( p \in P \colon \  \sum_{i=0}^{n_2-1} \psi(G^i(p))  >  (K -\epsilon) n_2 \  and \ \psi(G^i(p))<  n_2 \ for \ every \ i\leq n_2)\end{equation} $$ \geq (1-\delta_1)|P|.$$

So for  $P \in \mathcal{P}^0_j(G)$, with $j\leq - n_3$, we have
\begin{equation}\label{pum} m( p \in P \colon \  \sum_{i=0}^{n_2-1} \psi(G^i(p))  >  (K -\epsilon) n_2) \geq (1-\delta_1)|P|.\end{equation}

{\it Claim $A$:} Almost every point $x \in I\times \{j \}$, $j\leq -n_3$, visits at least once (and consequently infinitely many times) the set

\begin{equation}\label{tend} \bigcup_{j\geq-n_3} I\times \{j \}\end{equation}

Indeed, define a new random walk $\tilde{G}\colon I\times \mathbb{Z}\rightarrow I\times \mathbb{Z}$

$$\tilde{G}(x,n):=(\tilde{g}_n(x),n+\tilde{\psi}(x,n))$$
in the following way. Let $T$ be an integer larger than $n_2(K-\epsilon)$. If $n \geq -n_3$ then define $\tilde{g}_n\colon I \rightarrow I$ as an affine expanding map,  onto on  each element of $\mathcal{P}_n^{n_2}$, and $\tilde{\psi}(x,n)=T$.



For $(x,n)$, with $n< -n_3$, define  $\tilde{G}(x,n)=G^{n_2}(x,n)$. In this case
$$\tilde{\psi}(x,n)=\sum_{i=0}^{n_2-1} \psi(G^i(x,n)).$$

It is not difficult to see that the $\tilde{G}$-orbit of a point $(x,n)$, with $n< -n_3$, visits  the set in Eq. (\ref{tend}) at least once then  the $G$-orbit of $(x,n)$ visits the same set at least once.

To prove the claim, it is enough to show that $\tilde{G}$ is strongly transient. Indeed, let $P$ be an element of the Markov partition $\mathcal{P}^{k-1}_j(\tilde{G})$.  If $\pi_2(\tilde{G}^i(P))\geq -n_3$, for some $i\leq k$ then  $\pi_2(\tilde{G}^{k}(P))\geq -n_3$, so

\begin{equation} \label{stest1}\frac{1}{|P|} \int_P \tilde{\psi}\circ \tilde{G}^k \ dm=  \frac{1}{|P|} \int_P T \ dm \geq (K-\epsilon)n_2.\end{equation}
Otherwise $\pi_2(\tilde{G}^i(P))<  -n_3$ for every $i\leq k$. In particular $\tilde{G}^i=G^{i  n_2}$ on $P$, for every $i\leq k$. Note that $$\tilde{G}^{k}P= \bigcup_i Q_i,$$
where $\{ Q_i\}_i$ is  the family of all  interval $Q$ such that $Q \in \mathcal{P}^0_j(G)$ for some $j < -n_3$ and $Q\cap \tilde{G}^{k}P\not= \emptyset$   (this is a consequence of the Markovian property of $G$). By Eq. (\ref{pum})  we have

$$m(q \in Q_i\colon \tilde{\psi}(q) \geq (K-\epsilon)n_2)\geq(1-\delta_1)|Q_i|,$$
so by the distortion control in Eq. (\ref{distg}) we obtain

$$m(p \in P\cap \tilde{G}^{-k} Q_i\colon \tilde{\psi}(\tilde{G}^kp) \geq (K-\epsilon)n_2)\geq(1-\tilde{\delta}_1)|P\cap \tilde{G}^{-k} Q_i|,$$
consequently

\begin{equation}\label{stest2} \int_P \tilde{\psi}\circ \tilde{G}^k \  dm=\sum_i  \int_{P\cap \tilde{G}^{-k} Q_i}  \tilde{\psi}\circ \tilde{G}^k \  dm \end{equation}
$$ \geq \sum_i
((1-\tilde{\delta}_1)(K-\epsilon)n_2 + \tilde{\delta}_1 n_2 \min \psi) |P\cap \tilde{G}^{-k} Q_i| $$
$$\geq \sum_i (K -2\epsilon)n_2 |P\cap \tilde{G}^{-k} Q_i|= (K -2\epsilon)n_2 |P|$$

Eq. (\ref{stest1}) and (\ref{stest2}) imply that $\tilde{G}$ is strongly transient, so by Proposition \ref{largedeviationsst}, $\tilde{G}$ is transient. This concludes the proof of the claim.

{\it Claim $B$:} The $G$-orbit of almost every point of $I\times \mathbb{Z}$ eventually arrives at $\tilde{\Omega}_P$, for some $P \in \mathcal{P}^0_j$, with $j>  2|\min \psi|n_0$.

Since $F$ is  transient and $G$ is topologically conjugate to $F$ the set

$$\Omega:= \{ p \colon \  -n_3\leq \pi_2(p) \leq 2|\min \psi| n_0 \ and \ \lim_n \pi_2(G^n(p))=+\infty\}$$

is dense on
$$\bigcup_{j=-n_3}^{2|\min \psi| n_0 } I\times \{j\}.$$

This implies that for every non-empty open set $O \subset  I_j$, with $-n_3\leq j \leq 2|\min \psi| n_0$ we have
\begin{equation}\label{meio}
m((x,j) \in O \colon \exists \  k\geq 0 \ s.t. \ G^k(x,j) \in \tilde{\Omega}_P, \ with \ P \in \mathcal{P}^0_q(G), \ q >  2|\min \psi|n_0)> 0,
\end{equation}
where $\tilde{\Omega}_P$ is as in Eq. (\ref{rato}). Indeed, pick a point $p \in O\cap \Omega$.  By property $Ex$ and the definition of $\Omega$, there exists $k$ and $Q \in \mathcal{P}_j^k(G)$ such that $Q \subset O$, $P=G^k(Q)\in \mathcal{P}^0_q$, with $q>  2|\min \psi|n_0$. By Eq. (\ref{rato}) we have $m(\tilde{\Omega}_P)> 0$, so
$$m(O\cap G^{-k}\tilde{\Omega}_P)\geq m(Q\cap G^{-k}\tilde{\Omega}_P) > 0.$$

In particular there exists $\tilde{\delta}> 0$  such that for every interval $J  \subset I_j$, with $-n_3\leq j \leq 2|\min \psi| n_0$ and $|J|\geq \delta$, where $\delta$ is as in the LI property,  we have
\begin{equation}\label{meio2}
m((x,j) \in J \colon \exists \  k\geq 0 \ s.t. \ G^k(x,j) \in \tilde{\Omega}_P, \ with \ P \in \mathcal{P}^0_q(G), \ q >  2|\min \psi|n_0)\end{equation} $$> \tilde{\delta} |J|,
$$
It follows that there exists $\delta_3 > 0$ such that  for every $i$ and every $Q \in \mathcal{P}^{i-1}(G)$ such that $\pi_2(G^iQ)\geq  -n_3$ we have that
\begin{equation}\label{ubc} m(p \in Q\colon \ \exists k\geq 0 \ s.t. \ G^kp \in \tilde{\Omega}_P, \ with \ P \in \mathcal{P}^0_q(G), \ q >  2|\min \psi|n_0)\geq \delta_3 |Q|.
\end{equation}
Indeed, if $\pi_2(G^iQ)\leq 2|\min \psi| n_0$ we can apply Eq. (\ref{meio2}), BD and LI property.  Otherwise apply Eq. (\ref{rato}) and BD property. 

We will show Claim $B$ by contradiction. Suppose that it does not hold. Then there is a  set $W$ of  positive measure  whose $G$-orbit of its elements  never hits $\tilde{\Omega}_P$ for any $P\in \mathcal{P}^0_j$, with $j> 2 |\min \psi| n_0$. Pick a Lebesgue density point $p$ of $W$ whose $G$-orbit   visits
$$\bigcup_{j\geq -n_3} I\times\{j\}$$
infinitely many times, which is possible due Claim A. In particular  there exists a sequence $Q_k \in \mathcal{P}^{n_k-1}(G)$ such that $|Q_k|\rightarrow_n 0$, $p \in Q_k$, $\pi_2(G^{n_k}Q_k)\geq -n_3$ and

$$\lim_k \frac{m(Q_k\cap W)}{|Q_k|}=1.$$
That contradicts Eq. (\ref{ubc}). This concludes the proof of Claim $B$.

Note that Claim $B$ implies the following: almost every point in $I \times \{j\}$ belongs to the set
$$\Lambda_j := \bigcup_{k\geq 0} \Lambda_j^k,$$
where $$\Lambda_j^k:= \{p \in I \times \{j\}\colon \pi_2(G^n(p))-\pi_2(G^k(p))\geq (K-\epsilon)(n-k), \ for \ every \ n\geq k +n_0   \}.$$

Let $k_0$ be large enough such that for every $-n_3\leq j\leq 2|\min \psi|n_0$ we have

$$m(A\cap \bigcup_{k\leq k_0} \Lambda_j^k) \geq (1-\delta_1) |A|$$
for every interval $A\subset I\times \{j\}$ satisfying $|A|\geq \delta$, where $\delta >0$ is as in the property $LI$. Pick $n_4$ satisfying $n_4\geq k_0+n_0$ and
$$n_4 > \frac{-k_0\min \psi}{\epsilon}-k_0.$$

It is easy to see that if $p \in \bigcup_{k\leq k_0} \Lambda_j^k$ then
$$\pi_2(G^{n_4}p)-\pi_2(p)=\sum_{i=0}^{n_4-1}\psi(G^ip)\geq (K-2\epsilon)n_4.$$

In a argument similar to the proof of Claim $A$, consider  the random walk $\hat{G}$ defined in the following way: if $\pi_2(p)\leq -n_3$ define $\hat{G}(p)=G^{n_2}$. If $\pi_2(p)\geq  2|\min \psi|n_0$ define $\hat{G}(p)=G^{n_0}$. Finaly if $-n_3< \pi_2(p)<  2|\min \psi|n_0$ define $\hat{G}(p)=G^{n_4}$. The random walk $\hat{G}$ is $3\hat{K}$-strongly transient, for some $\hat{K} > 0$. The proof is quite similar to the proof of the strong transience of $\tilde{G}$, so we let it to the reader. So $\hat{G}$  is transient. It is easy to see that this implies that $G$ is transient. Finally Proposition \ref{bru} implies that
$$m( p \in P \colon \ \pi_2(\hat{G}^n(p))-\pi_2(p)  < 2\hat{K} n )\leq C\hat{\lambda}^n |P|,$$
for some  $\hat{\lambda} \in (0,1)$, which implies
$$m(Y^n_P)\leq C\hat{\lambda}^n |P|,$$
where
$$Y^n_P:=\{p \in P \colon \  \exists \ m \geq n \ s.t. \ \pi_2(\hat{G}^m(p))-\pi_2(p)  < 2\hat{K} m\}.$$
Let $n_5=\max\{n_0,n_4,n_2\}$.  Let $p \in P$ be such that
$$ \pi_2(G^i(p))-\pi_2(p)  < \frac{\hat{K}}{n_5}  i.$$
There exists $m$ and $j$  such that $\hat{G}^m(p)=G^j(p)$, with $i \geq  j$, $|i-j|\leq n_5$. Note that $$m\leq i\leq j + n_5\leq (m+1)n_5,$$
so we can find $i_0$ such that for every $i\geq i_0$ we have
$$ \frac{-n_5 \min \psi}{m} + \hat{K}\frac{m+1}{m}< 2\hat{K}.$$

 So
$$\pi_2(\hat{G}^m(p))-\pi_2(p)= \pi_2(G^j(p))- \pi_2(G^i(p))+ \pi_2(G^i(p))-   \pi_2(p)$$
$$\leq - n_5 \min \psi + \frac{\hat{K}}{n_5} i\leq- n_5 \min \psi + \hat{K}(m+1)<  2\hat{K} m, $$
where $$m\geq  \frac{i}{n_5}-1.$$
This implies
$$\{ p \in P\colon    \pi_2(G^i(p))-\pi_2(p)  < \frac{\hat{K}}{n_5}  i   \}\subset Y^{ \frac{i}{n_5}-1}_P,$$
so
$$m( p \in P\colon    \pi_2(G^i(p))-\pi_2(p)  < \frac{\hat{K}}{n_5}  i )\leq C\hat{\lambda}^{i/n_5}|P|$$
This completes the proof.
\end{proof}

Let $n > 0$ and $j$ be integers and $F$ be a deterministic random walk.  Then any connected component $C$  of $F^{-n}
\ int \ I_j$ is called a {\bf cylinder}. It follows from the Markovian property of $F$ that  a cylinder is a disjoint union of intervals in $\mathcal{P}^{n-1}$.  The {\bf lenght} $\ell(C)$ of the cylinder $C$ is $n$. If $C$ is a cylinder of lenght $n$ so that $F^i(C) \subset I_{j_i}$, for $i<n$, we will denote $C=C(j_0,j_1,\dots,j_n)$.

\begin{prop}\label{comeco} Let 
$F=(\{f_i\},\psi)\in Mk+LBD+LI+Ex+BD$. Assume that there exists $\epsilon > 0$ so that
 for $K
> 0$, we have
\begin{equation}\label{decaimento} m(\{ p \in I_n \colon \psi(p) < -K\}) \leq \frac{1}{K^{2+\epsilon}},\end{equation}
provided $n \geq n_0$.  Then
\begin{equation}\label{disj} \lim_k m( \{p \in I_{n_k}\colon \text{ there exists } i \leq k^2
\text{ so that } \psi(F^i(p)) < -k\})=0,\end{equation} uniformly for all
sequences satisfying $n_k
> k^3 + n_0$.
\end{prop} 
\begin{proof} For each $k$ and $i\leq k^2$, denote 
$$\Lambda_{n_k}^i=\{p \in I_{n_k}\colon \ \psi(F^j(p))\geq -k \text{ for every } j< i \text{ and } \psi(F^i(p)) < -k\}.  $$
The set in the l.h.s. of Eq. (\ref{disj}) is the   union of the sets  $\Lambda_{n_k}^i$.  The interval $I_{n_k}$ is the  union of the cylinders in $\mathcal{P}^{i-1}_{n_k}$. Let $Q \in \mathcal{P}^{i-1}_{n_k}$ and suppose that $Q\cap \Lambda_{n_k}^i\neq \emptyset$. Then $\pi_2(F^i(Q))\geq n_0$. By the property LI and Eq. (\ref{decaimento}) we get  
$$\frac{m(p \in F^i(Q)\colon \ \psi(p) < -k)}{m(F^i(Q))} \leq C\frac{1}{k^{2+\epsilon}}.$$
By the BD property 
$$ m(Q\cap \Lambda_{n_k}^i)\leq  \frac{C}{k^{2+\epsilon}} m(Q).$$ 
As a consequence
$$ m(I_{n_k}\cap \Lambda_{n_k}^i)\leq  \frac{C}{k^{2+\epsilon}}.
$$
So
$$ m(I_{n_k}\cap \cup_{i\leq k^2}\Lambda_{n_k}^i)\leq  \frac{C}{k^{\epsilon}}.$$

\end{proof}
\begin{rem} For a homogeneous random walk, the condition on $\psi$
is equivalent to $1_{I_0}\cdot\psi \in
L^{2+\epsilon}(m)$.\end{rem}

Let $F$ and $G$ be random walks which are topologically
conjugated by a homeomorphism $h$ that preserves states. For any $p
\in I \times \mathbb{Z}$ define
$$C_p:= \sup_{i\geq 0} dist_i(p).$$
For each $n_0 \in \mathbb{Z}\cup\{-\infty\}$ define
$$\Omega_{n_0+}(F):= \{p \colon \pi_2(F^n(p))\geq n_0, \text{ for
all } n \geq n_0\}.$$
In particular $\Omega_{-\infty+}(F)=I\times \mathbb{Z}$.
\begin{prop}\label{abs} Let $F$ and $G$ be random walks which are conjugated
by a homeomorphism $h$ which preserves states. Suppose that there
exists a $F$-forward invariant set $\Lambda$ so that
\begin{itemize}
\item[ ] \
\item[\it -H1:] $C_p:=\sup_{i\geq 0} \ dist_i(p)  < \infty$, for each
$p \in \Lambda$.\\
\end{itemize}
Then $h$ is absolutely continuous on $\cup_i F^{-i}\Lambda$ and
$h^{-1}$ is absolutely continuous on $\cup_i G^{-i}h(\Lambda)$.
Furthermore if also
\begin{itemize}
\item[ ] \
\item[\it -H2:]  There exists $C > 0$, $M > 0$ and $n_0 \in \mathbb{Z}\cup\{-\infty\}$ so that for every $n\geq n_0$ with  $n \in \mathbb{Z}$ and $P \in \mathcal{P}^0_n$,
 $$m(p \in P \cap \Lambda \colon \ C_p \leq C  ) \geq M |P|.$$
\end{itemize}
Then $h$ is absolutely continuous on $\cup_i
F^{-i}(\Omega_{n_0+}(F))$ and $h^{-1}$ is absolutely continuous on
$\cup_i G^{-i}(\Omega_{n_0+}(G))$. In particular when $n_0=-\infty$ we have that $h$ and $h^{-1}$ are absolutely continuous on $I\times \mathbb{Z}$.
\end{prop}
\begin{proof}For each $j \in \mathbb{N}$ denote
$$\Lambda_j:= \{ p \in \Lambda \colon  \sup_i \ dist_i(p) \leq
j\}.$$ Note that $\Lambda_i$ is forward invariant.

We claim that $h$ is absolutely continuous on $\Lambda_j$ and
$h^{-1}$ is absolutely continuous on $h(\Lambda_j)$. Indeed, for
each $p \in \Lambda_j$ and $k \in \mathbb{N}$, denote
$F^kp=(x_k,n_k)$. Denote by $J_k(x) \in \mathcal{P}^k$ the unique interval which
contains $x$  so that $F^k$ maps $J_k(x)$ diffeomorphically onto
$Q_k \subset I_{n_k}$. There is some ambiguity here if $x$ is in the boundary
of $J_k(x)$, but these points are countable, so they are
irrelevant for us.

If we use the analogous notation to $h(x)$ and $G$, we have
$h(J_k(x))=J_k(h(x))$ and, due the BD+LI property of
the random walks $F$ and $G$, there exist $C_1, C_2
> 0$ such that
\begin{equation}\label{lim_dist} C_1 e^{-dist_k(p)} \leq \frac{|h(J_k(x))|}{|J_k(x)|} \leq C_2
e^{dist_k(p)}.\end{equation} So, if $p \in \Lambda_j$ then
\begin{equation}\label{dis} C_1e^{-j} \leq \frac{|h(J_k(x))|}{|J_k(x)|} \leq
C_2 e^j, \ \text{ for all } k \in \mathbb{N}.\end{equation}

Let $A \subset \Lambda_j$ be a set with positive Lebesgue measure.
We claim that $h(A)$ also has positive Lebesgue measure. Indeed,
choose a compact set $K \subset A$ with positive Lebesgue measure.
Denote $U_k:= \cup_{x \in K} J_k(x)$. Since $|J_k(x)|\leq
\lambda^k$, we have that $\lim_k m(U_k)=m(K)$ and $\lim_k
m(h(U_k))=m(h(K))$. Since $U_k$ is a countable disjoint union of
intervals of the type $J_k(x)$, by Eq. (\ref{dis})
\begin{equation}\label{dist2}
C_1e^{-j} \leq  \frac{m(h(U_k))}{m(U_k)} \leq C_2 e^j, \ so \
C_1e^{-j}\leq  \frac{m(h(K))}{m(K)} \leq C_2
e^j,\end{equation} and we conclude that $h(K)$ also has positive
Lebesgue measure. An identical argument shows that, if $A \in
\Lambda_j$ has positive Lebesgue measure, then $h^{-1}A$ also has
positive Lebesgue measure. The proof of the claim is finished and
so $h$ and $h^{-1}$ are absolutely continuous on $\Lambda=\cup_j
\Lambda_j$ and $h(\Lambda)=\cup_j h(\Lambda_j)$.

Now it is easy to conclude that $h$ and $h^{-1}$ are absolutely continuous on
$\cup_i F^{-i}\Lambda$ and $\cup_i G^{-i}h(\Lambda)$.

Now assume $H2$. We claim that $\cup_i F^{-i}\Lambda$ has full
Lebesgue measure on $\Omega_{n_0+}(F)$. Indeed, Assume that $m(
\Omega_{n_0+}(F)\setminus \cup_i F^{-i}\Lambda) > 0$ and choose
 a Lebesgue density point $p$ of this set. Then
$$ \lim_k \frac{m(J_k(p)\cap \Omega_{n_0+}(F)\setminus \cup_i F^{-i}\Lambda )}{|J_k(x)|} =1.$$

Due the bounded distortion of $F$, if $F^k(p)=(x_k,n_k)$ and $F^k(J_k(x))=Q_k \subset I_{n_k}$, with $n_k\geq n_0$, where  $Q_k$ is a union of intervals in $\mathcal{P}^0_{n_k}$, then

$$ \limsup_k \frac{m(Q_k\cap \Lambda)}{|Q_k|}
\leq C(1- \liminf_k \frac{m(J_k(x)\cap \Omega_{n_0+}(F)\setminus
\cup_i F^{-i}\Lambda)}{|J_k(x)|}) =0,$$ which contradicts H2.

Since   $dist_k(p)$ is uniformly bounded with respect to $k$ and $p$ on the set $\{p \in P\cap \Lambda \colon \ C_p \leq C  \}$, we
can use an argument identical to the proof of Eq. (\ref{dist2}) to
conclude that
$$\frac{m(p \in P\cap \Lambda \colon \ C_p \leq C  )}
{m(h(p) \in h(P)\cap h(\Lambda) \colon \ C_p \leq C  )}\leq
C_1,$$ so $m(h(P\cap \Lambda\colon C_p\leq C)) \geq \tilde{C}M|h(P)|$, for all $P \in \mathcal{P}^0_n$, $n\geq n_o$ and
using an argument as above, we conclude that $\cup_i
G^{-i}h(\Lambda)$ has full Lebesgue measure on $\Omega_{n_0+}(G)$.
Since $h$ ($h^{-1}$) is absolutely continuous on $\cup_i
F^{-i}\Lambda$ ($\cup_i G^{-i}h(\Lambda)$) and
$$m(\Omega_{n_0+}(F)\setminus \cup_i
F^{-i}\Lambda)=m(h(\Omega_{n_0+}(F)\setminus \cup_i
F^{-i}\Lambda))=m(\Omega_{n_0+}(G)\setminus \cup_i
G^{-i}h(\Lambda))=0,$$
we have that $h$ and $h^{-1}$ are absolutely continuous on
$\Omega_{n_0+}(F)$ and $\Omega_{n_0+}(G)$. Now it is easy to prove that $h$ is absolutely continuous on $\cup_i F^{-i}\Omega_{n_0+}(F)$ and $h^{-1}$ is absolutely continuous on $\cup_i G^{-i}\Omega_{n_0+}(G)$.
\end{proof}

\begin{proof}[{\bf Proof of Theorem \ref{sttr}}] By Proposition \ref{homstr}, $G$ is transient. In particular for all $n_0 \in \mathbb{Z}$  the sets
$$\cup_i F^{-i}\Omega_{n_0+}(F) \text{ and }  \cup_i G^{-i}\Omega_{n_0+}(G)$$ have full Lebesgue measure. So  by Proposition
\ref{abs}, to prove that $h$ and $h^{-1}$ are absolutelly continuous, it is enough to find a forward invariant set satisfying
the assumptions H1 and H2 for some $n_0 \in \mathbb{Z}$. Indeed, fix $\delta > 0$ (we will choose $\delta$ later).
Consider the $F$-forward invariant set
$$\Lambda =\Lambda_\delta := \{ p \ \colon \
\liminf_k \frac{\pi_2(F^k(p))-\pi_2(p)}{k} \geq \frac{\delta}{3} \}.$$

We claim that $\Lambda$ satisfies H1. Indeed take $p \in \Lambda$.
Then, for $k \geq k_0(p)$ we have $n_k:=\pi_2(F^k(p)) \geq k
\delta/4$. So

\begin{equation} \label{distcont} dist_k(p)\leq  \sum_{i=0}^{k-1} |\log \frac{DF(F^{i+1}(p))}{ DG(h(F^{i+1}(p)))}|  \end{equation}
$$\leq \sum_{i=0}^{k_0-1} |\log \frac{ DF(F^{i+1}(p))}{ DG(h(F^{i+1}(p)))}|  +\sum_{i=k_0}^{k-1} |\log \frac{ DF(F^{i+1}(p))}{
DG(h(F^{i+1}(p)))}|$$
 $$\leq \sum_{i=0}^{k_0-1} |\log \frac{ DF(F^{i+1}(p))}{ DG(h(F^{i+1}(p)))}|  +
 \sum_{i=k_0}^{k-1}  \lambda^{n_i}$$
$$\leq \sum_{i=0}^{k_0-1} |\log \frac{ DF(F^{i+1}(p))}{ DG(h(F^{i+1}(p)))}|  +
 \sum_{i=k_0}^{\infty} \lambda^{i\delta/4 }$$
$$ \leq K_p + C(\delta).$$
Here $\lambda$ is as in Eq. (\ref{asymp}) and Eq. (\ref{asymp2}). 
To prove that $\Lambda$ satisfies H2, By Proposition \ref{largedeviationsst} for each $P \in \mathcal{P}^0_i$ we have
\begin{equation}\label{estexp} m( p \in P \colon \ \pi_2(F^k(p))-\pi_2(p)  < \delta k )\leq C\lambda^k |P|,\end{equation}
provided $\delta$ is small enough. From Eq. (\ref{estexp}) we obtain
\begin{equation}\label{estr2} \mu( p \in P \colon \ \pi_2(F^n(p))-\pi_2(p)  \geq  \delta n \text{ for all } n\geq n_0 )\geq  (1-C\lambda^{n_0}) |P|.\end{equation}
In particular, we have that, for every $n$,  \begin{equation}\label{estest} \pi_2(F^n(p))\geq \delta (n-n_0)  + \pi_2(p)+ n_0 \min \psi.\end{equation} in the set in Eq. (\ref{estr2}). Using the same argument as in Eq. (\ref{distcont}) we can easily obtain $H2$ from Eq. (\ref{estest}) and Eq. (\ref{estr2}), choosing $n_0$ large enough.
\end{proof}

\begin{proof}[{\bf Proof of Theorem \ref{sttrho}}] Since the mean drift is positive, by  the Birkhoff Ergodic Theorem $F$ is transient.  By Proposition \ref{homstr}, $G$ is also transient. Now the proof  goes exactly as the Theorem \ref{sttr}, except that to obtain Eq. (\ref{estexp}) we use Proposition   \ref{largedeviationssthomo} instead of Proposition \ref{largedeviationsst}.\end{proof}

\begin{proof}[{\bf Proof of Theorem \ref{sqr}}]  Let  $F$  be  either a $K$-strongly recurrent random walk or a homogeneous random walk with mean drift $K= \int \psi \ dm$. Let  $\epsilon < K$. 
By  Proposition \ref{prws}  there exists $\theta < 1$ such that 
for every $i$ we have
\begin{equation} \label{vital} m( p \in I_i \ \colon \ \frac{\pi_2(F^k(p))-\pi_2(p)}{k} \leq  \epsilon ) \leq  C\theta^{k}.\end{equation}
Using an argument as in the proof of Theorem \ref{sttr} we can conclude that
 \begin{equation} \label{equsei} m( p \in I_i \ \colon \ \frac{\pi_2(F^k(p))-\pi_2(p)}{k} \geq  \epsilon
 \ for \ k \geq k_0) \geq 1-  C\theta^{k_0}\end{equation}
 for every $i\geq 0$. By  Theorem \ref{sttr} and Theorem \ref{sttrho}     the conjugacy $h$ is absolutely continuous. Let $\delta= \sup_p dist_1(p)$. Then   Eq. (\ref{equsei}) implies that there exist $C> 0$ such that 
 \begin{equation} \label{bdinf} m( p \in I_i \ \colon \ dist_{k}(p) \geq \delta n+C \ for \ some \ k)\leq C\theta^n,\end{equation}
for $i\geq 0$. Denote $\Lambda_1:= \{ p \in I_i \ \colon h'(p) \leq  1\}$ and, for $n \geq 1$
 $$\Lambda_n:= \{ p \in I_i \ \colon \ {e}^{\delta(n-1)} < h'(p) \leq {e}^{\delta n}\}.$$ By Eq. (\ref{bdinf}) we have $m(\Lambda_n)\leq  C\theta^n$. Indeed, Let $J_k(p)$  be  as in the proof of Proposition \ref{abs}. By the Lebesgue differentiation Theorem for almost every $p$ we have
 $$\lim_k  \frac{|h(J_k(p))|}{|J_k(p)|} =h'(p).$$
 On the other hand, by Eq. (\ref{lim_dist})  we have that for almost every $p \in I_i$ outside the set in Eq. (\ref{bdinf})
 $$C_1  e^{-(n\delta+C) } \leq \frac{|h(J_k(p))|}{|J_k(p)|}  \leq C_2  e^{n\delta+C},$$
 so 
 \begin{equation}\label{lim_dist2} C_1  e^{-(n\delta+C) } \leq h'(p)  \leq C_2  e^{n\delta+C},\end{equation}
 in a subset of $I_i$ with measure larger than $1-  C\theta^n$. Of course this implies $m(\Lambda_n)\leq  C\theta^n$. Let $B \subset I_i$ be an arbitrary Lebesgue measurable set. Let $k_1$ be so that
 $$    \theta^{k_1+1}<  |B|\leq \theta^{k_1}.$$
First we prove Theorem \ref{sqr} assuming that  $e^\delta\theta < 1$. Since $h$ is absolutely continuous we have
 $$|h(B)|= \int_B h'\ dm$$
 $$ = \sum_{n=0}^{k_1} \int_{B\cap \Lambda_n} h' \ dm + \sum_{n=k_1+1}^{\infty} \int_{B\cap \Lambda_n} h' \ dm$$
 $$ \leq  \sum_{n=0}^{k_1} \theta^{k_1}e^{\delta n}  + \sum_{n=k_1+1}^{\infty} C(e^{\delta}\theta)^n$$
 $$ \leq C (e^{\delta}\theta)^{k_1}\leq  C |B|^{1+ \frac{\delta}{\ln \theta}}.$$
 Now if $B \subset J \in \mathcal{P}^n$ and  $F^n(J)=Q \subset I_i$, with $|Q|, |h(Q)|\geq C$ (due Property LI for $F$ and $G$), then due the bounded distortion of $F$ and $G$
 $$\frac{|h(B)|}{|h(J)|}\leq C \frac{|h(F^n(B)|}{|h(Q)|} \leq C \Big( \frac{ |F^n(B)|}{|Q|} \Big)^{1+ \frac{\delta}{\ln \theta}}\leq C\Big( \frac{ |B|}{|J|} \Big)^{1+ \frac{\delta}{\ln \theta}}.$$

 To prove a similar inequality to $h^{-1}$, define
 $$\tilde{\Lambda}_n:= \{ p \in I_i \ \colon \ {e}^{\delta(n-1)} < (h^{-1})'(p) \leq {e}^{\delta n}\}.$$
 of course
 $$h^{-1}\tilde{\Lambda}_n = \{ p \in I_i \ \colon \ {e}^{-\delta n} < h'(p) \leq {e}^{-\delta (n-1)}\},$$
 so by Eq. (\ref{bdinf}) and Eq. (\ref{lim_dist2}) we obtain
 $$m(h^{-1}\tilde{\Lambda}_n)\leq C \theta^n.$$
 In particular
$$m(\tilde{\Lambda}_n) = \int_{h^{-1}\tilde{\Lambda}_n} h'(x) \ dm \leq C (e^{-\delta}\theta)^n$$
Note that this argument gives us an exponential upper bound even if $\delta$ is large.

Now we can switch the roles of $F$ and $G$ to obtain the inequality to $h^{-1}$,
 which shows that $h$ is a mSQS-homeomorphism relative to the stochastic basis $\cup_n \mathcal{P}^n$.

To complete the proof when $e^\delta\theta \geq  1$ we do the following: find a continuous path of random walks $F_t$ with $F_0=F$ and $F_1=G$ and  so that for   every $t \in [0,1]$ we have that $F_t$ is an asymptotically small perturbation of $F$ and moreover  there exist $\epsilon > 0$ and  $\theta < 1$ such that  Eq. (\ref{vital}) holds for every random walk in this family.    Using the compactness of $[0,1]$ we can find a finite sequence of random walks $F_{t_0}=F, F_{t_1}, F_{t_2}, \dots F_{t_n}=G$ so that $F_{t_i}$ and $F_{t_{i+1}}$ are conjugated by a map $H_i$ such that
 
$$\delta_i := \sup_p \big|\frac{DF_{t_{i+1}}(H_i(p))}{DF_{t_i}(p)}\big|.$$
satisfy $e^{\delta_i}\theta < 1$.  So the conjugacy $H_i$  is mSQS with respect some dynamically defined stochastic basis. Composing these conjugacies we find a mSQS-conjugacy between $F$ and $G$.
 \end{proof}

\section{Stability of recurrence}

Let $F=(f,\psi)$ be a homogeneous random walk and let $G$ be an asymptotically small perturbation of $F$.  To avoid a cumbersome notation, in this section we make the convention that all inequalities holds only for large $n$.
Moreover in this section we assume that $\psi$ is unbounded. Recall that in this case we assume that asymptotically  small perturbations $G$ coincides with $F$ on negative states. The case where $\psi$ is bounded is similar.

The following is an easy consequence of the Central Limit Theorem for Birkhoff
sums (Proposition \ref{clt})

\begin{cor}\label{co}Let $a_n$ be a positive increasing sequence. Then
$$\mu( \frac{|S_n|}{\sqrt{n}} >  a_n) \leq Ce^{-\frac{\sigma^2 a_n^2}{2}} +
C\frac{1}{\sqrt{n}}.$$ Here
$$S_n(x) = \sum_{k=0}^{n-1} \psi(f^k(x)).$$
\end{cor}
\begin{proof} Use Proposition \ref{clt} and and note that the estimative $$\int_{-\infty}^v
e^{-\frac{u^2}{2}} \ du \leq C e^{-\frac{v^2}{2} }$$
holds for $v << 0$.

\end{proof}

Given $n \in \mathbb{N}$, split $[0,2n]\cap \mathbb{N}$ in $\sqrt{\log n}$ blocks
(called main blocks) , denoted $B_j$,  with length
$$\frac{n}{\log^{8j}n }, \ j=1,\dots, \sqrt{\log n},$$
and between the main blocks we put  little blocks $H_j$, called
holes, of length  $\log^4n$. These holes will warranty the
independence between the events in distinct main blocks. Put these
blocks in the following order:
$$\dots < B_{j+1} <  H_{j+1} <
B_j < H_{j}< \dots,$$ with $min \ B_{\sqrt{\log n}} = 0.$
 Note that we let  most of the second half of the interval $[0,2n]\cap \mathbb{N}$
  uncovered.

Define

$$S(j) = \sum_{i \in B_j} \psi\circ f^i$$
$$H(j) = \sum_{i \in H_j} \psi\circ f^i$$
Denote $|B_j|:= \max B_j - \min B_j$.

\begin{lem}\label{muitogrande} We have
$$\mu(\sum_{i=0}^{|B_j|} \psi\circ f^i \geq \frac{\sqrt{n}}{\log^{4j} n} \log^3n) \leq  C \frac{\log^{4j} n}{\sqrt{n}}.$$
\end{lem}
\begin{proof} This follows from  Corollary \ref{co}.
\end{proof}

\begin{prop}\label{aux1} For every $\epsilon > 0$ we have
$$\mu(S(j) > \frac{\sqrt{n}}{\log^{4j} n} \log^3n, \ for \ some \
j \leq \sqrt{\log n}) \leq C\frac{1}{\sqrt[2+\epsilon]{n}},$$
provided $n$ is large enough.
\end{prop}
\begin{proof} For $j \leq \sqrt{\log n}$  define
$$ \Lambda_j:=\{x \in I\colon \ S(j)(x) > \frac{\sqrt{n}}{\log^{4j} n} \log^3n   \}$$
$$ = \{x \in I\colon  \sum_{i < |B_j|} \psi\circ f^{i + \min B_j }(x)  > \frac{\sqrt{n}}{\log^{4j} n} \log^3n   \}  $$and for each $P \in \mathcal{P}^{\min B_j}$ denote $\Lambda_j(P):= \Lambda_j\cap P$.

Due Lemma \ref{muitogrande} and the bounded distortion of $f^{\min B_j}$ on $P$ we have
$$m(\Lambda_j(P)) \leq  C \frac{\log^{4j} n}{\sqrt{n}}|P|.$$
Summing on $j$ and $P$
$$m(\bigcup_{{ j }}\bigcup_{{ P }} \Lambda_j(P)) \leq \sqrt{\log n} \ \frac{\log^{4j} n}{\sqrt{n}} << C\frac{1}{\sqrt[2+\epsilon]{n}}.$$
\end{proof}

\begin{prop}\label{aux2} For every $\epsilon > 0$ and $d >0$  we have
\begin{equation}\label{mais} \mu(|\sum_{i \in H_j} \psi(f^i(x))| > \log^8 n, \ for \ some \
j \leq \sqrt{\log n}) \leq C\frac{1}{n^d},\end{equation} provided $n$
is large enough.
\end{prop}
\begin{proof} For   $i \in H_j-1$, with  $j \leq \sqrt{\log n}$, define
$$\Lambda_{i}^j:=\{x \in I\colon  |\psi(f^i(x))| > \log^4 n. \}.$$
By expanding and bounded distortion properties of $f$ and condition $GD$ we have
 that
$$\mu(\Lambda_i^j)\leq C\lambda^{\log^4 n}.$$
Since $|H_j|= \log^4 n$, if $x$ belongs to the set in Eq. (\ref{mais}) then $x \in \Lambda_{i}^j$, for some $i \in H_j-1$, with  $j \leq \sqrt{\log n}$. So
$$\mu(|\sum_{i \in H_j} \psi(f^i(x))| > \log^8 n, \ for \ some \
j \leq \sqrt{\log n})$$
$$ \leq \mu(\bigcup_{{ j \leq \sqrt{\log n}}} \ \bigcup_{{ i \in H_j-1}} \Lambda_i^j)$$

$$\leq  \sqrt{\log n}\ \log^4 n \ n^{\log \lambda \log^3n}$$
$$<< \frac{1}{n^d}, $$
where the last inequality holds for $n$ large enough.
\end{proof}

\begin{prop}[Independence between distant events] \label{inde}There exists $\lambda < 1$
so that the following holds: Let $C_1$ be a disjoint union of elements of $\mathcal{P}^{n-1}$ and  let  $C_2$ be a disjoint union of elements of $\mathcal{P}^{k-1}$. We 
have
$$\mu(C_1\cap f^{-(n+d)}C_2) = \mu(C_1
)\mu(C_2)(1+ O(\lambda^{d})).$$ Here $n=\ell(C_1)$.
\end{prop}
\begin{proof}  Let $J \in \mathcal{P}^{n-1}$. Since $F \in On$ we have  $f^n(J)=I$. Define the
measure $\rho(A):= \mu(f^{-n}A\cap J)/\mu(J)$. Note that by the
bounded distortion property of $f$, we have that $\log d\rho/dm$
is uniformly $\alpha$-Holder, that is,
\begin{equation}\label{uni_cone} |\log \frac{d\rho}{dm} (x) - \log \frac{d\rho}{dm}  (y)|\leq C|x-y|^\alpha,\end{equation}
where $C$ and $\alpha$ do not depend on $n$ and $C_1$.
Furthermore it is bounded by above by a constant which does not
depend on $n$. By the well-know theory of Ruelle-Perron-Frobenius
operators for Markov expanding maps (see for instance \cite{viana}), if $P$ is the
Perron-Frobenius-Ruelle operator of $f$, then there exists $\lambda
< 1$  so that
$$ P^d\frac{d\rho}{dm} = (1+ O(\lambda^d))\frac{d\mu}{dm}.$$
So
$$\frac{\mu(J\cap f^{-(n+d)}C_2)}{\mu(J)}$$
$$=\rho(f^{-d}C_2)= \int 1_{C_2}\circ f^d \  \frac{d\rho}{dm} dm $$
$$=\int 1_{C_2} \  P^d \frac{d\rho}{dm}  dm$$
$$=(1 + O(\lambda^d))\int 1_{C_2}  \ \frac{d\mu}{dm}  dm$$
$$=(1 + O(\lambda^d))\mu(C_2).$$
The constant $\lambda$ is the contraction of the  Ruelle-Perron-Frobenious operator  in certain cone of positive functions  and whose logarithm is $\alpha$-Holder continuous (see \cite{viana}). Since  all functions $\log \frac{d\rho}{dm}$ belongs to the very  same cone (Due Eq. (\ref{uni_cone})), $\lambda$ does not depend on $C_1$.  Since $C_1$ is a disjoint union of intervals $J\in \mathcal{P}^{n-1}$, we finished the proof.
\end{proof}

\begin{cor}\label{aux3} There exists $M > 0$ so that
$$\mu( S_j < \frac{\sqrt{n}}{\log^{4j} n}\ M \text{ for all } j \leq
\sqrt{\log n}) \leq C\big( \frac{2}{3} \big)^{\sqrt{\log n}}$$
\end{cor}
\begin{proof}Choose $M >0 $ so that
$$\frac{1}{\sqrt{2\pi}}\int_{-\infty}^{M}
e^{-\frac{u^2}{2}} \ du < \frac{2}{3}$$
 Consider 
$$\mathcal{C}_j:= \{ x \ s.t. \ \sum_{i=0}^{|B_j|} \psi\circ f^i(x) <  \frac{\sqrt{n}}{\log^{4j}
n}\ M\}.$$
Note that $\mathcal{C}_j$ is a the disjoint union of elements of $\mathcal{P}^{|B_j|-1}$. The Central Limit Theorem tells us that if $n$ is large enough then
$$\mu(\mathcal{C}_j) < \frac{2}{3}$$
for every $j \leq \sqrt{\log n}$.

Recall that between $B_j$ and $B_{j+1}$ there is a hole
$H_{j+1}$ with length $\log^4 n$.  Denote
$$\Lambda_j := \bigcap_{i=1}^j  f^{-\sum_{k=i+1}^{j}(|B_ k|+|H_k|)}\mathcal{C}_i$$
Note that $\Lambda_j$ is a disjoint union of elements of $\mathcal{P}^{|B_1|+ \sum_{k=2}^{j}(|B_ k|+|H_k|)-1}$ and
$$ \Lambda_j = \mathcal{C}_j\cap f^{-|B_j|-|H_j|}\Lambda_{j-1}.$$
Moreover
\begin{equation}\label{lambfinal} \Lambda_{ \sqrt{\log n}}=\{ x \ s.t. \  S_j < \frac{\sqrt{n}}{\log^{4j} n}\ M \text{ for all } j \leq \sqrt{\log n}\}.\end{equation}
By  Proposition \ref{inde} , we obtain
$$\mu(\Lambda_j)=(1+ O(\lambda^{|H_j|}))\mu(\mathcal{C}_j) \mu(\Lambda_{j-1}) $$
So by Eq. (\ref{lambfinal})
$$\mu( S_j < \frac{\sqrt{n}}{\log^{4j} n}\ M \text{ for all } j \leq
\sqrt{\log n}) $$
$$\leq \big( \frac{2}{3} \big)^{\sqrt{\log n}} (1 +
O(\lambda^{\log^4 n}))^{\sqrt{\log n}}  \leq C\big( \frac{2}{3}
\big)^{\sqrt{\log n}}$$
\end{proof}

\begin{prop}\label{sei} There exists $C > 0$ so that
$$\mu( x \in I \colon \text{ there exists }  i < \ell^3 \text{ so that }
\sum_{k=0}^{i} \psi\circ f^k(x) > \frac{\ell}{2}) \geq 1-
C\big(
\frac{2}{3} \big)^{\sqrt{3\log \ell}}$$
\end{prop}
\begin{proof} Let $M$ be as in Corollary \ref{aux3}. Denote $n = \ell^3$ and define
$$\mathcal{A}_\ell:=  \{x \colon \text{ there exists }  i < \ell^3 \text{ so
that } \sum_{k=0}^{i} \psi\circ f^k(x) > \frac{\ell}{2}\},$$
$$\mathcal{B}_\ell:= \{ x \colon \ |S_j| <  \frac{\sqrt{n}}{\log^{4j} n}
\log^3n, \text{ for  all } j \leq \sqrt{\log n}\},$$
$$\mathcal{C}_\ell:= \{x \colon S_j \geq  \frac{\sqrt{n}}{\log^{4j} n} \ M, \text{ for some } j \leq
\sqrt{\log n} \},$$
$$\mathcal{D}_\ell:= \{x \colon  |H_j(x)| \leq  \log^8 n, \text{ for all }
j \leq \sqrt{\log n} \}.$$ We claim that if $\ell$ is large then
$\mathcal{B}_\ell\cap \mathcal{C}_\ell \cap \mathcal{D}_\ell \subset \mathcal{A}_\ell$. Indeed, let $x \in
\mathcal{B}_\ell\cap \mathcal{C}_\ell \cap \mathcal{D}_\ell$. Then for some $j_0 \leq \sqrt{\log
n}$,
$$S_{j_0}(x) \geq \frac{\sqrt{n}}{\log^{4j_0} n} \ M.$$

We claim that, if $m = max \ B_{j_0}$, then
$$\sum_{0}^m \psi\circ f^i(x) > \frac{\ell}{2}.$$
Indeed, since  $x \in
\mathcal{D}_\ell$,
$$| \sum_{i\in H_{j}, \ j > j_0} \psi\circ f^i(x)|  \leq \sqrt{\log  n}
\log^8 n = o(\ell).$$ Moreover, since $ x \in \mathcal{B}_\ell$,
$$| \sum_{i \in B_{j}, \ j > j_0} \psi\circ f^i(x)| \leq \sum_{j > j_0}
 \frac{\sqrt{n}}{\log^{4j} n } \log^3 n \leq C \frac{\sqrt{n}}{\log^{4j_0 +1} n}. $$
So
$$\sum_{0}^m \psi\circ f^i(x) = \sum_{i \in B_{j_0}} \psi\circ
f^i(x) + \sum_{i \in B_{j}, \ j > j_0} \psi\circ f^i(x) + \sum_{i
\in H_{j}, \ j > j_0} \psi\circ f^i(x)$$
$$ \geq   \big(M- \frac{C}{\log n}\big) \frac{\sqrt{n}}{\log^{4j_0} n} + o(\ell) > C\ell^{\frac{6}{5}} - o(\ell)>  \frac{\ell}{2},$$
and we finished the proof of the claim. To finish the proof, note that by Proposition \ref{aux1}, Corollary \ref{aux3} and  Proposition \ref{aux2} $$\mu(\mathcal{A}_\ell) \geq \mu(\mathcal{B}_\ell\cap \mathcal{C}_\ell \cap \mathcal{D}_\ell) \geq 1-
C\frac{1}{\sqrt[2+\epsilon]{n}}-C\big(
\frac{2}{3} \big)^{\sqrt{\log n}} - C\frac{1}{n^d}  \geq 1- C\big(
\frac{2}{3} \big)^{\sqrt{\log n}}.$$
\end{proof}

Let $C> 0$ and $\lambda \in (0,1)$ be as in Eq. (\ref{asymp}) and Eq. (\ref{asymp2}). Define
$$Dist_n(p):= \sum_{i=0}^{n-1} C\lambda^{\pi_2(F^ip)}.$$
Of couse $dist_n(p)\leq Dist_n(p)$. 
\begin{prop} \label{fund_rec} There exist $\epsilon$ and $D$  so that for every $\ell \geq 0$,
$$\mu (\{ p \in I_\ell \colon \text{ there exists  }  i
 \text{ so that } F^i(p) \in \bigcup_{t \in [\min \psi, -\min \psi]} I_t  \text { and } dist_i(p) \leq D \}) \geq
 \epsilon$$
\end{prop}
\begin{proof} 
For $\ell\geq 0$ and $k$ define $B_k^{\ell}$ as the set of all $p \in I_\ell$  such that there exists $$ j \leq \sum_{i=0}^{k-1}
\frac{\ell^3}{2^{3i}}$$ satisfying  $$\pi_2(F^j(p)) \leq
\frac{\ell}{2^k} \text{ and } Dist_j(p) \leq  \sum_{i=0}^{k-1} C
\frac{\ell^3}{2^{3i}}\lambda^{\frac{\ell}{2^i}+\min \psi}.$$

 We are going to prove by ascending induction on $k\geq 0$ that there is  $C > 0$ so
that  for every $\ell\geq 0$ we have  
\begin{equation} \label{mest} \mu(B_k^{\ell}) \geq \prod_{i=0}^{k-1} \left( 1-
C\big(
\frac{2}{3} \big)^{\sqrt{\log \frac{\ell}{2^i}}}  \right), \end{equation}
for all $k\geq 1$ and $\mu(B^\ell_0)=1$. 

Note that $B^\ell_0=I_\ell$ so $\mu(B^\ell_0)=1$. Now assume the induction hypothesis  for some  $k\geq 0$. Take $p \in B_k^\ell$. Let $p \in
L=C(i_0,i_1,\dots,i_{j-1})$, where $j$ is the smallest integer as in the definition of
$B_k^\ell$. In particular  \begin{equation} \label{sanduba} \frac{\ell}{2^k} + \min \psi \leq  \pi_2(F^j(p))  \leq
\frac{\ell}{2^k}.\end{equation}

Note that $L \subset B_k^\ell$ and  $F^j(L)=I_r$,with $r:=\pi_2(F^j(p))$.  Applying  Proposition \ref{sei} to $-\psi$ we get 
\begin{equation}\label{oest} \mu( x \in I_r \colon \text{ there exists }  i <
\frac{\ell^3}{2^{3k}} \text{ so that } \sum_{n=0}^{i} \psi\circ
f^n(x) < - \frac{\ell}{2^{k+1}})\end{equation} $$ \geq 1-
C\big(
\frac{2}{3} \big)^{\sqrt{\log \frac{\ell}{2^k}}}.$$ Denote

$$D_L: = \{ x \in  L \colon \text{ there exists }  i <
\frac{\ell^3}{2^{3k}} \text{ so that } \sum_{n=0}^{i} \psi\circ
f^n(f^j(x)) < -\frac{\ell}{2^{k+1}}  \}$$ Due the bounded
distortion property for $F$, the estimative in Eq. (\ref{oest})
implies
\begin{equation} \label{estimative} \frac{\mu(D_L)}{|L|} \geq 1-
C\big(
\frac{2}{3} \big)^{\sqrt{\log \frac{\ell}{2^k}}}.\end{equation} 
We claim that $D_L \subset B^{\ell}_{k+1}$. Indeed, let  $x \in
D_L$.  Take the smallest $i$ so that $$\sum_{n=0}^{i}
\psi\circ f^n(f^j(x)) < -  \frac{\ell}{2^{k+1}}.$$

Then by Eq. (\ref{sanduba}) we have $\pi_2(F^{j+h}(p)) \geq \frac{\ell}{2^{k+1}}+\min \psi$, for every $0
\leq  h < i$, so

$$Dist_i(F^{j}(p)) \leq \sum_{h=0}^{i} C \lambda^{\pi_2(F^{j+h}(p))}
\leq C \frac{\ell^3}{2^{3k}} \lambda^{\frac{\ell}{2^{k+1}} +\min \psi}.$$

So $D_L \subset B^{\ell}_{k+1}$. Since  $B^{\ell}_{k}$ is a
disjoint union of cylinders $L$, the estimative in Eq.
(\ref{estimative}) implies that  Eq. (\ref{mest}) holds replacing $k$ by $k+1$. This concludes the induction step. 

Define

$$D:= \sum_{i=0}^{\infty}
C \frac{\ell^3}{2^{3i}}\lambda^{\frac{\ell}{2^i}+\min \psi} < \infty.$$
Let $k$
be so that $ 2^k \leq \ell \leq 2^{k+1}$. Now it is easy to check
that

$$\mu (\{ x \in I_\ell \colon \text{ there exists  }  i
 \text{ so that } F^i(p) \in I_0  \text { and } dist_i(p) \leq D \})$$
 $$ \geq C\mu(B^{\ell}_k) \geq \prod_{i=0}^{k-1} \left( 1-
C\big(\frac{2}{3} \big)^{\sqrt{\log \frac{\ell}{2^i}}}  \right) \geq C
\prod_{i=0}^{k-1} \left( 1- C\big(\frac{2}{3} \big)^{\sqrt{\log \frac{2^k}{2^i}}}
\right)$$ $$ \geq  \exp( -C\sum_{i=1}^{\infty} \big(\frac{2}{3} \big)^{\sqrt{i\log 2}}) >
\tilde{C}
> 0,$$
which finishes the proof.
\end{proof}

\begin{proof}[Proof of the Stability of Recurrence (Theorem \ref{strec})] Since $F$ is recurrent, its mean drift is zero. By Corollary \ref{sigma_nz} we can apply the Central Limit Theorem as in the introduction to conclude Eq. (\ref{recu1}) and Eq. (\ref{recu2}). Because $G$ coincides with $F$ on negative states,  the orbit by $G$  of almost every point $p$ satisfying $\pi_2(p) < 0$ will entry
$$\cup_{i\geq 0} I^i.$$
As a consequence the orbit by $G$  of almost every point $p$ visits this set infinitely  many times. 
Let $\ell \geq 0$.

By Proposition \ref{fund_rec}, there exist $D~>~0$ and $\epsilon >
0$ so that
$$A_\ell:= \{ p \in I_\ell \colon \text{ there exists } i \text{ so that } F^i(p)
\in \bigcup_{t =\min \psi}^{ -\min \psi} I_t  \text{ and } Dist_i(p) < D   \}$$ satisfies $\mu(A_\ell) >
\epsilon$, for all $\ell\geq 0$.

Consider a cylinder $C_F=C_F(\ell,k_1,\dots,k_{i-1},k_i) \subset
A_\ell$, with $C_F\not=\emptyset$, satisfying $|k_j| > -\min \psi$ for $0< j < i$, $\min \psi \leq k_i \leq -\min \psi$ and  $Dist_i(x) <
D$, for every $x \in C_F$. We claim that that corresponding
cylinder $C_G=C_G(\ell,k_1,\dots,k_{i-1},k_i)$ for the perturbed
random walk $G$ satisfies
$$\frac{1}{C} \leq \frac{|C_G|}{|C_F|} \leq C,$$
where $C$ depends only on $D$. Because we used $Dist_i(p)$ instead of $dist_i(p)$ in the definition of $A_\ell$, the set  $A_\ell$ is a disjoint union
of cylinders of this type, so we obtain that $B_\ell = H(A_\ell)$
satisfies $m(B_\ell) > C\epsilon > 0$, for all $\ell\geq 0$.

To prove that the set of points whose orbits returns infinitely
many times to $$\bigcup_{t = \min \psi}^{ -\min \psi} I_t$$  has full Lebesgue measure, it is enough to
prove that $\Lambda:=\cup_{j\geq 0, \ell } G^{-j}B_\ell$ has full
Lebesgue measure.

Indeed, assume by contradiction that $\Lambda$ is not full. Choose
a Lebesgue density point $p$ of the complement of $\Lambda$ and also satisfying $\limsup_k \pi_2(G^k(p))\geq 0$. Then
there exist a sequence of cylinders
$C_k\in \mathcal{P}^{k-1}$ so that $p \in C_k$ and
\begin{equation} \label{conv}    \frac{m(C_k\setminus \Lambda)}{|C_k|} \rightarrow_k
1. \end{equation}

But $G^k(C_k)=I_{\ell_k}$, with $\ell_k=\pi_2(G^k(C_k))$, and $m(I_{\ell_k} \cap B_{\ell_k}) \geq
C\epsilon |I_{\ell_k}|$. By the bounded distortion property

$$\frac{m(\Lambda \cap C_k)}{|C_k|} > \frac{m(G^{-k}B_{\ell_k} \cap C_k)}{|C_k|}
> \tilde{C} \epsilon,$$
which contradicts Eq. (\ref{conv}). Now we can use that $G$ is transitive and has bounded distortion to prove that $G$ is recurrent.
\end{proof}

\begin{proof}[\bf Proof of Proposition \ref{rigidity}] Since $F$ is recurrent, almost every point of $I^0$ returns to $I^0$ at least once. So the first return map $R_F\colon I^0 \rightarrow I^0$ is defined almost everywhere is $I^0$ and the same can be said about $R_G$. Of course, the absolutely continuous conjugacy $H$ also conjugates the expanding Markovian maps $R_F$ and $R_G$. Using the same argument used in Shub and Sullivan \cite{ss} and Martens and de Melo \cite{mm}, we can prove that $H$ is actually $C^1$ on $I^0$. Using the dynamics, it is easy to prove that $H$ is $C^1$ everywhere.
\end{proof}

\section{Stability of the multifractal spectrum}

\subsection{Dynamical defined intervals and root cylinders} When we are dealing with Markov expanding maps with {\em finite} Markov partitions, for each arbitrary interval $J$ we can find an element of $\cup_j \mathcal{P}^j$ which covers $J$ and has more or less the same size that $J$. Note that this is no longer true when the Markov partitions is infinite.  Since coverings by intervals are crucial in the study of the Hausdorff dimension of an one-dimensional set, this trick is very useful to estimate the dimension of dynamically defined sets, once we can replace an  arbitrary covering by intervals by another one with essentially the same metric properties but whose elements are themselves {\em dynamically defined} sets (cylinders).

Consider $j\geq 0$ and let $\{C_i\}_{i \in \Theta} \subset \mathcal{P}^j$ be a finite or countable family of cylinders  such that $W:=\bigcup_i \overline{C_i}$ is connected, $W \subset J \in \mathcal{P}^{j-1}$ and $F^{j}(int \ W)$ does not contain any point $d^n_i$ (as defined in property $Rb$). Then $W$ is called a dynamically defined interval (dd-interval, for short) of level $j$. Define the root cylinder of $W$ as the unique cylinder $C_{i_0}$ with the following property: if $\# \Theta =\infty$ then $W$ is a semi-open interval and $C_{i_0}$ will be the cylinder so that $\partial C_{i_0}\cap \partial W \neq \emptyset$. Otherwise $W$ is closed and let $C_{i_0}$ be the unique cylinder such that $F=\partial C_{i_0}\cap \partial W$ is the boundary of a semi-open dd-interval which contains $W$. The following Lemmas are an easy consequence of the regularity properties $Ra+Rb$ and it will be useful to recover the  trick described above for (certain) infinite Markov partitions. The proof is very simple.

\begin{lem}\label{rcia} For every $d \in (0,1)$  there exists $K > 1$  so that for every dd-interval $W:=\cup_i \overline{C_i}$  with root cylinder $C_{i_0}$ we have
\begin{equation}\label{comp_int_1} \frac{1}{K}\leq \frac{|W|^\alpha}{\sum_{i}|C_i|^\alpha}  \leq K\end{equation}
\begin{equation}\label{comp_int_2}  \frac{1}{K}\leq \frac{|C_{i_0}|^\alpha}{\sum_{i}|C_i|^\alpha}  \leq K\end{equation}
for every $1 \geq  \alpha \geq d$.
Indeed the constant $K$ depends only on $d$ and constants in the properties $Ra+Rb+Ex+BD$.
\end{lem}
\begin{proof} Due Property $Ra$, we can enumerate $C_i$ in such way that $C_0$ is the root cylinder of $W$ and $\partial C_{i+1}\cap \partial C_i \neq \emptyset$. Moreover if $j$ is the level of $W$ then $W \subset J \in \mathcal{P}^{j-1}$. Let  $F^j(J)=I_n$.  In particular  $F^j(C_i) \in \mathcal{P}^0_n$  and $F^j(W)=\cup_i \overline{F^j(C_i)}$ is a dd-interval of level $0$, with root cylinder $F^j(C_0)$ and $\partial F^j(C_{i+1})\cap F^j(\partial C_i) \neq \emptyset$. By property $Rb$ we have 
$$ \frac{|F^j(C_{i})|}{|F^j(C_{0})|}\leq C\lambda^i.$$
By the property $BD$ we have that 
$$ \frac{|C_{i}|}{|C_{0}|}\leq C\lambda^i,$$
so we obtain Eq. (\ref{comp_int_2}) since
$$|C_0|^\alpha \leq  \sum_{i}|C_i|^\alpha \leq |C_0|^\alpha \sum_{i=0}^\infty C\lambda^{d i}.$$
In particular for $\alpha=1$ we have
\begin{equation}\label{interval}    1 \leq \frac{|W|}{|C_{0}|} \leq C,  \end{equation} 
From Eq. (\ref{interval}) and Eq. (\ref{comp_int_2}) we can easily get Eq. (\ref{comp_int_1}) for every $d\leq \alpha\leq 1$. 
 \end{proof}

\begin{lem}\label{rci} Let $N$ be as in Properties $Ra+Rb$. For every $d \in (0,1)$ there exists $K > 1$ so that the following holds: For every interval  $J \subset I\times \mathbb{Z}$ there exists $m$ dd-intervals $W_j$, all of same level, with $m \leq 2N$, satisfying the following properties:
\begin{itemize}
\item[-] The interior of these dd-intervals are pairwise disjoint.
\item[-] The closure of the union of $W_j$ covers $J$:
$$J \subset \overline{\bigcup_j W_j}. $$
\item[-] We have
$$  \frac{1}{K} \leq \frac{\sum_{i=1}^m  |W_i|^{\alpha}}{|J|^\alpha} \leq K$$
for every $1\geq \alpha > d$.
\end{itemize}
Indeed the constant $K$ depends only on $d$ and constants in the properties $Ra+Rb+Ex+BD$.\end{lem}
\begin{proof} Let $\mathcal{P}^{-1}=\{I_n\}_n$. Define the  sequence of partitions $\mathcal{Q}^j$, $j\geq 0$, of  $I\times \mathbb{Z}$ in the following way: $\mathcal{Q}^{0}$ is the family of the connected components of 
$$I\times \mathbb{Z}\setminus \{c^n_i, d^n_i     \}_{i,n}$$ and an interval $Q$ belongs to $\mathcal{Q}^{j}$, $j\geq 1$ if there exists $P \in \mathcal{P}^{j-1}$ such that $Q$ is one of the connected components of 
$$P\cap F^{-j}\{c^n_i, d^n_i     \}_i,$$
where $n=\pi_2(F^jP)$. Here $c^n_i$, $d^n_i$ are as in Property $Ra$ and $Rb$. Note that each $P \in \mathcal{P}^{j-1}$ contains at most $2N$ intervals in $\mathcal{Q}^j$ and each $Q \in \mathcal{Q}^j$ is a dd-interval of level $j$. First we consider the case
\begin{equation} \label{cond_simples} \{ j \geq 0 \text{ such that }  \# \{ Q\in \mathcal{Q}^j  \colon \  J \cap Q \not= \emptyset \}\leq  2 \} =\emptyset.  \end{equation}
Let $n=\pi_2(J)$. Then $J$ intersects at least three connected components of  $$I_n\setminus \{c^n_i, d^n_i     \}_i,$$
so it contains one of the connected components of this set. In particular if  $$\{ W_j\}_j :=  \{Q \in \mathcal{Q}^0 \text{ such that } Q\cap J\not= \emptyset\},$$
then  by Property $Ra+Rb$ we have  $\max_j |W_j|\geq \delta$, so 
$$1 \leq \frac{\sum_{i=1}^m  |W_i|^{\alpha}}{|J|^\alpha} \leq \frac{2N}{\delta}.$$
If Eq. (\ref{cond_simples}) does not hold, let
$$j_0= \max \{ j \geq 0 \text{ such that }  \# \{ Q\in \mathcal{Q}^j  \colon \  J \cap Q \not= \emptyset \}\leq  2 \}. $$
Let $Q,R\in \mathcal{Q}^{j_0}$ be  such that $J \subset \overline{Q\cup R}$. Then
$$Q=\cup_i \overline{D}_i, \ R=\cup_i \overline{E}_i,$$
with $D_i, E_i \in \mathcal{P}^{j_0}$, $\partial D_i \cap \partial D_{i+1} \not= \emptyset$,  $\partial E_i \cap \partial E_{i+1} \not= \emptyset$ and $D_0$ and $E_0$ are the root intervals of $Q$ and $R$.  Let
$$i_Q:= \min \{ i\geq 0 \text{ such that } D_i \cap  J \not= \emptyset \}.$$
$$i_R:= \min \{ i\geq 0 \text{ such that } E_i \cap  J \not= \emptyset \}.$$
$$Q_{i_Q}=\cup_{i\geq i_Q} \overline{D}_i, \ R_{i_R}=\cup_{i \geq i_R} \overline{E}_i,$$
Note that $Q_{i_Q}$ and $R_{i_R}$ are dd-intervals of level $j_0$.  Without lost of generality, suppose that $|D_{i_Q}|\geq |E_{i_R}|$. 
 Then $\overline{Q_{i_Q}\cup R_{i_R}} \supset J$ is an interval.  Let $K$ be  the constant given  by Proposition \ref{rcia} for $\alpha > d$.  Then
\begin{equation}\label{aux_33} (1+ C_1)|D_{i_Q}|\leq  |D_{i_Q}|+|D_{i_Q+1}|\leq |Q_{i_Q}|\leq |Q_{i_Q}| +|R_{i_R}|\leq 2 K |D_{i_Q}|,\end{equation}
where the first inequality follows from Eq. (\ref{prop_rb}). 
We have three cases.

\noindent {\it Case 1.} Suppose that $Q\not=R$ and  the  intervals  are in the order 
$$        D_{i_Q} < D_{i_{Q}+1} < \cdots < E_{i_{R}+1} < E_{i_{R}}$$
then  $|J|\geq C_1|D_{i_Q}|$, otherwise
$$J \subset D_{i_Q}\cup D_{i_Q+1},$$
which  contradicts  $J\cap E_{i_{R}}\not= \emptyset$.
So 
$$   1 \leq \frac{|Q_{i_Q}|^\alpha +|R_{i_R}|^\alpha }{|J|^\alpha} \leq \frac{(4K|D_{i_Q}|)^\alpha}{C_1^\alpha|D_{i_Q}|^\alpha}\leq \frac{4K}{C_1}.$$

\noindent {\it Case 2.} Suppose that $Q\not=R$ and  the  intervals  are in the order
$$        \cdots  < D_{i_Q+1} < D_{i_{Q}} <   E_{i_{R}} < E_{i_{R}+1}<\cdots$$
Then $i_Q=i_R=0$ and there exists $y \in \partial D_{0}\cap \partial E_{0}$. By Properties $Ra+Rb+BD$   there exist $ [d,y],[y,e] \in \mathcal{P}^{j_0+1}$, with $[d,y]\subset D_0$, $[y,e]\subset E_0$ such that 
$$ C_3 |D_0|\leq |d-y|, |e-y| \leq C_2 |D_0|.$$
since $J$ intersects $D_0$ and $Q_0$ and at least three intervals in $\mathcal{Q}^{j_0+1}$ intersect $J$, we have that either $[d,y]$ or $[y,e]$ is contained on $J$.  So $|J|\geq C_2 |D_0|$. We conclude
$$   1 \leq \frac{|Q_{i_Q}|^\alpha +|R_{i_R}|^\alpha }{|J|^\alpha} \leq \frac{(4K|D_{i_Q}|)^\alpha}{C_2^\alpha|D_{i_Q}|^\alpha}\leq \frac{4K}{C_2}.$$

\noindent {\it Case 3.} Suppose that   $R=Q$, that is $J \subset Q_{i_Q}$. By Properties $Ra+Rb+BD$  
$$C_5 |D_{i_Q}|\leq  |D_{i_Q+1}|\leq C_4 |D_{i_Q}|,$$ 
Using $Ra+Rb+BD$ again, for every interval  $S \subset  D_{i_Q} \cup D_{i_Q+1}$ such that $S \in \mathcal{Q}^{j_0+1}$ we have
$$|S|\geq C_6 |D_{i_Q}|.$$
 Since  at least three intervals in $\mathcal{Q}^{j_0+1}$ intersect $J$ there is $S \subset  D_{i_Q} \cup D_{i_Q+1}$ with  $S \in \mathcal{Q}^{j_0+1}$  such that $S\subset J$. So $|J|\geq C_6 |D_{i_Q}|.$   We conclude 
 $$   1 \leq \frac{|Q_{i_Q}|^\alpha }{|J|^\alpha} \leq \frac{(4K|D_{i_Q}|)^\alpha}{C_6^\alpha|D_{i_Q}|^\alpha}\leq \frac{4K}{C_6}.$$ 
\end{proof}

\subsection{Dimension of dynamically defined sets}
Let $f \in Mk+ BD + Ex$ and denote by $\mathcal{P}^0$ its Markov partition.
Let $$\mathcal{I}:= \{ C_i\}_i \subset \cup_n \mathcal{P}^n$$ be a finite or countable family of
disjoint cylinders. Define the induced Markov map $f_{\mathcal{I}}\colon
\cup_i C_i \rightarrow I$ by
$$f_{\mathcal{I}}(x) = f^{\ell(C_i)-1}(x), \ if \ x \in C_i.$$
We can also define an induced drift function
$\Psi\colon \cup_{i} C_i \rightarrow \mathbb{Z}$
in the following way: Define, for $x \in C \in \mathcal{P}^n_0$,
$$\Psi_{\mathcal{I}}(x):= \sum_{i=0}^{n-1} \psi(f^i(x)).$$
Under the same conditions on $x$, define $N_{\mathcal{I}}(x)=n$.
The maximal invariant set of $f_{\mathcal{I}}$ is
$$\Lambda(\mathcal{I}):= \{x \in I \colon \ f^j(x) \in \bigcup_i C_i, \ for \ all \ j\geq 0   \}.$$
Denote by $HD(\mathcal{I})$ the Hausdorff dimension of the maximal
invariant set of $f_{\mathcal{I}}$.

We are going to use the following result

\begin{prop}[Theorem 1.1 in \cite{mu2}] We have
$$HD(\mathcal{J})= \sup \{HD(\mathcal{I})\colon \mathcal{I} \subset \mathcal{J}, \ \mathcal{I} \ finite \}.$$
\end{prop}

Before to give the proof of Proposition  \ref{ms} we need to introduce some tools which
are useful to estimate the Hausdorff dimension.

Let $\mathcal{I}$ be a finite collection of disjoint cylinders. Then there exists $\beta$ such that
$$\sum_{C \in \mathcal{I}} |C|^\beta =1,$$ we will call $\beta$ the {\bf virtual
Hausdorff dimension } of $f_{\mathcal{I}}$, denoted
$VHD(\mathcal{I})$. The virtual Hausdorff dimension is a nice way
to estimate $HD(\mathcal{I})$: indeed if $f_{\mathcal{I}}$ is linear on each interval of the Markov partition then these values coincide. When the distortion is positive,  these values remain related, as expressed in the following result (which is
included, for instance,  in the proof of Theorem 3, Section 4.2 of
\cite{pt}). 

\begin{prop}\label{vhd} Let $\mathcal{I}$ be a finite family of disjoint
cylinders. Then

$$|HD(\mathcal{I})-VHD(\mathcal{I})| \leq  \frac{d}{\log \lambda -
d},$$ where $$d := \sup_{C \in \mathcal{I}} \sup_{x, y \in C} \log
\frac{Df_{\mathcal{I}}(y)}{Df_{\mathcal{I}}(x)} \text{ and }
\lambda := \inf_{C \in \mathcal{I}} \inf_{x \in C} |Df_{\mathcal{I}}(x)|.$$
\end{prop}

Recall that if $\mathcal{I}$ is finite then $f_{\mathcal{I}}$ has an invariant probability measure $\mu_{\mathcal{I}}$ supported on its maximal invariant set $\Lambda(\mathcal{I})$ such that for any subset $S \subset \Lambda(\mathcal{I})$ satisfying  $\mu_{\mathcal{I}}(S)=1$ we have $HD(S)=HD(\mathcal{I})$ (see for instance \cite{pu}).

Note that for a homogeneous random walk $F$
$$\Omega_+^k(F)= \{k\}\times\{x \in I\ s.t.  \sum_{i=0}^{j} \psi(f^j(x))+ k\geq 0, \ for \ j\geq 0    \}$$
and
$$\Omega_{+\beta}^k(F)=$$
$$\{k\}\times \{  x \in  I\ s.t. \
\sum_{j=0}^{n-1}\psi(f^j(x)) + k \geq 0\, \ for \ all \ n\geq 0 \ and \ \underline{\lim}_{\ n} \frac{1}{n}\sum_{j=0}^{n-1}\psi(f^j(x)) \geq\beta\}.$$
Define $\pi_1(x,n):=x$. The following is an easy consequence of this observation:
\begin{lem}\label{aux} If $F$ is a homogeneous random walk then
$\pi_1(\Omega_{+}^0(F)) \subset \pi_1(\Omega_{+}^k(F))$ and $\pi_1(\Omega_{+\beta}^0(F)) \subset \pi_1(\Omega_{+\beta}^k(F))$, for all $k\geq 0$. Furthermore $$HD(\Omega_{+}^0(F)) =HD(\Omega_{+}^k(F))$$  and
$$HD(\Omega_{+\beta}^0(F)) =HD(\Omega_{+\beta}^k(F)).$$
\end{lem}

\begin{prop}\label{sim} Let $F$ be a homogeneous random walk.
Then there exists a sequence of finite families of cylinders $$\mathcal{\mathcal{F}}_s \subset \cup_i \mathcal{P}^i_0$$ so that
\begin{itemize}
\item[ ] \
\item[-]$\Lambda(\mathcal{F}_s) \subset  \Omega_+^0(F),$\\
\item[-] Denote  $\beta_n:= \int \Psi_{\mathcal{F}_s}\ d\mu_{\mathcal{F}_s }$. Then  $\beta_n > 0$. \\
\item[-] $\lim_{s\rightarrow \infty} HD(\mathcal{F}_s)=HD(\Omega_+^0(F)).$ \\
\end{itemize}
\end{prop}
\begin{proof}
Denote $d=HD \ \Omega_{+}^0(F)\leq 1$. Given any $s \in \mathbb{N}^\star$, $m_{d_s}(\Omega_{+}(F))=\infty$, where $d_s:=d(1-1/s) < 1$. Here $m_{D}$ denotes the $D$-dimensional Haussdorf measure. By Theorem 5.4 in \cite{f}, for each positive number $M$ we can find a
compact subset $\Lambda_s \subset \Omega_{+}^0(F)$ satisfying $m_{d_s}(
\Lambda_s) = M$. We may assume that $\Lambda_s$ does not have isolated points.  We will specify $M$ later.

In particular, for each $\epsilon$ small enough the following holds:
\begin{itemize}
\item[]\ \\
\item[i.]{\em For  every}  family of intervals $\{ J_i  \}_i$ which covers $\Lambda_s$, with $|J_i| < \epsilon$ we have
$$\frac{M}{2}\leq \sum_i |J_i|^{d_s}.$$
\item[ii.] {\em There exists } a family of intervals $\{ J_i  \}_i$, with $|J_i|\leq
\epsilon$, which covers $\Lambda_s$ and
$$\sum_i |J_i|^{d_s} \leq 2M.$$
Furthermore we can assume that $\partial J_i \subset \Lambda_s$.\\
\end{itemize}
Assume that  $d_s\geq d/2$. By Lemma \ref{rcia} and Lemma \ref{rci}, there exists some $K$ such that we can replace the special  covering  $\{ J_i\}$ in ii. by a
new covering by dd-intervals $\{W_i^\ell\}_{i, \ \ell}$, with root cylinders $R_i^\ell$, where
\begin{equation}\label{p1} J_i\cap \Lambda_s \subset \overline{ \bigcup_\ell W_i^\ell},\end{equation}
\begin{equation}\label{p2} W_i^\ell:=\bigcup_k \overline{C^{i\ell}_k}, \ for \ each  \ \ell \leq m_{i\ell} \leq 2N,\end{equation}

\begin{equation}\label{p3} \frac{1}{K} \leq \frac{\sum_\ell |R^{\ell}_i|^{d_s}}{|J_i|^{d_s}}\leq K,\end{equation}

\begin{equation}\label{p4} \frac{1}{K} \leq \frac{\sum_k |C^{i\ell}_k|^{d_s}}{|R^{\ell}_i|^{d_s}}\leq K,\end{equation}
Indeed we can replace  $W_i^\ell$ by a dd-subinterval of it, if necessary, in such way that   $R_i^\ell \cap \Lambda_s \neq \phi$ and Eq. (\ref{p1}), Eq. (\ref{p2}), Eq. (\ref{p3}) and Eq. (\ref{p4}) hold, except  perhaps the lower bound in Eq. (\ref{p3}), since the new root cylinder could be  smaller than the original one.  The above estimates, together with  the fact that $\{ W^\ell_i \}$ covers $\Lambda_s$ (up to a countable set) gives
\begin{equation} \label{eqsb}\frac{M}{2K^2}\leq \sum_{i,\ell,k} |C^{i\ell}_k|^{d_s} \leq 2K^2M.\end{equation}
The lower bound in Eq. (\ref{eqsb}) follows from  i. Since these intervals are cylinders, if necessary we can replace this family of cylinders  by a subfamily of disjoint cylinders which covers $\Lambda_s$ up to a countable number of points and such that each cylinder intersects $\Lambda_s$. Indeed we can choose a finite subfamily $\mathcal{F}_s:=\{C_r\}_r$  satisfying
\begin{equation} \label{eqsb1}\frac{M}{3K^2}\leq \sum_r |C_r|^{d_s} \leq 2K^2M.\end{equation}

Let's call this finite subfamily $\mathcal{F}_s$.  Note that, since $C_r\cap \Lambda_s\neq \emptyset$ we have that
$$\sum_{t=0}^{\ell} \psi(f^t(x)) \geq 0$$ for every $x \in C_r$ and $\ell \leq \ell(C_r)$. If 
$$\sum_{t=0}^{\ell(C_r)} \psi(f^t(x)) =0$$
for every $C_r$, choose a very small cylinder $\tilde{C}$ satisfying 
$$\tilde{C}\cap \bigcup_r C_r = \emptyset$$ 
and such that
$$\sum_{t=0}^{\ell} \psi(f^t(x)) \geq  0$$ for every $x \in \tilde{C}$ and $\ell < \ell(\tilde{C})$, and 
$$\sum_{t=0}^{\ell(\tilde{C})} \psi(f^t(x)) > 0$$
on $\tilde{C}$, and moreover
\begin{equation} \label{neqsb}\frac{M}{3K^2}\leq |\tilde{C}|^{d_s} + \sum_r |C_r|^{d_s} \leq 3K^2M.\end{equation}
Add $\tilde{C}$  to the family $\mathcal{F}_s$.  Then, if $\mu_s$ is the geometric invariant measure of $f_{\mathcal{F}_s}$, we have
$$\int \Psi_{\mathcal{F}_s} \ d\mu_s > 0.$$
We can find such $\tilde{C}$ because $F \in On+GD$  implies that there is at least a point  $x_0$ such that 
$$\min_{k\geq 0}  \sum_{i=0}^k \psi(f^i(x_0)) > 0.$$
and $x_0 \notin \Lambda_s$. By  Proposition  \ref{vhd} and Eq. (\ref{neqsb})

$$| HD(\Lambda(f_{\mathcal{F}_s}))- d_s|\leq -\frac{C}{\log \epsilon}.$$

Since $\epsilon$ can be taken arbitrary, we can choose $\mathcal{F}_s$ such that $$HD(\Lambda(f_{\mathcal{F}_s}))\rightarrow_s d.$$

\end{proof}

\begin{cor} If $F$ is a homogeneous random walk we have that
$$HD(\Omega_+(F))= \lim_{\beta\rightarrow 0^+} HD(\Omega_{+\beta}(F))=\sup_{\beta > 0} HD(\Omega_{+\beta}(F)).$$
\end{cor}
\begin{proof} Due Lemma \ref{aux}, it is enough to prove the Corollary for $k=0$. Of course $\Omega_{+\beta}^0(F)\subset \Omega_{+}^0(F)$ and $\beta_0 \leq \beta_1$ implies $\Omega_{+\beta_1}^0(F) \subset \Omega_{+\beta_0}^0(F)$, so
$$\lim_{\beta\rightarrow 0^+} HD(\Omega_{+\beta}^0(F))=\sup_{\beta > 0} HD(\Omega_{+\beta}^0(F)) \leq HD(\Omega_{+}^0(F)).$$
To obtain the opposite inequality, let $\mathcal{F}_s$ be as in Proposition \ref{sim}. Denote
$$\gamma_s := \int \Psi_{\mathcal{F}_s} \ d\mu_{\mathcal{F}_s}, \text{ and } W_n := \int N_{\mathcal{F}_s} \ d\mu_{\mathcal{I}}$$
and $\beta_s:= \gamma_s/W_s$. Then by the Birkhoff Ergodic Theorem there is subset $T_s \subset \Lambda(\mathcal{I}_n)$ such that $\mu_{\mathcal{F}_s}(T_s)=1$ and
$$\lim_k \frac{1}{k} \sum_{i=0}^{k-1} \psi(f^i(x)) = \lim_k \frac{ \sum_{j=0}^{k-1} \Psi_{\mathcal{I}_n}(f^j_{\mathcal{F}_s}(x))}{ \sum_{j=0}^{k-1} N_{\mathcal{I}_n}(f^j_{\mathcal{F}_s}(x))}=\frac{\gamma_s}{W_s}=\beta_s > 0.$$
for every $x \in T_s$. Since the Hausdorff dimension of $\mu_{\mathcal{F}_s}$ is equal to $HD(\mathcal{F}_s)$, we have that $HD(T_s)=HD(\mathcal{F}_s)$. Note also that
$$T_s \subset \Omega_{+\beta_s}^0,$$
which implies $HD(\mathcal{F}_s)\leq HD(\Omega_{+\beta_s}^0)$, so by the choice of $\mathcal{F}_s$, we conclude that
$$ HD(\Omega_{+}^0)= \lim_s \ HD(\mathcal{F}_s) \leq \overline{\lim}_s \  HD(\Omega_{+\beta_s}^0)\leq \sup_{\beta > 0} \  HD(\Omega_{+\beta}^0).$$
\end{proof}

\begin{proof}[{\bf Proof of Theorem \ref{multi}.}] Define
$$\Gamma_n(F) := \{ x \in \Omega_{+\beta}^k(F) \ s.t. \  \pi_2(F^i(x,k))\geq \frac{\beta}{2} i, \ for \ all \ i\geq n\}.$$
Of course $$\Omega_{+\beta}^k(F) = \bigcup_n \Gamma_n(F).$$
To prove the Theorem, it is enough to verify that $HD (\Gamma_n(F))=HD (\Gamma_n(G))$.  Indeed,  for every $\epsilon > 0$ and $\alpha \in (HD(\Gamma_n(F)),1)$ there exists  a covering of $\Gamma_n(F)$ by intervals $A_i$ so that
$$\sum_j |A_j|^\alpha \leq \epsilon.$$
Note that we can assume that $\partial A_j \subset \Gamma_n(F)$. Since $G$ is an asymptotically small perturbation of $F$, it is easy to see that $G$ also satisfies the properties $Ra+Rb$, replacing the points $c^n_i$ and $d^n_i$ by $h(c^n_i)$ and $h(d^n_i)$, and modifying the constant .  Indeed can choose  constants in the definitions of the properties $Ex+BD+Ra+Rb$ which works for both random walks, so we can take $K > 0$ in the statements of  Lemma \ref{rci} and Lemma \ref{rcia} in such way that it works for both random walks.

In particular (as in the proof of Proposition \ref{sim}) for each $A_j$ we can find  at most $2N$ dd-intervals  $$W_j^\ell:=\bigcup_k \overline{C^{j\ell}_k}, \ with \ \ell \leq m_j \leq 2N$$
which  satisfy
$$A_i\cap \Gamma_n(F) \subset \overline{\bigcup_{\ell} W^\ell_i},$$
and
$$\sum_{k, \ell} |C^{j \ell}_k|^\alpha\leq K |A_j|^\alpha.$$
Furthermore we can assume that the root $R_{j}^\ell$ of $W_j^\ell$ satisfies
\begin{equation}\label{abc}\frac{1}{K} \leq \frac{|R_{j}^\ell|^\alpha}{\sum_{k} |C^{j \ell}_k|^\alpha}\leq K\end{equation}
and  $R^{j}_\ell \cap \Gamma_n(F)\neq \emptyset$.

The constant $K$ does not depend on $\alpha$, $j$ or $\ell$. In particular the union of all cylinders $C^{j\ell}_k$  covers $\Gamma_n(F)$ up to a countable set and
\begin{equation}\label{hdq} \sum_{j,k,\ell} |C^{j\ell}_k|^\alpha \leq K\epsilon.\end{equation}
Note that if $x \in \Gamma_n(F)$ then
$$dist_i(x)\leq r_n:=Cn + C\lambda^n$$
for every $i \in \mathbb{N}$. So
$$e^{-r_n}  \leq \frac{|\mathcal{P}^i_F(x)|}{|\mathcal{P}^i_G(h(x))|} \leq e^{r_n}.$$

There is a point in the cylinder $R_j^\ell$ which belongs to $\Gamma_n(F)$, so
\begin{equation}\label{hdd} e^{-\alpha r_n}\leq \frac{|R_j^\ell|^\alpha}{|h(R_j^\ell)|^\alpha} \leq e^{\alpha r_n}.\end{equation}
Note that $h(W_j^\ell)=\bigcup_k \overline{h(C^{j\ell}_k)}$ is a dd-interval for $G$ and $h(R_j^\ell)$ is its root cylinder. So, using Eq. (\ref{abc})
\begin{equation}\label{fgh} \frac{1}{K}\leq \frac{|h(R_j^\ell)|^\alpha}{\sum_{i}|h(C^{j\ell}_i)|^\alpha}  \leq K\end{equation}
But the union of the cylinders $h(C^{j\ell}_k)$ covers $\Gamma_n(G)$ up to a countable set and Eq. (\ref{abc}),  Eq. (\ref{hdq}), Eq. (\ref{hdd}) and  Eq. (\ref{fgh})  gives
$$\sum_{j,k,\ell} |h(C^{j\ell}_k)|^\alpha \leq K^3 e^{\alpha r_n}\epsilon.$$
Since $\alpha > HD(\Gamma_n(F))$ and $\epsilon$ is arbitrary we obtain that $HD(\Gamma_n(G)) \leq HD(\Gamma_n(F))$. Switching the roles of $F$ and $G$ in the above argument gives the opposite inequality.
\end{proof}

\begin{lem}\label{ol} Let $G \in On+Ra+Rb$ be a random walk. For every $\alpha > 0$ there exist $\epsilon$ and $C$
 so that 
\begin{equation}\label{best} \sum_{P \in \mathcal{P}^n_\ell} |P|^{1-\epsilon} \leq C (1+
\alpha)^n,\end{equation} for all $n$ and $\ell$.
\end{lem}
\begin{proof} Indeed, denote  \begin{equation} \label{part23} \mathcal{P}^{n}_\ell=\{ Q^j \}_j \text{  and } \mathcal{P}^{n+1}_\ell=\{ Q^j_k \}_{j,k},\end{equation} in such way that $Q_k^j \subset Q_j$. To avoid cumbersome notation we are omitting explicit indexing on $n$ and $\ell$. Since $G \in BD+Ra+Rb$, it is possible to order $Q_k^j$ so that  there exist $C$ and $\lambda < 1$ satisfying

\begin{equation} \label{bound23}   \frac{|Q^j_k|}{|Q^j|}\leq C\lambda^k,\end{equation}
for every $j, k, n$.   As a consequence the family of functions 
$$h_{j,\ell,n}(\epsilon)= \sum_k
\frac{|Q^j_k|^{1-\epsilon}}{|Q^j|^{1-\epsilon}}$$ is an
equicontinuous  set of functions in a small neighborhood of
 $0$. In particular, since $h_{j,\ell,n}(0)=1$, there exists $\epsilon_0$ so
 that, for every  $ \epsilon < \epsilon_0$ and every $j$, $\ell$ and $n$
\begin{equation} \label{uni} \sum_k \frac{|Q^j_k|^{1-\epsilon}}{|Q^j|^{1-\epsilon}} \leq
1+\alpha.\end{equation}
So
$$\sum_{P \in \mathcal{P}^{n+1}_\ell } |P|^{1-\epsilon} = \sum_{j, k } |Q^j_k|^{1-\epsilon} \leq (1+\alpha)  \sum_j |Q_j|^{1-\epsilon}= (1+\alpha)\sum_{P \in \mathcal{P}^{n}_\ell } |P|^{1-\epsilon}.$$
\end{proof}

From now on we are going to assume that $F=(f,\psi) \in  On$ is a homogeneous random walk with negative mean drift and $G$ is an asymptotically small perturbation of $F$.

\begin{lem}\label{bom} Let $G \in On+Ra+Rb$ be a  random walk that is an asymptotically small perturbation of a homogeneous random walk $F \in On+Ra+Rb$  with negative mean drift. Then for every $\alpha > \int \psi \ d\mu$, there exists $C> 0$,  $\sigma < 1$
 so that for any $n_1 \geq n_0$, with  $n_0$ large
enough, \begin{equation}\label{nzero} m\{p \in I_ {n_1}\colon \
\pi_2(G^k(p))\geq n_0, \text{ for } k\leq n,  \text{ and } \pi_2(G^n(p))
- n_1 \geq \alpha n \}\leq C\sigma^{n}.\end{equation}
\end{lem}
\begin{proof} Denote $$\Lambda_{n_0,n_1}^n(G):=\{p \in I_ {n_1}\colon \
\pi_2(G^k(p))\geq n_0 \text{ for all } k \leq n \text{ and }
\pi_2(G^n(p)) - n_1  \geq \alpha n \}.$$ The statement for $F$ is
consequence of the large deviations estimative (see, for instance
\cite{broise}) for every $K > 0$  there exists $C_K > 0$, $\gamma_K \in (0,1)$ such that 

$$m\{ p \in I\colon | \frac{\sum_{k=0}^{n-1} \psi(f^k(p))}{n} - \int \psi \ d\mu | \geq
K \}\leq C_K\gamma_K^n$$  
Pick $K= \alpha-\int \psi \ d\mu$ and  $\tilde{\sigma}=\gamma_K$.  Then for every $n_1$

$$m\{p \in I_ {n_1}\colon
\frac{\sum_{k=0}^{n-1} \psi(f^k(\pi_1(p)))}{n}=\pi_2(F^n(p)) -n_1 \geq \alpha n \}\leq C \tilde{\sigma}^{n},
$$
which implies (of course)
\begin{equation}\label{unp} m(\Lambda_{n_0,n_1}^n(F))\leq C \tilde{\sigma}^{n}.
\end{equation}
We are going to use this estimative to obtain Eq. (\ref{nzero})
for the perturbation of $F$.

Indeed, for every $\delta > 0$, there is $n_0$ so that
if $\pi_2(x)\geq n_0$  then
\begin{equation}\label{distor} 1-\delta \leq \frac{|DF(x)|}{|DG(H(x))|} \leq 1+
\delta,\end{equation}
 Here $H$ is the topological conjugacy between $F$
and $G$ which preserves states. Note that $\Lambda_{n_0,n_1}^n(F)$
is a disjoint union of elements $Q_i \in \mathcal{P}^n(F)$, so
$\Lambda_{n_0,n_1}^n(G)$ is a disjoint union of the intervals
$H(Q_i)$. Due the property BD of $F$ and $G$, Eq. (\ref{unp}) and Eq. (\ref{distor}), we have

\begin{equation} m(\Lambda_{n_0,n_1}^n(G))= \sum_i |H(Q_i)|\leq \sum_i C (1+\delta)^n |Q_i| \leq  C (1+\delta)^n \tilde{\sigma}^n.\end{equation}
Choose    $n_0$ large enough such that $\sigma:=(1+\delta)\tilde{\sigma}< 1$ . \end{proof}

We would like to replace $n_0$  by an arbitrary state in Eq.
(\ref{nzero}). The following Lemma will be useful for this task:

\begin{lem}\label{absa} Let $p_n$ and $q_n$ sequences of non-negative real
numbers such that
\begin{enumerate}
\item $p_0+ q_0 \leq 1$,
 \item There exists $\epsilon > 0$ and
$\ell\geq1$ such that $s_n:= p_n + q_n \leq
(1-\epsilon)^{\ell} p_{n-\ell} + q_{n-\ell}$ for every $n\geq \ell$ and $q_n \leq
C(1-\epsilon)^n + \sum_{k=1}^{n} (1-\epsilon)^kp_{n-k}$ , for every
$n$.
\end{enumerate}
Then there exists $C> 0$ and $\delta=\delta(\epsilon) > 0$ such that $s_n \leq C (1-\delta)^n$,
for every $n \in \mathbb{N}$.
\end{lem}
\begin{proof} If $n \geq \ell$, we have $s_n \leq (1-\epsilon)p_{n-\ell}+ q_{n-\ell}=
(1-\epsilon)s_{n-\ell} + \epsilon q_{n-\ell}$. It follows by
induction that if $n= i\ell + r$, with $r < \ell$, then  $$s_n
\leq (1-\epsilon)^i s_r +
\sum_{k=0}^{i-1}\epsilon(1-\epsilon)^{k\ell} q_{n-(k+1)\ell}$$
$$\leq C(1-\epsilon)^{n/\ell} s_r +
\sum_{k=0}^{n-\ell}\epsilon(1-\epsilon)^{k} q_{n-\ell -k}$$ Since
$q_{n-\ell} \leq C(1-\epsilon)^{n-\ell} + \sum_{k=1}^{n-\ell}
(1-\epsilon)^kp_{n-\ell-k}$, we obtain
$$s_n\leq C(1-\epsilon)^{n/\ell} s_r + C\epsilon (1-\epsilon)^{n/\ell} + \sum_{k=1}^{n-\ell} \epsilon(1-\epsilon)^{k}(p_{n-\ell-k}
 + q_{n-\ell-k})$$
 $$\leq (1-\epsilon)^{n/\ell}C(s_r + \epsilon)+ \sum_{k=1}^{n-\ell} \epsilon(1-\epsilon)^{k}
s_{n-\ell-k},$$ for every $n\geq \ell$.

We claim that there exists $\delta < 1$ and $K$ so that $s_n \leq
K(1-\delta)^n$, for every $n$.  Indeed, fix $\delta < 1$, For each
$n$, define  $K_n := s_n/(1-\delta)^n$. Note that

\begin{equation} \label{final} s_n \leq (1-\epsilon)^{n/\ell}C(s_r + \epsilon)+
\sum_{k=1}^{n-1} \epsilon(1-\epsilon)^{k} s_{n-\ell-k}$$ $$ \leq
(1-\epsilon)^{n/\ell}C(s_r + \epsilon)+ \sum_{k=1}^{n-\ell}
\epsilon(1-\epsilon)^{k}
K_{n-\ell-k}(1-\delta)^{n-\ell-k}\end{equation}
$$\leq  \big[ \big(
\frac{(1-\epsilon)^{1/\ell}}{1-\delta} \big)^n C(\max_{ j<\ell} s_j + \epsilon)+
\max_{i< \ n-\ell}K_i \ \frac{\epsilon}{(1-\delta)^\ell}
\sum_{k=1}^{n-\ell} \big(\frac{1-\epsilon}{1-\delta}\big)^k
\big](1-\delta)^n$$

Choose $\delta > 0$ close enough to $0$ so that

$$\sigma_1:= \frac{(1-\epsilon)^{1/\ell}}{1-\delta} < 1, \ and $$

$$\sigma_2:=  \frac{\epsilon}{(1-\delta)^\ell}
\sum_{k=1}^{\infty} \big(\frac{1-\epsilon}{1-\delta}\big)^k < 1.$$
Then by Eq. (\ref{final}) we have $K_{n} \leq \sigma_2\max_{i< \
n-\ell}K_i + C\sigma_1^n$, for every $n > \ell$,  which easily
implies that $\max_i K_i < \infty$.

\end{proof}

Define $$\Omega_{+}^{n_1,n}:= \{p \in I_{n_1}\colon \pi_2(G^{k}(p)) \geq 0,
\text{ for } 0 \leq k \leq n \}.$$

\begin{lem}\label{estimativa} Let $G \in On+Ra+Rb$ be a  random walk that is an asymptotically small perturbation of a homogeneous random walk $F \in On+Ra+Rb$  with negative mean drift. Then there exists $\delta < 1$ so that for every $n_1\geq
0$ there exists $C=C(n_1)$ satisfying
$$m(\Omega_{+}^{n_1,n}(G))\leq C(1-\delta)^n.$$
\end{lem}
\begin{proof}Take $n_0$ as in Lemma \ref{bom} and fix $n_1 \geq 0$. Define the sets and sequences

$$s_n:= m(\Omega_{+}^{n_1,n}) $$

$$p_n:= m(B^n), \text{ where } B^n:=\{ p \in  \Omega_{+}^{n_1,n}\colon \ \pi_2(G^{n}(p))\in [0,n_0] \}, \ and $$

$$q_n:= m(C^n), \text{ where } C^n:=\{p \in  \Omega_{+}^{n_1,n}\colon \  \pi_2(G^{n}(p)) > n_0  \}.$$

To  prove Lemma \ref{estimativa}, it is enough to verify that
these sequences satisfy the assumptions of Lemma \ref{absa}.
Indeed, of course $p_0 + q_0 \leq 1$. To prove the other
assumptions, take $i \in [0,n_0]$. Since $G$ is topologically
transitive, there are $\ell_i \in \mathbb{N}$ and intervals $J_i
\subset I_i$ so that $\pi_2(G^{\ell_i}(J_i)) < 0$. Denote
$\ell=max_{\ 0\leq i\leq n_0} \ell_i$ and $r= min_{\ 0\leq i\leq
n_0} |J_i|/|I_i|$.

Clearly $\Omega^{n_1,n}_{+}=B^n\cup C^n \subset B^{n-\ell}\cup
C^{n-\ell} $. Let $J \subset B^{n-\ell}$ be an interval so that
$G^{n-\ell}(J)=I_i$, with $0\leq i\leq n_0$. Note that
$B^{n-\ell}$ is a disjoint union of such intervals. By the bounded
distortion control for $G$,

\begin{equation}\label{ind} \frac{m(J\cap \Omega_{+}^{n_1,n})}{m(J)} \leq
 1- \frac{m(J\cap G^{-(n-\ell)}J_i)}{m(J)}\leq (1-\frac{r}{c})\end{equation}
Choose $\epsilon_0$ satisfying $(1-r/c)\leq (1-\epsilon_0)^\ell$.
 Then  Eq. (\ref{ind}) implies $$m(B^{n-\ell}\cap \Omega_{+}^{n_1,n})\leq
(1-\epsilon_0)^\ell m(B^{n-\ell})$$ and we obtain

$$s_n = m(B^{n-\ell}\cap \Omega_{+}^{n_1,n}) + m(C^{n-\ell}\cap
\Omega_{+}^{n_1,n}) \leq (1-\epsilon_0)^\ell p_{n-\ell} +
q_{n-\ell}.$$

It remains to prove that $q_n \leq \sum_{k=1}^{n}
(1-\epsilon)^kp_{n-k}$. There are two kind of points $p$ in $C^n$:

{\em Type 1.} For every $j\leq n$ we have   $\pi_2(G^j(p))\geq n_0$
(in particular $n_1\geq n_0$). We are going to estimate the
measure of the set of these points, denoted $\Theta_1^n$. It
follows from Lemma \ref{bom}, choosing $\alpha=\int\psi\ d\mu/2 < 0$, that

\begin{equation}\label{acima} m(\{p \in I_{n_1}\colon \ \pi_2(G^k(p))\geq n_0, \ for \ k\leq n
\text{ and } \pi_2(G^n(p)) \geq n_1  +\alpha n\}) \leq
C\sigma^n.\end{equation}
Note that if $n\geq (n_0-n_1)/\alpha$ then $n_1 + \alpha n\leq n_0$. Then  the set in the l.h.s. of Eq. (\ref{acima}) contains 
$\Theta^{n}_{1}$. In particular 
$$m(\Theta_1^n)\leq C_{n_1}\sigma^n,$$
for some $\sigma < 1$ which does not depend on $n_1$.

{\em Type 2.} For some $j \leq n$ we have   $\pi_2(G^j(p)) \leq n_0$.
Denote the set of these points by $\Theta^n_2$. Denote by
$\Theta_{2,k}^n$ the set of points $p$ so that $k\geq 1$ is the
smallest natural satisfying  $\pi_2(G^{n-k}p)\leq n_0$. Clearly
$\Theta^n_2$ is a disjoint union of these sets. We are going to
estimate their measure. Note that  $\Theta_{2,k}^n \subset
B^{n-k}$. The set $B^{n-k}$ is a disjoint union of intervals $L$
so that $\pi_2(G^{n-k}L)=I_i$, for some $i \leq n_0$. To estimate

$$\frac{m(\Theta_{2,k}^n\cap L)}{|L|} $$ note that $L \subset B^{n-k},$
and $\Theta_{2,k}^n\cap L$ is the set of points $p \in L$ so that
$\pi_2(G^{n-k+j}p) > n_0$, for every $0< j \leq k$. Define

$$L_y := \{ p \in L\colon \psi(G^{n-k}p)=y\}.$$

Firstly note that for  $y \leq n_0 -i$ we have \begin{equation}
\label{estum} |L_y\cap \Theta_{2,k}^n|=0,\end{equation} since  $p
\in L_y\cap \Theta_{2,k}^n$ satisfies $\pi_2(G^{n-k+1}p)= i+
\psi(G^{n-k}p)= i + y
> n_0$. In particular for $y < 0$ we have $|L_y\cap \Theta_{2,k}^n|=0$, which implies, due the bounded
distortion control

$$\frac{m(L\cap \Theta_{2,k}^n)}{|L|}\leq \frac{\sum_{y \geq 0}
|L_y|}{|L|} \leq (1-\delta),$$ for some $\delta < 1$ which does
not depends on $k$, $L$ or $n_1$, which implies
\begin{equation}\label{estquatro} m(\Theta^n_{2,k}) \leq
(1-\delta)m(B^{n-k})= (1-\delta)p_{n-k}.\end{equation}

Furthermore, using again the distortion control and the regularity
condition $GD$(big jumps are rare) we have
\begin{equation}\label{estdois} \frac{\sum_{y > -\alpha (k-1)} |L_y \cap \Theta_{2,k}^n|}{|L|}\leq
\frac{\sum_{y > -\alpha (k-1)} |L_y|}{|L|}\leq
C\gamma^{k},\end{equation} for some $C\geq 0$ and $\gamma < 1$.

To estimate $|L_y\cap \Theta_{2,k}^n|/|L_y|$, in the case $n_0-i
\leq y \leq -\alpha (k-1)$, recall that $G^{n-k+1}L_y=I_{i+y}$,
with $i+y
> n_0$. By Lemma \ref{bom}, we have
$$ m\{p \in I_ {i+y}\colon \ \pi_2(G^m(p))\geq n_0, \text{ for } m\leq
k-1,  \text{ and } \pi_2(G^{k-1}(p))
 \geq i+y + \alpha (k-1) \}\leq C\sigma^{k}.$$
 Since $i+y+ \alpha (k-1) \leq n_0$, this implies that

 $$ m\{p \in I_ {i+y}\colon \ \pi_2(G^m(p))\geq n_0, \text{ for every } m\leq
k-1 \}\leq C\sigma^{k}.$$

The points in $L_y\cap\Theta^n_{2,k}$ are exactly the points whose
$(n-k+1)$th-iteration belongs to the set in the estimate above.
Using the bound distortion control we have

$$ \frac{|L_y\cap \Theta^n_{2,k}| }{|L_y|} \leq C\sigma^k,$$
so

\begin{equation}\label{esttres}
\frac{|\sum_{n_0-i \leq y \leq -\alpha (k-1)}L_y\cap
\Theta^n_{2,k}| }{|L|} \leq  C  \frac{|\sum_{n_0-i \leq y \leq
-\alpha (k-1)}L_y\cap \Theta^n_{2,k}| }{\sum_{n_0-i \leq y \leq
-\alpha (k-1)} |L_y|} \leq   C\sigma^k.
\end{equation}

Choose $\epsilon < \epsilon_0$ so that $ min\{ max\{ C\sigma^k,
C\gamma^k\}, 1-\delta\} \leq (1-\epsilon)^k$, for every $k \geq
0$, and put together Eq. (\ref{estum}), Eq. (\ref{estquatro}), Eq.
(\ref{estdois}) and Eq. (\ref{esttres}), to get  $m(L\cap
\Theta_{2,k}^n)\leq (1-\epsilon)^k |L|$. Since $B^{n-k}$ is  a
disjoint union of such intervals $L$, we obtain
$$m(\Theta_{2,k}^n)\leq (1-\epsilon)^k m(B^{n-k})=
(1-\epsilon)^kp_{n-k}$$ and now we can conclude with
$$q_n = m(\Theta^n_{1}) + \sum_{k} m(\Theta^n_{2,k}) \leq C_{n_1}\sigma^n + \sum_k (1-\epsilon)^kp_{n-k}.$$
\end{proof}

Now we are ready to prove Theorem \ref{menor}.

\begin{proof}[{\bf Proof of Theorem \ref{menor}.}] There are three cases:

{\bf $F$ is transient with $ M > 0$.} If $M > 0$ then  the random walk $F$ is transient and it is easy to see (using for instance Proposition \ref{ldt}) that $m(\Omega_+(F)) > 0$. Since  the conjugacy  with an asymptotically small perturbation $G$ is absolutely continuous (Theorem \ref{abscont}), we conclude that  $m(\Omega_+(G)) > 0$.

{\bf $F$ is recurrent  ($ M = 0$).} If $M=0$ then $F$ is recurent \cite{guivarch}  and its asymptotically small perturbations are   recurrent by Theorem \ref{strec}. In particular almost every point visits negative states infinitely many times, so $m(\Omega_+(G)) = 0$. It remains to prove that $HD \ \Omega_+(G)=1$. By Theorem \ref{omega} it is enough to verify that  $HD \ \Omega_+(F)=1$.  Indeed, it is easy to show using the Central Limit Theorem that if
$$\int \psi \ d\mu=0$$
then there exist $C > 0$  and for each $n$, subsets $\mathcal{A}_n \subset \mathcal{P}^n_0$ so that
$$\sum_{i=0}^{n-1}\psi(f^i(x)) > 0$$
for all $x \in J \in \mathcal{A}_n$ and
\begin{equation}\label{tamanho} 1\geq m(\bigcup_{J \in \mathcal{A}_n} J) > C >0.\end{equation}
here $C$ does not depend on $n$. Of course we can assume that $\mathcal{A}_n$ is finite. Property $Ex$ implies that  there exists $\theta \in (0,1)$ such that 
$$\sup_{J \in \mathcal{A}_n}  |J| \leq \theta^n.$$
Consider the function $$h(\epsilon):= \sum_{J \in \mathcal{A}_n} |J|^{1-\epsilon}.$$
Then by Eq. (\ref{tamanho}) if $0\leq \epsilon <  1$ we have
$$h'(\epsilon):= \sum_{J \in \mathcal{A}_n} -\log |J| |J|^{1-\epsilon}\geq -Cn \log \theta.$$
In particular if 
$$\tilde{\epsilon}:=  \frac{C-1}{Cn \log \theta}$$
then $h(\tilde{\epsilon})\geq 1$. Since $h(0)\leq 1$ there exist $\epsilon_n = 1- O(1/n)$ such that $h(\epsilon_n)=1$. But $VHD(\mathcal{A}_n)=\epsilon_n$, so $$\big| VHD(\mathcal{A}_n) -1\big| \leq  \frac{C}{n}.$$
By property $BD$ that there exists $C_1 > 0$ such that for every $n$
$$d_n := \sup_{C \in \mathcal{A}_n} \sup_{x, y \in C} \log
\frac{Df_{\mathcal{A}_n}(y)}{Df_{\mathcal{A}_n}(x)}\leq C_1$$
and since $\mathcal{A}_n \subset \mathcal{P}^n$, by property  $Ex$ we have  that there exists $\theta \in (0,1)$ such that for every $n$ 
$$\lambda_n := \inf_{C \in \mathcal{I}} \inf_{x \in C} |Df_{\mathcal{I}}(x)| \geq \frac{1}{\theta^n}.$$
we can apply   Proposition \ref{vhd} to obtain
$$\big| HD \ \Lambda(\mathcal{A}_n) - VHD(\mathcal{A}_n)\big| =  O(\frac{1}{n}).$$
so
$$ HD(\mathcal{A}_n) = 1- O(\frac{1}{n}).$$
If $\mu_{\mathcal{A}_n }$ is the geometric invariant measure of $f_{\mathcal{A}_n}$ then
$$\int \psi_{\mathcal{A}_n} \ d\mu_{\mathcal{A}_n} > 0.$$
So by the Birkhoff Ergodic Theorem

\begin{equation} \label{conv1} \lim_{n\rightarrow \infty} \sum_{i=0}^{n-1}\psi(f^i(x)) = + \infty \end{equation}
in a set $S_n \subset \Lambda(\mathcal{A}_n)$ satisfying $\mu_{\mathcal{A}_n}(S_n)=1$, so $HD \ S_n = 1 - O(1/n)$. In particular the set $S$ of points satisfying Eq.(\ref{conv1}) has Hausdorff dimension $1$. We can decompose $S$ in subsets $B_j$ defined by
$$B_j :=\{ x \in S\colon   {\min}_{n} \sum_{i=0}^{n-1} \psi(f^i(x)) \geq -j \}.$$
Clearly $\sup_j HD \ B_j=1$.

By properties $GD+On$, for each $j$ there are $k_j$ and $J_j\not= \emptyset  \in \mathcal{P}^{k_j}$ so that for all $x \in J_j$ we have

$$\sum_{i=0}^{\ell-1}\psi(f^i(x)) \geq 0$$
for every $\ell\leq k_j$ and

 $$\sum_{i=0}^{k_j} \psi(f^i(x)) \geq  j.$$
 Then $$(J_j \cap f^{-k_j}B_j)\times \{0\}$$
 belongs to $\Omega_+(F)$, for every $j$. This implies $HD \ \Omega_+(F) \geq HD \ B_j$ so

 $$ HD \ \Omega_+(F) \geq \sup_j HD \ B_j =1.$$

{\bf $F$ is transient with  $ M < 0$.} By Lemma
\ref{estimativa}, there is some $\delta \in (0,1)$, which does not depend on
$n_1$, so that

\begin{equation} \label{useum} m(\Omega^{n_1,n}_+)\leq
C(1-\delta)^n.\end{equation}

By Lemma \ref{ol}, there exists $\epsilon$ so that

\begin{equation}\label{usedois} \sum_{P \in \mathcal{P}^n, \ P \subset I_k} |P|^{1-\epsilon}
\leq C(1-\delta)^{-n/2}.\end{equation}

Denote by $\{J_i^n\}_i \subset \mathcal{P}^n$ the family of
disjoint intervals so that $\Omega^{n_1,n}_+ = \cup_i J_i^n$. We
claim that there exists $C > 0$ satisfying

\begin{equation} \sum_i  |J_i^n|^{1-\epsilon/4} \leq C(1-\delta)^n.\end{equation}
Since $sup_i \ |J^n_i|\rightarrow_n 0$, this proves that $HD \
\Omega_+^{n_1,\infty} \leq 1-\epsilon/4$.

Indeed,

$$\sum_i  |J_i^n|^{1-\epsilon/4} = \sum_{ |J_i|> (1-\delta)^{2n/\epsilon}}  |J_i^n|^{1-\epsilon/4} +
\sum_{ |J_i|\leq (1-\delta)^{2n/\epsilon}}
|J_i^n|^{1-\epsilon/4}$$
$$\leq  (1-\delta)^{n/2} \sum_{i}  |J_i^n| +
(1-\delta)^{3n/2} \sum_{i} |J_i^n|^{1-\epsilon}$$
$$\leq C(1-\delta)^{n/2},$$
where in the last line we made use of Eq. (\ref{useum}) and Eq.
(\ref{usedois}). The proof is complete.
\end{proof}

\section{Applications to one-dimensional renormalization theory}\label{applications}

\subsection{(Classic) infinitely renormalizable maps}\label{apl} Denote $I=[-1,1]$. Consider a real analytic unimodal
maps $f\colon I \rightarrow I$, with negative Schwarzian
derivative and even order critical point at $0$. The map $f$ is called
infinitely renormalizable if there exists an sequence of natural
numbers $n_0 < n_1 < n_2 < \dots$ and a nested sequence of
intervals
$$I=I_0 \supset I_1 \supset I_2 \supset \cdots $$
so that
\begin{itemize}
\item$f^{n_k}\partial I_k \subset \partial I_k$,
\item $f^{n_k}I_k
\subset I_k$,
 \item $f^{n_k}\colon I_k \rightarrow I_k$ is a unimodal map.
\end{itemize}

We say that $f$ has bounded combinatorics if there exists $C
>0$ so that $n_{k+1}/n_k \leq C$, for all $k$. Two infinitely renormalizable
maps $f$ and $g$ have the same
combinatorics if there exists a homeomorphism $h\colon I
\rightarrow I$ such that $f\circ h = h \circ g$.

The following result is a deep result in renormalization theory:

\begin{prop}[\cite{mc2}]\label{convren} Let $f$ and $g$ be two infinitely
renormalizable unimodal maps with the same bounded combinatorics
and same even order. Then for every $r > 0$ there exists $C > 0$
and $\lambda < 1$ so that
$$||\frac{1}{|I_k^f|}\ f^{n_k}(|I_k^f| \cdot ) -  \frac{1}{|I_k^g|}\ g^{n_k}(|I_k^g| \cdot )||_{C^r}
 \leq C\lambda^k.
 $$
\end{prop}
Here $|I^f_k|$ denotes the length  of $I_k^f$.

\begin{figure}
\centering
 \psfrag{f}{$f$}
 \psfrag{g}{$f^2$}
 \psfrag{h}{$f^4$}
\psfrag{p1}[][]{$p_1$}
 \psfrag{p1l}[][]{$p_1'$}
 \psfrag{p2}[][]{$p_{2}$}
 \psfrag{p2l}[][]{$p_{2}'$}
\psfrag{p3}[][]{$p_{3}$}
 \psfrag{p3l}[][]{$p_{3}'$}
\includegraphics[width=0.70\textwidth]{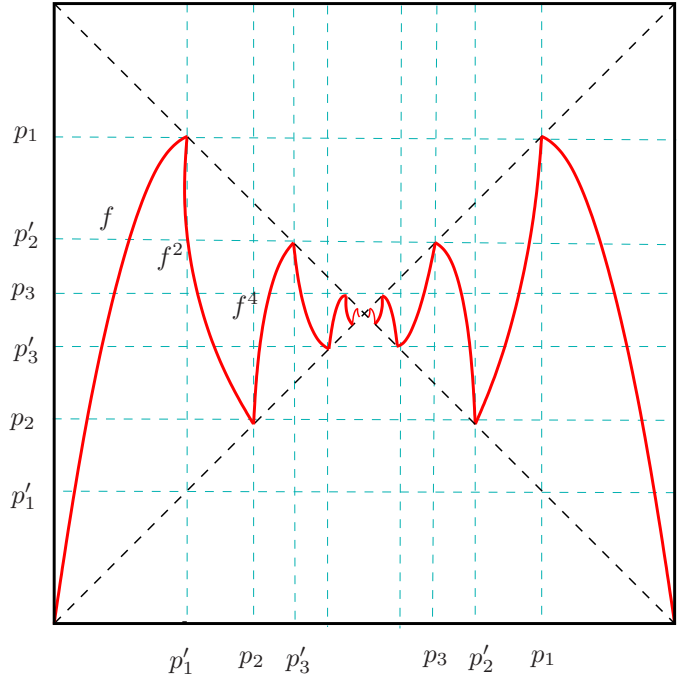}
\caption{The "Bat" map: the induced map $F$ for a  Feigenbaum unimodal map}
\end{figure}

\begin{proof}[{\bf Proof of Theorem \ref{apl1}.}] Let $f$ be an infinitely renormalizable map with bounded
combinatorics. We are going to define an induced map $F\colon
I\rightarrow I$, following  Y. Jiang (see \cite{jianga},
\cite{jiangb}):  Let $p_k$ be the periodic point in $\partial
I_k$. Define $E$ as the set
$$\{1,-1,-p_k, p_k,f(p_k), -f(p_k),\dots, f^{n_k-1}(p_k),-f^{n_k-1}(p_k)\}-\{f(p_k),-f(p_k)  \}.$$
The set $E$ cuts $I_{k-1}\setminus I_k$ in $m_k$ intervals.
Denote these intervals $M_{k-1,i}$, with $i=1,\dots,m_k$. For each $x \in M_{k-1,i}$, define $n(x)\geq 1$ as the minimal positive integer so that $$I_k\subset f^{n(x)n_{k-1}} M_{k-1,1}.$$ Note that $f^{n(x)n_{k-1}}$ does not have critical points on $ M_{k-1,i}$.  Define
the induced map $F$, which is defined everywhere in $I$, except
for a countable set of points:

$$F(x):=f^{n(x)}(x), \ for \ x \in I_k\setminus I_{k+1}.$$

See in Fig. 2 the induced map for an infinitely renormalizable
maps satisfying $n_{i+1}=2n_i$ for all $i$ (the so called
Feigenbaum maps). The map $F$ is Markovian with respect to the
partition
$$\mathcal{P}:=\{ M_{k,i}  \}_{k \in \mathbb{N},i \leq m_k}.$$
Furthermore, if $f$ and $g$ have the same bounded combinatorics
and even order, then by Proposition \ref{convren}, the
corresponding induced maps $F$ and $G$ satisfies

$$||\ \frac{1}{|I_k^f|}\ F(|M_{k,i}^f| \cdot +  |I_k^f|-|M_{k,i}^f|)-
\frac{1}{|I_k^g|}\ G(|M_{k,i}^g| \cdot + |I_k^g|-|M_{k,i}^g|)\
||_{C^r([0,1])} \leq C\lambda^k.$$

Define $L_k$ as, say, the right component of $I_k\setminus
I_{k+1}$ and $\gamma_k \colon  I \rightarrow L_k$ as the unique
bijective order preserving affine map between this two intervals.
We are going to define a random walk $\mathcal{F} \colon  I \times \mathbb{N} \rightarrow I \times \mathbb{N}$ from the map $F$ in the following
way:

\begin{equation}\label{bat_map}
\mathcal{F}(x,k):=
\begin{cases}
(\gamma^{-1}_i \circ F\circ \gamma_k (x),i) & \text{if $F\circ \gamma_k (x) \in L_i$;}\\
(\gamma^{-1}_i \circ (-F)\circ \gamma_k (x),i) & \text{if $F\circ \gamma_k (x) \in -L_i$.}\\
\end{cases}
\end{equation}
It is easy to see that we can extend  $\mathcal{F}\colon  I \times \mathbb{Z} \rightarrow I \times \mathbb{Z}$ to a strongly transient deterministic random walk with non-negative drift. Indeed if $k < 0$ define
$$\mathcal{F}(x,k):=
\begin{cases}
(\gamma^{-1}_i \circ F\circ \gamma_0 (x),k+i) & \text{if $F\circ \gamma_k (x) \in L_i$;}\\
(\gamma^{-1}_i \circ (-F)\circ \gamma_0 (x),k+i) & \text{if $F\circ \gamma_k (x) \in -L_i$.}\\
\end{cases}
$$
Furthermore if $g$ is another infinitely
renormalizable map with the  combinatorics of $f$
 then by Proposition \ref{convren} and Proposition \ref{how} we can define the corresponding random walk $\mathcal{G}\colon  I \times \mathbb{N} \rightarrow I \times \mathbb{N}$ and extend this to a random walk  $\mathcal{G}\colon  I \times \mathbb{Z} \rightarrow I \times \mathbb{Z}$  defining $\mathcal{G}(x,k)=\mathcal{F}(x,k)$ if $k < 0$. Then $\mathcal{G}$ is an asymptotically small
 perturbation of $\mathcal{F}$.  So we can apply  Theorem \ref{sqr} to conclude
that there is a conjugacy between $F$ and $G$ which is strongly
quasisymmetric with respect to the nested sequence of partitions
defined by the random walk $\mathcal{F}$. We can now easily
translate this result in terms of the original unimodal maps $f$
and $g$ saying that the continuous conjugacy $h$ between $f$ and
$g$ is a strongly quasisymmetric mapping with respect to
$\mathcal{P}$.
\end{proof}

\begin{rem}\label{quotient}{\rm An interesting case is when the unimodal map $f$ is a periodic point to the renormalization operator: there exists $n_0$ and $\lambda$, with $|\lambda|< 1$ so that
$$\frac{1}{\lambda}f^{n_0}(\lambda x)=f(x).$$
In this case, if we take $n_k=kn_0$, then the induced map $F$ will satisfy the functional equation
\begin{equation}\label{fe}  F(\lambda x)=\lambda F(x).\end{equation}
Define the relation $\sim$ in the following way: $$x\sim y \text{ iff there exists $i \in \mathbb{Z}$ so that } x=\pm \lambda^i y.$$
By Eq. (\ref{fe}), $F$ preserves this relation, so we can take the quotient of $F$ by the relation $\sim$.  Note that
$$L_0 = \mathbb{R}^\star/\sim.$$
It is easy to see that  $q=F/\sim\colon L_0\rightarrow L_0$ is a Markov expanding map. Now define $\psi\colon L_0 \rightarrow \mathbb{Z}$ as $\psi(x)= k$, if $f(x) \in I_{k}\setminus I_{k+1}$. Then $\mathcal{F}$ is exactly the homogeneous random walk defined by the pair $(q,\psi)$. }
\end{rem}

\subsection{Fibonacci maps} \label{aplf}
The Fibonacci renormalization is the simplest way to generalize
the concept of classical renormalization as described in Section
\ref{apl}. Actually we could prove all the results stated for
Fibonacci maps to a wider class of maps: maps which are infinitely
renormalizable in the generalized sense and with periodic
combinatorics and bounded geometry, but we will keep ourselves in the simplest case to
avoid more technical definitions and auxiliary results with its
long proofs.

Consider the class of real analytic maps $f$ with $Sf < 0$ and
defined in a disjoint union of intervals $I^0_1 \sqcup I^1_1$,
where $-I_1^0=I_1^0$, so that
\begin{itemize}
\item[ ]\
\item[{\it -}] The map $f\colon I^1_1 \rightarrow I^0_0:=f(I^1_1)$ is a
diffeomorphism. Furthermore $I^1_1$ is compactly contained in
$I^0_0$.\\

\item[{\it -}] The map $f\colon I^0_1 \rightarrow I^0_0$ is an even map
which has as $0$ as its unique critical point of even order.\\
\end{itemize}

We say that $f$ is {\bf Fibonacci renormalizable} if
$$f(0) \in I_1^1, \ f^2(0) \in I^1_0 \ and \  f^3(0) \in I^1_0.$$
In this case, the Fibonacci renormalization of $f$ is defined as
the first return map to the interval $I_1^0$ restricted to the
connected components of its domain which contain the points $f(0)$
and $f^2(0)$. This new map is denoted $\mathcal{R}f$: it could be
Fibonacci renormalizable again and so on, obtaining an infinite
sequence of renormalizations $\mathcal{R}f$,  $\mathcal{R}^2f$,
$\mathcal{R}^3f$, $\dots$.

We will denote the set of infinitely renormalizable maps in
the Fibonacci sense with a critical point of order $d$ by
$\mathcal{F}_d$. A map $f \in \mathcal{F}_d$ will be called a {\bf Fibonacci map}.

As in the original map $f$, the $n$-th renormalization $f_n:=
\mathcal{R}^nf$ of $f$ is a map defined in two disjoint intervals,
denoted $I^0_n$ and $I_n^1$, where  $-I^n_0=I^n_0$. Indeed $f_n$
on $I_n^0$ is a unimodal restriction of the $S_n$-th iteration of
$f$, where $\{ S_n \}$ is the Fibonacci sequence
$$S_0 = 1, \ S_1=2, \ S_2=3, \ S_3= 5, \ \dots \ , S_{k+2} = S_{k+1}
+ S_k, \dots$$ and $f_n$ on $I_n^1$ is the restriction of the
$S_{n-1}$-th iteration of $f$.

\begin{figure}
\centering \psfrag{f}{$f$}
\psfrag{un}[][][0.8]{$u_n$}
 \psfrag{unl}[][][0.8]{$u_n'$}
 \psfrag{un1}[][][0.8]{$u_{n+1}$}
 \psfrag{un1l}[][][0.8]{$u_{n+1}'$}
  \psfrag{unm1}[][][0.8]{$u_{n-1}$}
 \psfrag{unm1l}[][][0.8]{$u_{n-1}'$}

\psfrag{unm2}[][][0.8]{$u_{n-2}$}

 \psfrag{pn}[][][0.8]{$p_n$}
 \psfrag{pnl}[][][0.8]{$p_n'$}
 \psfrag{pn1}[][][0.8]{$p_{n+1}$}
 \psfrag{pn1l}[][][0.8]{$p_{n+1}'$}

\includegraphics[width=\textwidth]{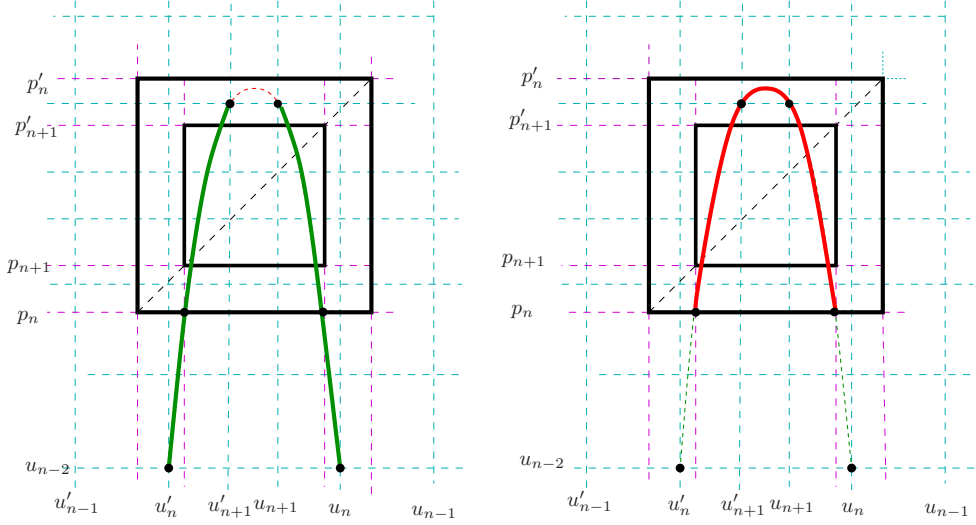}
\caption{On the left figure the solid curves represents the part of the $f^{S_n}$ used in the definition of the induced map. On the right figure the solid curve is the part of $f^{S_n}$ which coincides with the  $n$-th Fibonacci renormalization on its central domain. }
\end{figure}

\begin{figure}
\centering
\psfrag{f}{$f$}
\psfrag{un}[][][0.8]{$u_n$}
 \psfrag{unl}[][][0.8]{$u_n'$}
 \psfrag{un1}[][][0.8]{$u_{n+1}$}
 \psfrag{un1l}[][][0.8]{$u_{n+1}'$}
  \psfrag{unm1}[][][0.8]{$u_{n-1}$}
 \psfrag{unm1l}[][][0.8]{$u_{n-1}'$}

\psfrag{unm2}[][][0.8]{$u_{n-2}$}

\psfrag{pn}[][][0.8]{$p_n$}
\psfrag{pnl}[][][0.8]{$p_n'$}
\psfrag{pn1}[][][0.8]{$p_{n+1}$}
\psfrag{pn1l}[][][0.8]{$p_{n+1}'$}
\psfrag{pnm1}[][][0.8]{$p_{n-1}$}

\includegraphics[width=\textwidth]{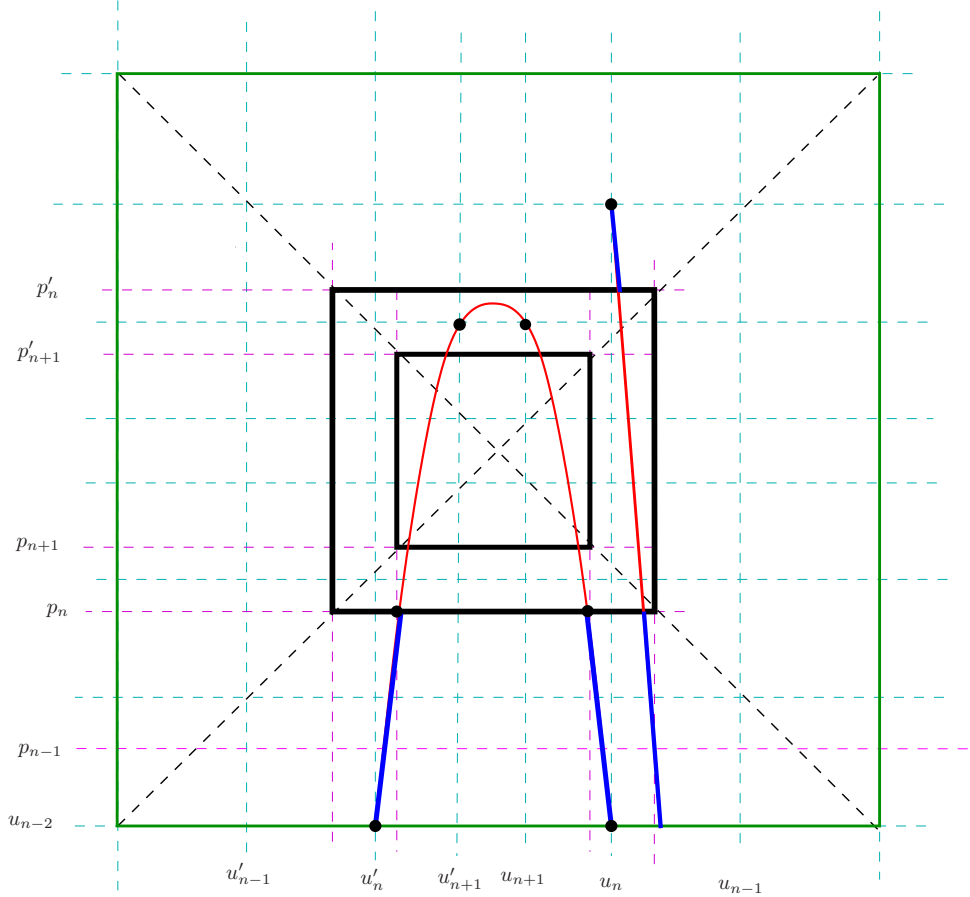}
\caption{The solid curves inside the medium square is the graph of the $n$-th Fibonacci renormalization $f_n$. The solid curves inside the largest square is the graph of an extension of $f_n$ which has the same maximal invariant set. }
\label{figure:extension}
\end{figure}

Denote by $p_k$ the sequence of points $p_k \in \partial I^k_0$ so that

$$f_k(p_{k+1})=p_k$$
and denote $I^k_0=[p_k, p_k']$.

It is possible to define a sequence $u_k$ of points satisfying
\begin{itemize}
\item[ ] \
\item[{\it 1.}]  $\dots <  \ p_{k+1} <  u_k <   p_k < \dots <  p_0, $ \\
\item[{\it 2.}] $f^{S_k}$ is monotone on $[0,u_k]$, \\
\item[{\it 3.}] $f^{S_k}(u_{k+1})=u_k$, \\
\item[{\it 4.}] $f^{S_k}(u_k)=u_{k-2}$. \\
\end{itemize}

We are going to define an induced map for an infinitely
renormalizable map in the Fibonacci sense in the following way:
Firstly, define $f_{-1}\colon I^0_0 \setminus I^1_0$ as an
$C^3$ monotone extension of $f_0$ on $I^1_1$ which has negative
Schwarzian derivative and bounded distortion. Define $F\colon I^0_0
\rightarrow \mathbb{R}$ as
$$F(x) := f^{S_i}(x) \ \ if \ \ x \in [u_i, -u_i]\setminus [u_{i+1},-u_{i+1}]$$
for each $i\geq 0$.

Define $L_i$ as, say, the right component of $[u_i, -u_i]\setminus [u_{i+1},-u_{i+1}]$ and $\gamma_i \colon  I \rightarrow L_i$ as the unique
bijective order preserving affine map between these two intervals.

We are ready to define the map $\mathcal{F}\colon I \times (\mathbb{N}\setminus \{0\}) \rightarrow I \times \mathbb{N}$ as
$$
\mathcal{F}(x,k):=
\begin{cases}
(\gamma^{-1}_i \circ F\circ \gamma_k (x),i) & \text{if $F\circ \gamma_k (x) \in L_i$,}\\
(\gamma^{-1}_i \circ (-F)\circ \gamma_k (x),i) & \text{if $F\circ \gamma_k (x) \in -L_i$.}\\
\end{cases}
$$

If the order of the critical point is even and larger than two then  there is a very special Fibonacci map $f^\star$, called the Fibonacci fixed point (see, for instance \cite{smfib}), whose induced map $F^\star$ satisfies (choosing a good $u_0$)
\begin{equation}\label{func_eq} F^\star(\lambda x) = \pm  \lambda F^\star(x)\end{equation}
for some $\lambda \in (0,1)$. In this case we can use the argument in Remark \ref{quotient} to conclude that the corresponding map $\mathcal{F}^\star\colon I \times (\mathbb{N}\setminus \{0\}) \rightarrow I \times \mathbb{N}$ can be extended to a homogeneous random walk $\mathcal{F}^\star\colon I \times \mathbb{Z}  \rightarrow I \times \mathbb{Z}$. For an arbitrary Fibonacci map $f$, we can extend $\mathcal{F}\colon I \times (\mathbb{N}\setminus \{0\}) \rightarrow I \times \mathbb{N}$ to a random walk  $\mathcal{F}\colon I \times \mathbb{Z}  \rightarrow I \times \mathbb{Z}$ defining $\mathcal{F}(x,k)=\mathcal{F}^\star(x,k)$ for $k\leq 0$. Then $\mathcal{F}$  is  not homogeneous,  however due Proposition \ref{how} and the following result  $\mathcal{F}$ is an asymptotically small perturbation of $\mathcal{F}^\star$:

\begin{prop}[see \cite{smfib}] For each even integer larger than two the following holds: for every Fibonacci map $f$, denote
$$g_i =    \alpha_{i}^{-1} \circ f^{S_i}\circ \alpha_{i+1}\colon I \rightarrow I,$$
where $\alpha_i\colon I \rightarrow [u_i^{f},-u_i^{f}]$ is an bijective affine map so that $\alpha_i^{-1}(f_{i+1}(0)) > 0$ and consider the correspondent maps $g_i^\star$ for $f^\star$. Then
$$|| g_i - g_i^{\star}||_{C^r} \leq K_r\rho^{i}$$
for some $\rho < 1$ and every $r \in \mathbb{N}$.
\end{prop}

The {\bf real Julia set} of $f$, denoted $J_{\mathbb{R}}(f)$, is
the maximal invariant of the map $$f\colon I_1^0 \sqcup I_1^1
\rightarrow I_0^0,$$ in other words,
$$J_{\mathbb{R}}(f_j):= \cap_i f^{-i}_j I_j^0. $$

Denote
$$\Omega_{+}^j(F):= \{  (x,i) \ s.t. \
\pi_2(F^n(x,i)) \geq j\, \ for \ all \ n\geq 0\}.$$

\begin{prop}\label{ji} There exists some $k_0$ so that
$$ \Omega_{+}^{j+1}(F)  \subset  J_{\mathbb{R}}(f_j) \subset  \Omega_{+}^{j-1}(F).$$  In particular
\begin{equation}\label{esthaus}
HD \ \Omega_{+}^{j+1}(F)  \leq HD \ J_{\mathbb{R}}(f_j) \leq   HD \     \Omega_{+}^{j-1}(F),
\end{equation}
and, for the Fibonacci fixed point, since $\Omega_{+}^{j+1}(F)$ is an affine copy of $\Omega_{+}^{j-1}(F)$ we have
\begin{equation}
HD \ \Omega_{+}^{j}(F) = HD \ J_{\mathbb{R}}(f).
\end{equation}
for all $j\geq 0$.
\end{prop}
\begin{proof} Denote by $F_\ell$ the restriction of $F$ to  $\cup_{i\geq \ell}L_i$. Then the maximal invariant set of $F_\ell$
$$\Lambda(F_\ell):= \cap_{i \in \mathbb{N}} F^{-i}\mathbb{R}$$
 is  $\Omega_{+}^{\ell}(F)$. Consider the extension of $f_j$ described in Fig. (\ref{figure:extension}). Let's call this extension $\tilde{f}_j$. An easy analysis of its graph shows that $f_j$ and $\tilde{f}_j$ have the same maximal invariant set. We claim that $\tilde{f}_{j+1}$ is just a map induced by $\tilde{f}_{j}$. Indeed, the restriction of $\tilde{f}_{j+1}$ to $[u_{j+1},u_{j+1}']$ coincides with $\tilde{f}_j^2$ on the same interval. On the rest of $\tilde{f}_{j+1}$-domain $\tilde{f}_{j+1}$ coincides with $\tilde{f}_{j}$.

\begin{figure}
\centering \psfrag{f}{$f$}

\psfrag{f2}[][][0.6]{$f_n$}
 \psfrag{f1}[][][0.7]{$f_{n+1}$}
 \psfrag{f0}[][][0.7]{$\ \ f_{n+2}$}
 \psfrag{f1p}[][][0.7]{$f_{n+3}$}
 \psfrag{f2p}[][][0.7]{$f_{n+4}$}

\psfrag{2}[][][0.6]{$u_n$}
 \psfrag{1}[][][0.6]{$u_{n+1}$}
 \psfrag{0}[][][0.6]{$u_{n+2}$}
 \psfrag{1p}[][][0.6]{$u_{n+3}$}
 \psfrag{2p}[][][0.6]{$u_{n+4}$}

 \psfrag{2t}[][][0.6]{$u_n'$}
 \psfrag{1t}[][][0.6]{$u_{n+1}'$}
 \psfrag{0t}[][][0.6]{$u_{n+2}'$}
 \psfrag{1pt}[][][0.6]{$u_{n+3}'$}
 \psfrag{2pt}[][][0.6]{$u_{n+4}'$}
\includegraphics[width=0.730\textwidth]{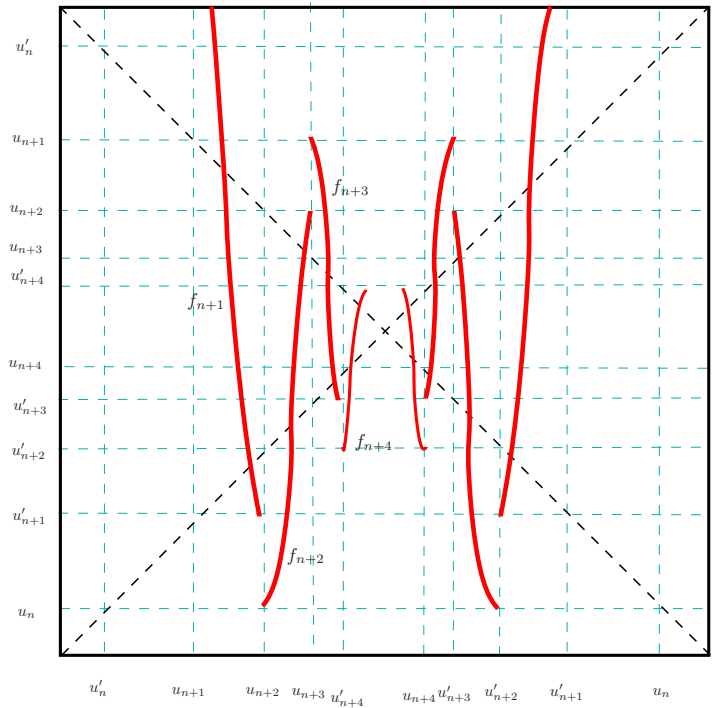}
\caption{Induced map $F$ for a  Fibonacci map}
\end{figure}

By consequence, for $i \geq j$ the map  $\tilde{f}_i$ is induced by $\tilde{f}_j$ and, since $F_{j+1}$ restricted to $L_i$ is equal to $\tilde{f}_i$, we obtain that $F_{j+1}$ is a map induced by $\tilde{f}_{j}$. In particular
$$\Lambda(F_{j+1}) \subset \Lambda(\tilde{f}_j)=J_{\mathbb{R}}(f_j).$$

To prove that $\Lambda(\tilde{f}_j) \subset \Lambda(F_{j-1})$, we are going to prove that
\begin{equation}\label{contido} x \in \Lambda(\tilde{f}_j) \ implies \ F_{j-1}(x) \in \Lambda(\tilde{f}_j).\end{equation}

If $x$ belongs to the interval $I^1_j \subset L_{j-1}$, where $\tilde{f}_j$ coincides with $F_{j-1}$, then $F_{j-1}(x) \in \Lambda(\tilde{f}_j)$. Otherwise  $x \in I^0_j \subset \cup_{i\geq j}L_i$, so $x \in \Lambda(\tilde{f}_j)\cap \  L_i$, for some $i\geq j$, then $F_{j-1}$ is an iteration of  $\tilde{f}_j$ on $L_i$, so $F_{j-1}(x) \in \Lambda(\tilde{f}_j)$. This finishes the proof of Eq. (\ref{contido}). Since $\Lambda(\tilde{f}_j)$ is invariant by the action of $F_{j-1}$ we have $\Lambda(\tilde{f}_j) \subset \Lambda(F_{j-1})$.\end{proof}

\begin{proof}[{\bf Proof of Theorem \ref{juliathm}}] Consider the homogeneous random walk $F^\star= (g,\psi)$ induced by $f^\star$. Denote
$$M = \int \psi \ d\mu,$$
where $\mu$ is the absolutely continuous invariant measure of $g$. Using Thorem \ref{menor}, there are three cases:
\vspace{4mm}

{\bf 1. ${ \mathbf M < 0}$.} In this case $\mathbf F^\star$ is transient and we have that  $HD \ \Omega_+(F) < 1$ for every asymptotically small perturbation of $F^\star$, in particular when $F$ is a random walk induced by a Fibonacci map $f$. By Proposition \ref{ji}, $HD \ J_{\mathbb{R}}(f) < 1$.
\vspace{4mm}

{\bf 2. $\mathbf M =0$.} Then   $F^\star$ is recurrent \cite{guivarch} so  every asymptotically small perturbation $G$ of $F^\star$ is recurrent and    $m(\Omega_+(G))=0$ but  $HD \ \Omega_+(G)=1$. By Proposition \ref{ji} we obtain $m(J_\mathbb{R}(f)) =0$ and $HD \ J_\mathbb{R}(f)=1$.
\vspace{4mm}

{\bf 3. $\mathbf M > 0$.} In this case $\mathbf F^\star$ is transient with $m(\Omega_+(F^\star)) > 0$ and the conjugacy between $F^\star$ and any asymptotically small perturbation of it is absolutely continuous on $\Omega_+^i(F^\star)$. In particular $m(\Omega_+(F)) > 0$ for every random walk $F$ induced by a Fibonacci map $f$ so $m(J_\mathbb{R}(f)) > 0$ by Proposition \ref{ji}.\end{proof}

A  map $f\colon I \rightarrow I$ is called a unimodal map if $f$ has a unique critical point, with even order $d$, which is a maximum, and $f(\partial I) \subset \partial I$. We will assume that $f$ is real analytic, symmetric with respect the critical point and  $Sf < 0$. If the critical value is high enough, then $f$ has a reversing fixed point $p$. Let $I_0^0:=[-p,p]$. Consider the map of first return $R$ to $f$: if $x \in I_0^0$ and $f^r(x) \in I_0^0$, but $f^n(x) \not\in I_0^0$ for $i < r$, define $$R(x):=f^r(x).$$
If there exists exactly two connected components $I_1^0$ and $I_1^1$ of the domain of $R$ containing points in the orbit of the critical point, and furthermore the map $$R\colon I_1^0 \cup I_1^1 \rightarrow I_0^0$$ is a Fibonacci map,  then we will called $f$ an {\bf unimodal Fibonacci map}. The class of all unimodal Fibonacci maps will be denoted $\mathcal{F}^{uni}_d$.

\begin{proof}[{\bf Proof of Theorem \ref{deep}}]We will use the notation in the proof of Theorem \ref{juliathm}. Since $m(J_{\mathbb{R}}(f)) > 0$, we conclude that the mean drift $M$ of $F^\star$ is positive. By Proposition \ref{prws}
any asymptotically small perturbation $G$ of $F^\star$ has the following property: there exists $\lambda \in [0,1)$, $C >0$  and $K > 0$  so that for every $P \in \mathcal{P}^0(G)$
$$m(p \in P \colon \ \sum_{i=0}^{n-1}\psi(G^i(p)) < K n  )\leq C\lambda^n |P|.$$
This implies that
$$m(p \in I_j \colon \ \sum_{i=0}^{\ell}\psi(G^i(p)) \geq  K \ell \text{ for every } \ell\geq n  )\geq (1-C\lambda^n).$$
So if $j= n|min \psi|$ we obtain
$$m(\Omega_+^j(G))\geq 1-C\lambda^{C_1 j}.$$
here $C_1 > 0$. If $G$ is a random walk induced by a Fibonacci map $g$ then this implies that for $j$ large
$$m(L_j\setminus J_{\mathbb{R}}(g))= m((-L_j)\setminus J_{\mathbb{R}}(g)) \leq C\lambda^{C_1 j}|L_j|.$$
Since $$[-u_{j+1},u_{j+1}]=\bigcup_{i\geq j} L_i\cup(-L_i),$$
we conclude that
\begin{equation}\label{soalguns} m([u_{j+1},-u_{j+1}]\setminus J_{\mathbb{R}}(g)) \leq C\lambda^{C_1 j}|u_{j+1}|.\end{equation}
For every $\delta$, choose $j$ so that $|u_{j+2}|\leq \delta \leq |u_{j+1}|$. Because $|u_{j+2}|> \theta |u_{j+1}| $, where $\theta \in (0,1)$ does not depend on $j$, we have that $|u_{j}|\geq C\theta^j$. Together with Eq. (\ref{soalguns}) this implies
$$m([-\delta,\delta]\setminus J_{\mathbb{R}}(g)) \leq C\lambda^{C_1 j}|u_{j+1}|\leq C|u_{j+1}|^{1+ \alpha}\leq C|\delta|^{1+ \alpha}.$$\end{proof}

\begin{proof}[\bf Proof of Theorem \ref{apl2}]We will prove each one of the following implications:

{\bf (1) implies  (2):} From the proof of Theorem \ref{juliathm}, if $m(J_{\mathbb{R}}(f))> 0$ for some $f \in \mathcal{F}_d$ the mean drift $M$ of the homogeneous random walk $\mathcal{F}^\star$ of $f^\star$ is positive. So $\mathcal{F}^\star$ (and all its asymptotically small perturbations) is transient (to $+\infty$). In terms of the original Fibonacci map $f$, this means that almost every orbit  in $J_{\mathbb{R}}(f)$ accumulates in the post-critical set: So $f$ has a wild attractor.

{\bf (2) implies (3):} if there exists a wild attractor for $f$ then $m(J_{\mathbb{R}}(f))> 0$. From the proof of Theorem \ref{juliathm} we obtain that the mean drift $M$ of $\mathcal{F}^\star$  is positive. So there exists a absolutely continuous conjugacy between $\mathcal{F}^\star$ and any asymptotically small perturbation of $\mathcal{F}^\star$. This implies that any two maps $f_1, f_2 \in \mathcal{F}_d$ admits a continuous and absolutely continuous conjugacy  $$h\colon J_{\mathbb{R}}(f_1) \rightarrow J_{\mathbb{R}}(f_2).$$
Now consider two arbitrary maps $g_1, g_2 \in \mathcal{F}_d^{uni}$. Then we already know that there exists an absolutely continuous conjugacy
$$h\colon J_{\mathbb{R}}(R_{g_1}) \rightarrow J_{\mathbb{R}}(R_{g_2})$$ between the induced Fibonacci maps $R_{g_1}$ and $R_{g_2}$ associated to $g_1$ and $g_2$. Of course $h$ is just the restriction of a topological conjugacy between $g_1$ and $g_2$. By a Blokh and Lyubich result \cite{bl} (see also page 332 in \cite{ms}), every map of $\mathcal{F}_d^{uni}$ is ergodic with respect the Lebesgue measure. Since $g_1$ and $g_2$ have wild attractors, this implies that the orbit of almost every point $x \in I$ hits  $J_{\mathbb{R}}(R_{g_1})$ at least once. Let $n(x)$ be a time when this happens.

So consider a  arbitrary measurable set $B \subset I$ so that $m(B)>0$. Then for at least one $n_0 \in \mathbb{N}$  the set
$$B_{n_0}:=\{x \in B\colon \ n(x)=n_0  \}$$
has positive Lebesgue measure. This implies that $f^{n_0}B_{n_0}$ has positive Lebesgue measure, so $m(h(f^{n_0}B_{n_0})) > 0$. Now it is easy to conclude that $m(h(B_{n_0}))> 0$ and $h(B) > 0$. Switching the places of $g_1$ and $g_2$ in this argument we can conclude that $h$ is absolutely continuous on $I$.

Finally note that the eigenvalues of the periodic points are not constant on the class $\mathcal{F}_d^{uni}$.

{\bf (3) implies (4):} By the argument in Martens and de Melo \cite{mm}, if a Fibonacci map does not have a wild attractor then any continuous absolutely continuous conjugacy with other Fibonacci map is $C^1$: in particular the conjugacy preserves the eigenvalues of the periodic points. So if (3) holds then we can use the same argument in the proof of the previous implication to conclude that every Fibonacci map has a wild attractor.

{\bf (4) implies (5):} The proof goes exactly as  the proof of  (2)$\Rightarrow$ (3).

{\bf (5) implies (1):} The proof goes exactly as  the proof of  (3)$\Rightarrow$ (4).

\end{proof}

\end{document}